\documentclass[a4paper,dvipsnames,fleqn]{article}

\usepackage{amsmath}
\usepackage{amssymb}
\usepackage{amsthm}
\usepackage{mathtools}
\usepackage{arydshln}
\usepackage{bm}
\usepackage[colorlinks=true,allcolors=purple]{hyperref}
\usepackage{cleveref}
	\crefformat{equation}{#2(#1)#3}
\usepackage{dsfont}
\usepackage{enumitem}
\usepackage[OT2,OT1]{fontenc}
\usepackage{ifthen}
\usepackage[cal=cm,scr=boondoxo]{mathalfa}
\usepackage{pifont}
\usepackage{stmaryrd}
\usepackage{textcomp}
\usepackage{tikz}
	\usetikzlibrary{arrows}
	\usetikzlibrary{patterns}
\usepackage{tikz-cd}
\usepackage[textwidth=3cm,colorinlistoftodos]{todonotes}
\usepackage{xfrac}

\numberwithin{equation}{section}			
\swapnumbers						

\newcommand\cyr{
\renewcommand\rmdefault{wncyr}%
\renewcommand\sfdefault{wncyss}%
\renewcommand\encodingdefault{OT2}%
\normalfont
\selectfont}
\DeclareTextFontCommand{\textcyr}{\cyr}

\makeatletter						
\newcommand\captionof[1]{\def\@captype{#1}\caption}	
\makeatother

		\allowdisplaybreaks
\newcounter{Enum}				
\newenvironment{Enumerate}{\begin{enumerate}[label={\rm({\roman*})}]}{\end{enumerate}}

\newcommand{\descriptionlabelsave}{}		
\newenvironment{Itemize}{%
	\renewcommand{\descriptionlabelsave}{\descriptionlabel}\renewcommand{\descriptionlabel}{$\triangleright$}%
	\begin{description}[leftmargin=15pt,itemindent=-5.2pt]}{%
	\end{description}\renewcommand{\descriptionlabel}{\descriptionlabelsave}}

\newcounter{StepsCount}				
\newenvironment{Elist}{%
	\begin{list}{\ding{\value{StepsCount}}}{\usecounter{StepsCount} \leftmargin=0pt \labelwidth=12pt \itemindent=\labelwidth%
	\itemsep=5pt\listparindent=\parindent} \setcounter{StepsCount}{191}}{\end{list}}
\newcounter{StepsRefCount}

\newenvironment{Ilist}{
	\begin{list}{$\triangleright$}{\leftmargin=0pt \labelwidth=11pt \itemindent=\labelwidth%
	\itemsep=5pt\listparindent=\parindent}}{\end{list}}

\newenvironment{IIlist}{
	\begin{list}{$\triangleright\triangleright$}{\leftmargin=0pt \labelwidth=17pt \itemindent=\labelwidth%
	\itemsep=5pt\listparindent=\parindent}}{\end{list}}

\theoremstyle{plain}
	\newtheorem{lemma}{Lemma}[section]
	\newtheorem{proposition}[lemma]{Proposition}
	\newtheorem{theorem}[lemma]{Theorem}
	\newtheorem{corollary}[lemma]{Corollary}
	\newcommand{\GenericTheoremName}{}\newtheorem{generictheorem}[lemma]{\GenericTheoremName}
\theoremstyle{definition}
	\newtheorem{definition}[lemma]{Definition}
	\newcommand{\GenericDefinitionName}{}\newtheorem{genericdefinition}[lemma]{\GenericDefinitionName}
\theoremstyle{remark}
	\newtheorem{remark}[lemma]{Remark}
	\newtheorem{example}[lemma]{Example}
	\newcommand{\GenericRemarkName}{}\newtheorem{genericremark}[lemma]{\GenericRemarkName}
\newenvironment{Lemma}{\begin{lemma}}{\par\noindent\rule{5em}{1pt}\end{lemma}}
\newenvironment{Corollary}{\begin{corollary}}{\par\noindent\rule{5em}{1pt}\end{corollary}}
\newenvironment{Definition}{\begin{definition}}{\par\noindent\rule{5em}{1pt}\end{definition}}

\newenvironment{GenericDefinition}[1]
	{\renewcommand{\GenericDefinitionName}{#1}\begin{genericdefinition}}{\par\noindent\rule{5em}{1pt}\end{genericdefinition}}
\newenvironment{Remark}{\begin{remark}}{\par\noindent\rule{5em}{0.5pt}\end{remark}}

\newcommand{\mc}[1]{{\mathcal{#1}}}			
\newcommand{\ms}[1]{{\mathscr{#1}}}			
\newcommand{\bb}[1]{{\mathbb{#1}}}			

\newcommand{\Side}[1]{\hfill{#1}\kern10pt}		

\newcommand{\smmatrix}[4]{\Bigl(			
\begin{smallmatrix}
\hspace*{-0.2ex} #1 \hspace*{0.2ex} & \hspace*{0.2ex} #2 \hspace*{-0.2ex}
\\[0.5ex]
\hspace*{-0.2ex} #3 \hspace*{0.2ex} & \hspace*{0.2ex} #4 \hspace*{-0.2ex}
\end{smallmatrix}
\Bigr)}
\newcommand{\smfrac}[2]{{\textstyle\frac{{#1}}{{#2}}}}	

\newcommand{\Dummy}{\text{\textvisiblespace\kern1pt}}	
\newcommand{\BigO}{{\rm O}}				

\newcommand{\DS}{\mid\mkern3mu}				
\newcommand{\DP}{{\mathrel{\mathop:}\kern5pt}}		
\newcommand{\DF}{\colon}				
\newcommand{\DE}{\mathrel{\mathop:}=}			
\newcommand{\DI}{\mathrel{\mathop:}\Leftrightarrow}	
\newcommand{\DD}{\mkern4mu\mathrm{d}}			
\newcommand{\CAS}{&\text{if}\ }				
\newcommand{\CASO}{&\text{otherwise}\ }			

\DeclareMathOperator{\GL}{GL}				

\DeclareMathOperator{\Ind}{ind}				
\newcommand{\PL}[3]{{#3}^{#1}(\log{#3})^{#2}}		
\newcommand{\pL}[2]{(\log{#2})^{#1}}

\begin{document}

\begin{flushleft}
	{\Large\bf An upper bound for the Nevanlinna matrix\\[2mm] of an indeterminate moment sequence}
	\\[5mm]
	\textsc{
	Raphael Pruckner
	\,\ $\ast$\,\ 
	Jakob Reiffenstein
	\,\ $\ast$\,\ 
	Harald Woracek
		\hspace*{-14pt}
		\renewcommand{\thefootnote}{\fnsymbol{footnote}}
		\setcounter{footnote}{2}
		\footnote{This work was supported by the project P~30715-N35 of the Austrian Science Fund (FWF). 
			The second and third authors were supported by the joint project I~4600 of the Austrian
			Science Fund (FWF) and the Russian foundation of basic research (RFBR).}
		\renewcommand{\thefootnote}{\arabic{footnote}}
		\setcounter{footnote}{0}
	}
	\\[6mm]
	{\small
	\textbf{Abstract:}
		The solutions of an indeterminate Hamburger moment problem can be parameterised using the Nevanlinna matrix of 
		the problem. The entries of this matrix are entire functions of minimal exponential type, and any growth less
		than that can occur. 

		An indeterminate moment problem can be considered as a canonical system in limit circle
		case by rewriting the three-term recurrence of the problem to a first order vector-valued recurrence. We give a
		bound for the growth of the Nevanlinna matrix in terms of the parameters of this canonical system. In most
		situations this bound can be evaluated explicitly. It is sharp in the sense that for well-behaved
		parameters it coincides with the actual growth of the Nevanlinna matrix up to multiplicative	constants. 
	\\[3mm]
	}
\end{flushleft}


%
%
%
\section{Introduction}

Let $(s_n)_{n=0}^\infty$ be a sequence of real numbers. The \emph{Hamburger moment problem} is the task of describing the set 
\[
	\mc M\big((s_n)_{n=0}^\infty\big)\DE
	\bigg\{\,\mu\mkern7mu\Big|\mkern10mu
		\parbox[c]{49mm}{\small $\mu$ positive Borel measure on $\bb R$\\ 
		$s_n=\int_{\bb R}t^n\DD\mu(t)$ for $n=0,1,2,\ldots$}
	\bigg\}
	.
\]
This is a classical problem of analysis and was treated extensively in work of H. Hamburger, M. Riesz, R. Nevanlinna, and many
others. Standard references are, e.g., \cite{shohat.tamarkin:1943,akhiezer:1961,schmuedgen:2017}. 
It is known that one of the following alternatives takes place for the set $\mc M((s_n)_{n=0}^\infty)$: 
it is either empty, or contains exactly one element, or contains infinitely many elements. 
In the last case $(s_n)_{n=0}^\infty$ is called an \emph{indeterminate moment sequence}, and this is the case we are concerned 
with in the present paper. 

For an indeterminate moment sequence $(s_n)_{n=0}^\infty$ the set $\mc M\big((s_n)_{n=0}^\infty\big)$ can be parameterised:
There exist four entire functions $A,B,C,D$, such that the formula 
\begin{equation}
\label{B90}
	\int_{\bb R}\frac 1{t-z}\DD\mu(t)=\frac{A(z)\tau(z)+B(z)}{C(z)\tau(z)+D(z)}
\end{equation}
establishes a bijection between $\mc M((s_n)_{n=0}^\infty)$ and the set of all functions $\tau$ that are analytic in the open
upper half-plane $\bb C^+$ and have nonnegative imaginary part (formally including the function constant equal to $\infty$). The
matrix 
\[
	W(z)\DE\begin{pmatrix} A(z) & B(z) \\ C(z) & D(z)\end{pmatrix}
\]
is called the \emph{Nevanlinna matrix} of the sequence $(s_n)_{n=0}^\infty$. 

The Nevanlinna matrix admits an operator theoretic interpretation 
(and this viewpoint could be used to prove \cref{B90}). 
To the moment sequence there is an
associated sequence $(p_n)_{n=0}^\infty$ of orthogonal polynomials that satisfies a three-term recurrence of the form
\begin{align}
\label{B95}
	zp_n(z)=b_np_{n+1}(z)+a_np_n(z)+b_{n-1}p_{n-1}(z),\qquad n=0,1,2,\ldots
\end{align}
with certain parameters $a_n \in \bb R$ and $b_n>0$ for $n=0,1,2,\ldots$ (and formally $b_{-1}=0$). 
The \emph{Jacobi operator} is the closure of the operator 
\[
	(Ju)_n= 
	\begin{cases} 
		b_0u_1+a_0u_0 \CAS n=0,
		\\
		b_nu_{n+1}+a_nu_n+b_{n-1}u_{n-1}\CAS n \geq 1
	\end{cases}
\]
defined on $\mc D\DE \{u \in \ell^2\DS u_n=0 \text{ for almost all } n \}$, and \cref{B95} is the formal eigenvalue equation 
for $J$. 
The Jacobi operator is closed and symmetric, and has deficiency index $(1,1)$. Hence, the self-adjoint extensions of the Jacobi operator (we write again $J$ for simplicity) are 
described by M.G. Krein's resolvent formula. 
The Nevanlinna matrix of $(s_n)_{n=0}^\infty$ is precisely the $u$-resolvent matrix of $J$ 
for a certain generating element $u$. 

The entries $A,B,C,D$ of the Nevanlinna matrix all have the same growth \cite{berg.pedersen:1994}, and a classical theorem of
M. Riesz in \cite{riesz:1923a} states that they are of minimal exponential type. Moreover, it is known that any growth smaller
than that may occur, e.g.\ \cite{borichev.sodin:1998}. 
Revealing more refined information about the growth of the Nevanlinna matrix is an intricate problem. It is of great interest
also for spectral theoretic reasons: due to the above interpretation of $W$ as a resolvent matrix, the set of zeros of $D$
coincides with the spectrum of a particular self-adjoint extension of $J$. Hence, if the growth of $W$ is known, information 
about the distribution of eigenvalues can be obtained using standard tools from complex analysis. 

In the present paper we give an upper bound for 
\[
	M(r)\DE\max_{|z|=r}\|W(z)\|	
\]
of which there are two versions: First, a general formulation where any (possibly rough) data is admitted, and second, a much
more explicit bound for data satisfying mild regularity conditions. If, in addition, the data decays sufficiently fast, then 
the upper bound coincides with $M(r)$ up to multiplicative constants. 

So far, we talked about two different (equivalent) objects, i.e., the moment sequence and the Jacobi operator. Our method
actually relies on a third object -- a \emph{canonical system} of differential equations. It is obtained from rewriting the
three-term recurrence \eqref{B95} as a first order vector difference equation and interpreting this as a discrete differential
equation. From the fundamental solution of the canonical system, we obtain its \emph{monodromy matrix} 
which again coincides with the Nevanlinna matrix $W$. 

The canonical systems occurring in the context of moment problems are represented by a \emph{Hamiltonian}
which reflects the discrete nature of the difference equation. It is determined by two sequences, its \emph{lengths} $l_j>0$ 
and \emph{angles} $\phi_j\in\bb R$, which we call the \emph{Hamiltonian parameters} (details are given below, cf.\ \Cref{B22}). 
This model for an indeterminate moment sequence is well suited for the study of various properties. For example, the moment
sequence is indeterminate if and only if $(l_j)_{j=1}^\infty$ is summable. Our results will almost exclusively 
be formulated in terms of the Hamiltonian parameters.

Let us briefly review some earlier results on the growth of Nevanlinna matrices.
The history of the subject starts probably with a theorem of M.S. Liv\v sic in \cite{livshits:1939} that gives a
lower bound for $M(r)$ in terms of the moment sequence itself. This bound is easy to handle, but will usually not give the 
correct size. 
In recent work of C. Berg and R. Szwarc \cite{berg.szwarc:2014} it is shown that the order and type of $W$ 
\[
	\rho\DE\limsup_{r\to\infty}\frac{\log\log M(r)}{\log r},\qquad 
	\tau\DE\limsup_{r\to\infty}\frac{\log M(r)}{r^\rho}
	,
\]
coincide with those of a certain entire function built in a complicated way from the coefficients 
of the orthonormal polynomials.
A theorem which takes the Jacobi parameters as input is due to Yu.M. Berezanskii \cite{berezanskii:1956} and was generalised in \cite{berg.szwarc:2014}. It states that the order of $W$
coincides with the convergence exponent of the sequence $(b_n)_{n=0}^\infty$, but the assumptions are very restrictive, involving regularity of $b_n$ and smallness of $a_n$.
Bounds for the order of $W$ in terms of the Hamiltonian parameters, which under certain conditions give the correct
value, are obtained in \cite{pruckner.romanov.woracek:jaco}. 

In our present theorems we start with the Hamiltonian parameters, and give bounds for the type of $W$ with respect to a general
comparison function (e.g., a proximate order in the sense of Valiron). 
This improves significantly upon \cite{pruckner.romanov.woracek:jaco} in several ways: 
we work on a much more refined scale of measuring growth, we obtain type estimates, and in some situations our bound improves the 
earlier results even on the rough scale of order. Several cases occur, which are presented in \Cref{B7}. The path that leads to
this result is divided into two main sections: First, we prove a very general, albeit complicated, upper bound in \Cref{B8}. 
Second, in \Cref{B9} we use J. Karamata's theory of regularly varying functions to evaluate this general bound. 
For the convenience of the reader, the significantly simpler case of usual order and type is covered separately in \Cref{B66}.
We show in \Cref{B38} that the upper bound is attained if the lengths and angle differences are themselves close to regularly
varying and decay sufficiently fast (corresponding to order less than $\frac 12$). 
In this case the growth of $\log M(r)$ is fully determined up to multiplicative constants. 

The proof of our foundational theorem is based on a somewhat tricky application of a recent result from
\cite{pruckner.woracek:sinqA}, approximating the target Hamiltonian with a finite dimensional one. Building upon that, a
detailed (and partly tedious) analysis of functions follows. We recommend the reader to get an overall picture before diving
into the actual estimates. 

\subsection*{Canonical systems with Hamburger Hamiltonian}

A two-dimensional canonical system is an equation of the form
\[
	y'(t)=zJH(t)y(t), \qquad t \in (a,b)\text{ a.e.},
\]
where 
\begin{Ilist}
\item $-\infty\leq a<b\leq\infty$,
\item $H\in L^1_{\rm loc}\big((a,b),\bb R^{2\times 2}\big)$, and $H(t)\geq 0$ and $H(t)\neq 0$ for $t\in(a,b)$ a.e.,
\item $J\DE\smmatrix 0{-1}10$ and $z\in\bb C$, 
\end{Ilist}
and the solution $y\DF(a,b)\to\bb C^2$ is required to be locally absolutely continuous. 
The function $H$ is called the \emph{Hamiltonian} of the system. 

Canonical systems that occur from indeterminate moment sequences are those whose Hamiltonian has the
following particular -- discrete -- form. Here we denote 
\[
	\xi_\phi\DE\binom{\cos\phi}{\sin\phi},\qquad \phi\in\bb R
	.
\]

\begin{Definition}
\label{B22}
	Let $(l_j)_{j=1}^\infty$ be a summable sequence of positive numbers and $(\phi_j)_{j=1}^\infty$ be a sequence of real
	numbers. Set $x_0 \DE 0$, $x_j \DE \sum_{i=1}^j l_i$ for $j \geq 1$, and $L\DE\sum_{j=1}^\infty l_j <\infty$. 
	The \emph{Hamburger Hamiltonian} with \emph{lengths} $l_j$ and \emph{angles} $\phi_j$ is the function $H$ 
	defined on the interval $(0,L)$ as 
	\[
		H(t)\DE\xi_{\phi_j}\xi_{\phi_j}^T\quad\text{for }j\in\bb N\text{ and }
		x_{j-1}\leq t<x_j
		.
	\]
\end{Definition}

\noindent
A Hamburger Hamiltonian $H$ thus can be pictured as 
\begin{center}
\begin{tikzpicture}[x=1.2pt,y=1.2pt,scale=0.8,font=\fontsize{8}{8}]
	\draw[thick] (10,30)--(215,30);
	\draw[dotted, thick] (215,30)--(270,30);
	\draw[thick] (10,25)--(10,35);
	\draw[thick] (70,25)--(70,35);
	\draw[thick] (120,25)--(120,35);
	\draw[thick] (160,25)--(160,35);
	\draw[thick] (190,25)--(190,35);
	\draw[thick] (210,25)--(210,35);
	\draw[thick] (270,25)--(270,35);
	\draw (40,44) node {${\displaystyle \xi_{\phi_1}\xi_{\phi_1}^T}$};
	\draw (95,44) node {${\displaystyle \xi_{\phi_2}\xi_{\phi_2}^T}$};
	\draw (140,44) node {${\displaystyle \xi_{\phi_3}\xi_{\phi_3}^T}$};
	\draw (177,43) node {${\cdots}$};
	\draw (-20,30) node {\large $H\!:$};
	\draw (10,18) node {${\displaystyle x_0}$};
	\draw[dashed,stealth-stealth] (11,26)--(69,26);
	\draw (40,21) node {${l_1}$};
	\draw (70,18) node {${x_1}$};
	\draw[dashed,stealth-stealth] (71,26)--(119,26);
	\draw (95,21) node {${l_2}$};
	\draw (120,18) node {${x_2}$};
	\draw[dashed,stealth-stealth] (121,26)--(159,26);
	\draw (140,21) node {${l_3}$};
	\draw (160,18) node {${x_3}$};
	\draw (195,18) node {${\cdots}$};
	\draw (270,18) node {${\displaystyle L}$};
\end{tikzpicture}
\end{center}
Since a Hamburger Hamiltonian in the sense of the above definition is even integrable on the whole interval $(0,L)$, 
there exists a unique $2\times 2$ matrix-valued solution $W\DF[0,L]\times\bb C\to\bb C^{2\times 2}$ of the initial value problem 
\[
	\left\{
	\begin{array}{l}
		\frac{\partial}{\partial t} W(t;z)J=zW(t;z)H(t),\qquad t\in(0,L)\text{ a.e.},
		\\[2mm]
		W(0;z)=I.
	\end{array}
	\right.
\]
We refer to $W$ as the \emph{fundamental solution} of $H$, and to the matrix $W_H(z)\DE W(L;z)$ as its 
\emph{monodromy matrix}. 

\subsection*{A notational convention}

We frequently compare functions up to multiplicative constants or asymptotically, and throughout the paper use the following 
notation.

\begin{GenericDefinition}{Notation}
	Let $X$ be a set and $f,g\DF X\to(0,\infty)$.
	\begin{Enumerate}
	\item We write ``$f\lesssim g$\/'' (or ``$f(x)\lesssim g(x)$\/'') to say that 
		there exists a constant $C>0$ such that $f(x)\leq C\cdot g(x)$ for all $x\in X$. 
		We write ``$f\gtrsim g$\/'' if $g\lesssim f$, and ``$f\asymp g$\/'' if $f\lesssim g$ and $f\gtrsim g$. 
	\item Assume that $X$ is directed. We say that ``$f\lesssim g$ \emph{for sufficiently large} $x$\/'', if there exists 
		$x_0 \in X$ such that $f|_Y\lesssim g|_Y$ where $Y\DE\{x\in X\DS x\succeq x_0\}$.
	\item Assume again that $X$ is directed. We write ``$f\sim g$\/'' if $\lim_{x\in X}\frac{f(x)}{g(x)}=1$, 
		and ``$f\ll g$\/'' if $\lim_{x\in X}\frac{f(x)}{g(x)}=0$,
		and ``$f\approx g$\/'' if $\lim_{x\in X}\frac{f(x)}{g(x)}$ exists in $(0,\infty)$.
	\item Assume that $X$ is a subset of a topological space. We say that ``$f\lesssim g$ \emph{locally}\/'', if 
		$f|_K\lesssim g|_K$ for every compact subset $K$ of $X$. Analogous wording applies to ``$\asymp$\/''.
	\end{Enumerate}
\end{GenericDefinition}

\noindent
For the convenience of the reader, we include an appendix where the definition and some basic results about regularly varying 
functions in Karamata sense are recalled.

\section{An upper bound for the monodromy matrix}

In this section we give a generic upper bound for $\log\|W_H(z)\|$.
The functions $\ms d_l,\ms d_\phi$ used below play the role of well-behaved comparison functions for the lengths and
for the differences of consecutive angles of the Hamiltonian. 

\begin{Definition}
\label{B70}
	Let $\ms d_l\DF[1,\infty)\to(0,\infty)$ and $\ms d_\phi\DF[1,\infty)\to(0,1]$ be measurable with 
	$\ms d_l\asymp 1\asymp\ms d_\phi$ locally. Then we denote
	\begin{align}
		\label{B44}
		\ms k(R)\DE &\, 
		\sup\Big\{t\in[1,\infty)\DS \sup_{1\leq s\leq t}\frac 2{R(\ms d_l\ms d_\phi)(s)}\leq 1\Big\}
		,
		\\
		\label{B45}
		\ms h(R)\DE &\,
		\sup\Big\{t\in[1,\infty)\DS \sup_{1\leq s\leq t}\frac{\ms d_\phi(s)}{R\ms d_l(s)}\leq 1\Big\}
		,
	\end{align}
	where $R\in[\frac 2{(\ms d_l\ms d_\phi)(1)},\infty)$ and $R\in[\frac{\ms d_\phi(1)}{\ms d_l(1)},\infty)$, respectively.
	Further, we set 
	\[
		g(t,R)\DE
		\begin{cases}
			\log\big(R(\ms d_l\ms d_\phi)(t)\big) \CAS 1\leq t<\ms k(R),
			\\[1mm]
			R^{\frac 12}(\ms d_l\ms d_\phi)^{\frac 12}(t) \CAS \ms k(R)\leq t<\ms h(R),
			\\[1mm]
			R\ms d_l(t) \CAS \ms h(R)\leq t,
		\end{cases}
	\]
	where $(t,R)\in[1,\infty)\times[\frac 2{(\ms d_l\ms d_\phi)(1)},\infty)$, and 
	\begin{equation}
	\label{B48}
		\ms g(t,R)=\int_1^t g(s,R)\DD s
		.
	\end{equation}
\end{Definition}

\noindent
Note here that $\ms d_\phi(t)\leq 1$, and hence $\ms k(R)\leq\ms h(R)$ for all $R\geq\frac 2{(\ms d_l\ms d_\phi)(1)}$. 

\begin{theorem}
\label{B8}
	Let $(l_j)_{j=1}^\infty$ be a summable sequence of positive numbers, and $(\phi_j)_{j=1}^\infty$ a sequence of real
	numbers. Denote by $H$ the Hamburger Hamiltonian with these lengths and angles, and let $W_H$ be its monodromy matrix. 

	Let $\psi\in\bb R $ and let $\ms d_l,\ms d_\phi,\ms c_l,\ms c_\phi\DF[1,\infty)\to(0,\infty)$ be measurable and 
	nonincreasing, with $\ms d_\phi\leq 1$ and $\ms c_\phi\leq\ms c_l$. 
	Assume that
	\begin{align}
		\label{B13}
		& \forall j\in\bb N\DP 
		l_j\leq\ms d_l(j)\ \wedge\ |\sin(\phi_{j+1}-\phi_j)|\leq\ms d_\phi(j)
		,
		\\
		\label{B14}
		& \forall N\in\bb N\DP 
		\sum_{j=N+1}^\infty l_j\leq\ms c_l(N)\ \wedge\ \sum_{j=N+1}^\infty l_j\sin^2(\phi_j-\psi)\leq\ms c_\phi(N)
		.
	\end{align}
	Let $\ms g(t,R)$ be as in \cref{B48}, and set 
	\begin{align}
	\label{B49}
		L&(t,R) \DE 
		1+\log^+ R+\log^+\frac{\ms c_l(\lceil t \rceil)}{\ms c_\phi(\lceil t \rceil)}
		+\log^+\frac{\ms d_l(1)}{\ms d_\phi(1)}
		\\
	\nonumber
		&+\log^+\frac{\ms d_l(\min \{\lceil t \rceil,\lfloor\ms h(R) \rfloor\})}{\ms d_\phi(\min 
		\{\lceil t \rceil,\lfloor\ms h(R) \rfloor\})}
		+\sum_{j=1}^{ \min \{\lceil t \rceil,\lfloor\ms h(R) \rfloor\}-1 } 
		\bigg| \log\bigg( \raisebox{4pt}{$\frac{\ms d_\phi(j)}{\ms d_l(j)}$}\mkern-3mu\bigg/\mkern-5mu
		\raisebox{-5pt}{$\frac{\ms d_\phi(j+1)}{\ms d_l(j+1)}$}\bigg) \bigg|
		.
	\end{align}
	Then we have, for all $R\geq\frac 2{(\ms d_l\ms d_\phi)(1)}$,
	\begin{equation}
	\label{B32}
		\log\Big(\max_{|z|=R}\|W_H(z)\|\Big)\leq 9\cdot
		\inf_{t\geq 1}\Big(\max\Big\{\ms g\big(t,R\big),R(\ms c_l\ms c_\phi)^{\frac 12}(t)\Big\}+L(t,R)\Big)
		.
	\end{equation}
\end{theorem}

\begin{Remark}
\label{B17}
	The terms appearing in the upper bound in \cref{B32} have the following meaning:
	\begin{Enumerate}
	\item $\ms g(t,R)$ estimates the contribution of the first $\lceil t \rceil$ intervals of $H$;
	\item $R(\ms c_l\ms c_\phi)^{\frac 12}(t)$ estimates the contribution of the remaining intervals of $H$;
	\item $L(t,R)$ is usually a remainder term.
	\end{Enumerate}
\end{Remark}

\noindent
We mention two possible scenarios. If lengths and angle differences are nonincreasing but decay slowly, choosing
$\ms d_l(j),\ms d_\phi (j)$ such that equality holds in \cref{B13} makes $\ms g(t,R)$ a rather precise bound for the 
contribution of the first $\lceil t \rceil$ intervals. On the other hand, if lengths and angle differences behave 
irregularly, then $\ms g(t,R)$ likely overestimates the contribution of the first $\lceil t \rceil$ intervals.
However, choosing $\ms c_l,\ms c_\phi$ such that equality holds in \cref{B14}, the bound \cref{B32} is still a good estimate for 
$\log \|W_H(z)\|$ as long as the decay is sufficiently fast.

If $\lim_{t\to\infty}(\ms c_l\ms c_\phi)(t)>0$, the bound \cref{B32} is $\gtrsim R$ and thus trivial: we know
that $W_H$ is of minimal exponential type by the classical Krein-de~Branges formula. 
Hence we may safely assume that $\lim_{t\to\infty}(\ms c_l\ms c_\phi)(t)=0$ whenever this is convenient.

The following lemma hints at a way to evaluate the upper bound in \Cref{B8}.

\begin{lemma}
\label{B46}
	Let $\ms d_l\DF[1,\infty)\to(0,\infty)$ and $\ms d_\phi\DF[1,\infty)\to(0,1]$ be measurable with 
	$\ms d_l\asymp 1\asymp\ms d_\phi$ locally, and let $\ms c_l,\ms c_\phi$
	be continuous and nonincreasing with $\lim\limits_{t\to\infty}(\ms c_l\ms c_\phi)(t)=0$. 
	Then the equation
	\begin{equation}
	\label{B57}
		\ms g(t,R)=R(\ms c_l\ms c_\phi)^{\frac 12}(t),
	\end{equation}
	has a unique solution $T(R)$, and 
	\[
		\min_{t\geq 1}\max\Big\{\ms g\big(t,R\big),R(\ms c_l\ms c_\phi)^{\frac 12}(t)\Big\}
		=\ms g(T(R),R)=R(\ms c_l\ms c_\phi)^{\frac 12}(T(R))
		.
	\]
	In addition, $\lim_{R \to \infty} T(R)=\infty$.
\end{lemma}	
\begin{proof}
	For each fixed $R$ the function $t\mapsto\ms g(t,R)$ is continuous and increasing while 
	$t\mapsto R(\ms c_l\ms c_\phi)^{\frac 12}(t)$ is
	continuous and nonincreasing. Moreover, $\ms g(1,R)=0<R(\ms c_l\ms c_\phi)^{\frac 12}(1)$ and 
	$\lim_{t \to \infty} R(\ms c_l\ms c_\phi)^{\frac 12}(t)=0$. By the intermediate value theorem, \cref{B57} has a unique 
	solution $T(R)$.

	If $\liminf_{R \to \infty} T(R)<\infty$, choose a sequence $(R_n)_{n \in \bb N}$ such that 
	$(T(R_n))_{n \in \bb N}$ is bounded. Using the crude estimate
	\[
		\ms g(t,R)\leq t\log\big(R(\ms d_l\ms d_\phi)(1)\big)+\sqrt 2\, t +t \lesssim t \log R
		,
	\]
	this leads to the contradiction 
	\[
		R(\ms c_l\ms c_\phi)^{\frac 12}(T(R_n))=\ms g(T(R_n),R) \lesssim 
		T(R_n) \cdot \log R \ll R \asymp  R(\ms c_l\ms c_\phi)^{\frac 12}(T(R_n)).
	\]
\end{proof}

\noindent 
In many situations, we are going to first determine $T(R)$, and then show that $L(T(R),R)$ is small. 
The bound \cref{B32} is then asymptotically equal to $R(\ms c_l\ms c_\phi)^{\frac 12}(T(R))$.

Before we go into the proof of \Cref{B8}, let us discuss a particular situation. 

\begin{corollary}
\label{B66}
	Let $(l_j)_{j=1}^\infty$ be a summable sequence of positive numbers, and $(\phi_j)_{j=1}^\infty$ a sequence of real
	numbers. Denote by $H$ the Hamburger Hamiltonian with these lengths and angles, and let $W_H$ be its monodromy matrix. 

	Assume that we have $\delta_l,\delta_\phi,\gamma_l,\gamma_\phi\geq 0$, such that 
	\begin{align*}
		& \forall j\in\bb N\DP 
		l_j\lesssim j^{-\delta_l}\ \wedge\ |\sin(\phi_{j+1}-\phi_j)|\lesssim j^{-\delta_\phi}
		,
		\\
		& \forall N\in\bb N\DP 
		\sum_{j=N+1}^\infty l_j\lesssim N^{-\gamma_l}\ \wedge\ 
		\sum_{j=N+1}^\infty l_j\sin^2(\phi_j-\psi)\lesssim N^{-\gamma_\phi}
		,
	\end{align*}
	i.e., we are in the situation of \Cref{B8} with
	\[
		\ms d_l(t)=ct^{-\delta_l},\ \ms d_\phi(t)=ct^{-\delta_\phi},\ 
		\ms c_l(t)=ct^{-\gamma_l},\ \ms c_\phi(t)=ct^{-\gamma_\phi}
		,
	\]
	where $c$ is some positive constant. Set $\delta\DE\delta_l+\delta_\phi$, 
	$\gamma\DE\frac 12(\gamma_l+\gamma_\phi)$, and assume that $\delta>0$ and $\gamma>0$. Then we have, for sufficiently
	large $R$, the following bounds for $T(R)$ (cf. \Cref{B46}), $\log \max_{|z|=R} \|W_H(z)\|$ and the order $\rho_H$ of 
	$W_H$.
	\[
		\begin{array}{lll|l|l|l}
			\multicolumn{3}{l|}{\text{Data satisfies}} 
			& T(R) \asymp 
			& \log \max\limits_{|z|=R} \|W_H(z)\|\lesssim
			& \rho_H \leq 
			\raisebox{-6pt}{\rule{0pt}{17pt}}
			\\
			\hline
			\hline
			\multicolumn{3}{l|}{\delta<1+\gamma}
			& \big(\frac R{\log R}\big)^{\frac 1{1+\gamma}} 
			& R^{\frac 1{1+\gamma}}(\log R)^{\frac\gamma{1+\gamma}}
			\raisebox{-9pt}{\rule{0pt}{25pt}}
			& \frac 1{1+\gamma}
			\\
			\hline
			\multicolumn{3}{l|}{\delta=1+\gamma}
			& R^{\frac 1\delta} 
			& R^{\frac 1\delta}
			& \frac{1}{\delta}
			\raisebox{-9pt}{\rule{0pt}{25pt}}
			\\
			\hline
			1+\gamma<\delta,
			& \delta>2,
			& 
			& R^{\frac {\delta-1}{\gamma \delta}}
			& R^{\frac 1\delta}
			& \frac{1}{\delta}
			\raisebox{-9pt}{\rule{0pt}{25pt}}
			\\
			\cline{2-6}
			& \delta=2,
			& 
			& \big( \frac {R^{\frac 12}}{\log R} \big)^{\frac 1\gamma }
			& R^{\frac 12}\log R
			& \frac{1}{2}
			\raisebox{-9pt}{\rule{0pt}{25pt}}
			\\
			\cline{2-6}
			& \delta<2,
			& \delta_l\leq 1+\gamma
			& R^{\frac 1{2-\delta+2\gamma}}
			& R^{\frac{2-\delta+\gamma}{2-\delta+2\gamma}}
			& \frac{2-\delta+\gamma}{2-\delta+2\gamma}
			\raisebox{-12pt}{\rule{0pt}{31pt}}
			\\
			\cline{3-6}
			& 
			& 1+\gamma<\delta_l
			& R^{\frac {\delta_l-1}{\gamma (\delta_l-\delta_\phi)}}
			& R^{\frac{1-\delta_\phi}{\delta_l-\delta_\phi}}
			& \frac{1-\delta_\phi}{\delta_l-\delta_\phi}
			\raisebox{-12pt}{\rule{0pt}{31pt}}
			\\
			\hline
		\end{array}
	\]
\end{corollary}
\begin{proof}
	To start with we observe
	\begin{align*}
		& \ms k(R)=\Big(\frac{c^2R}2\Big)^{\frac 1\delta}
		,\quad \ms h(R)=
		\begin{cases}
			R^{\frac 1{\delta_l-\delta_\phi}} \CAS \delta_l>\delta_\phi,
			\\
			\infty \CASO,
		\end{cases}
		\\
		& R(\ms c_l\ms c_\phi)^{\frac 12}(t)=Rt^{-\gamma}
		.
	\end{align*}
	The main effort is to evaluate $\ms g(t,R)$, and this is where the case distinctions come from. 
	\begin{Ilist}
	\item Range $t\in[1,\ms k(R)]$:
		\begin{align}
		\label{B58}
		\begin{split}
			\ms g(t,R)= &\, \int_1^t\log\big(c^2Rs^{-\delta}\big)\DD s=t\log(c^2R)-\delta\int_1^t\log s\DD s
			\\
			= &\, t\log(c^2R)-\delta\big(t\log t-t+1\big)=t\log\big(c^2Rt^{-\delta}\big)+\delta t-\delta
			.
			\end{split}
		\end{align}
		In particular,
		\begin{align*}
			\ms g(\ms k(R),R)= &\, 
			\ms k(R)\log\big(\underbrace{R\ms k(R)^{-\delta}}_{=2}\big)+\delta\ms k(R)-\delta
			\\
			= &\, (\log 2+\delta)\ms k(R)-\delta\asymp\ms k(R)
			.
		\end{align*}
	\item Range $t\in[\ms k(R),\ms h(R)]$:
		\[
			\int_{\ms k(R)}^t s^{-\frac\delta 2}\DD s=
			\begin{cases}
				\frac 1{1-\frac \delta 2}\big(t^{1-\frac\delta 2}-\ms k(R)^{1-\frac\delta 2}\big)
				\CAS \delta<2
				,
				\\
				\log\frac t{\ms k(R)} \CAS \delta=2
				,
				\\
				\frac 1{\frac \delta 2-1}\big(\ms k(R)^{1-\frac\delta 2}-t^{1-\frac\delta 2}\big)
				\CAS \delta>2
				.
			\end{cases}
		\]
		Let $\sigma>\frac 1\delta$ and $\sigma\leq\frac 1{\delta_l-\delta_\phi}$ in case $\delta_l>\delta_\phi$. 
		Then 
		\[
			\ms g(R^\sigma,R)\asymp R^{\frac 1\delta}+R^{\frac 12}\cdot
			\begin{cases}
				R^{\sigma(1-\frac\delta 2)} \CAS \delta<2
				,
				\\
				\log R \CAS \delta=2
				,
				\\
				R^{\frac 1\delta (1-\frac\delta 2)} \CAS \delta>2
				.
			\end{cases}
		\]
		Observe that for $\delta <2$
		\[
			\frac 1\delta=\frac 12+\frac 1\delta\big(1-\frac\delta 2\big)
			<
			\frac 12+\sigma\big(1-\frac\delta 2\big)
			.
		\]
		Hence, 
		\begin{align}
		\label{B15}
			\ms g(R^\sigma,R)\asymp
			\begin{cases}
				R^{\frac 12+\sigma(1-\frac\delta 2)} \CAS \delta<2
				,
				\\
				R^{\frac 12}\log R \CAS \delta=2
				,
				\\
				R^{\frac 1\delta} \CAS \delta>2
				.
			\end{cases}
		\end{align}
		\item Range $t\in[\ms h(R),\infty)$:
			\\
			This case can occur only if $\delta_l>\delta_\phi$ (since otherwise $\ms h(R)=\infty$). 
			\[
				\int_{\ms h(R)}^t s^{-\delta_l}\DD s=
				\begin{cases}
					\frac 1{1-\delta_l}\big(t^{1-\delta_l}-\ms h(R)^{1-\delta_l}\big)
					\CAS \delta_l<1
					,
					\\
					\log\frac t{\ms h(R)} \CAS \delta_l=1
					,
					\\
					\frac 1{\delta_l-1}\big(\ms h(R)^{1-\delta_l}-t^{1-\delta_l}\big)
					\CAS \delta_l>1
					.
				\end{cases}
			\]
			Let $\sigma>\frac 1{\delta_l-\delta_\phi}$. Then 
			\[
				\ms g(R^\sigma,R)\asymp
				\begin{cases}
					R^{\frac{1-\delta_\phi}{\delta_l-\delta_\phi}}+R^{1+\sigma(1-\delta_l)}
					\CAS \delta<2,\delta_l<1
					,
					\\
					R^{\frac{1-\delta_\phi}{\delta_l-\delta_\phi}}+R\log R
					\CAS \delta<2,\delta_l=1
					,
					\\
					R^{\frac{1-\delta_\phi}{\delta_l-\delta_\phi}}+ 
					R^{1+\frac{1-\delta_l}{\delta_l-\delta_\phi}}
					\CAS \delta<2,\delta_l>1
					,
					\\
					R^{\frac 12}\log R+R^{1+\frac{1-\delta_l}{\delta_l-\delta_\phi}}
					\CAS \delta=2\quad\text{($\Rightarrow\delta_l>1$)}
					,
					\\
					R^{\frac 1\delta}+R^{1+\frac{1-\delta_l}{\delta_l-\delta_\phi}}
					\CAS \delta>2\quad\text{($\Rightarrow\delta_l>1$)}
					.
				\end{cases}
			\]
			Observe that 
			\[
				\delta_l<1,\sigma>\frac 1{\delta_l-\delta_\phi}\ \Rightarrow\ 
				\frac{1-\delta_\phi}{\delta_l-\delta_\phi}<1+\sigma(1-\delta_l)
			\]
			since 
			\begin{align*}
				\frac{1-\delta_\phi}{\delta_l-\delta_\phi}< &\, 1+\sigma(1-\delta_l)
				\\
				\Leftrightarrow &\, 1-\delta_\phi<\delta_l-\delta_\phi+\sigma(1-\delta_l)(\delta_l-\delta_\phi)
				\\
				\Leftrightarrow &\, 1-\delta_l<\sigma(1-\delta_l)(\delta_l-\delta_\phi)
				\\
				\Leftrightarrow &\, 1<\sigma(\delta_l-\delta_\phi)
			\end{align*}
			Next, 
			\[
				\delta\geq 2\ \Rightarrow\ \frac{1-\delta_\phi}{\delta_l-\delta_\phi}\leq\frac 1\delta
			\]
			since 
			\begin{align*}
				\frac{1-\delta_\phi}{\delta_l-\delta_\phi}\leq &\, \frac 1\delta
				\\
				\Leftrightarrow &\, (1-\delta_\phi)(\delta_l+\delta_\phi)\leq\delta_l-\delta_\phi
				\\
				\Leftrightarrow &\, 0\leq\delta_\phi(\delta-2)
			\end{align*}
			It follows that 
			\begin{align}
			\label{B43}
				\ms g(R^\sigma,R)\asymp
				\begin{cases}
					R^{1+\sigma(1-\delta_l)} \CAS \delta<2,\delta_l<1
					,
					\\
					R\log R \CAS \delta<2,\delta_l=1
					,
					\\
					R^{\frac{1-\delta_\phi}{\delta_l-\delta_\phi}} \CAS \delta<2,\delta_l>1
					,
					\\
					R^{\frac 12}\log R \CAS \delta=2
					,
					\\
					R^{\frac 1\delta} \CAS \delta>2
					.
				\end{cases}
			\end{align}
	\end{Ilist}
	Denote by $\tau (R)$ the function (whose form depends on the considered case) given in the second column of the table in
	\Cref{B66}. Using \cref{B58}, \cref{B15}, and \cref{B43}, we see that 
	\begin{equation}
	\label{B64}
		\ms g(\tau (R),R) \asymp R\tau (R)^{-\gamma} = R (\ms c_l\ms c_\phi)^{\frac 12}(\tau (R))
	\end{equation}
	and hence $T(R) \asymp \tau (R)$. Moreover, the upper bound \cref{B32} for $\log \max_{|z|=R} \|W_H(z)\|$ is $\asymp$
	to all terms in \cref{B64}.

	It remains to check that $L(T(R),R)$ is small. First, observe that the function
	$\frac{\ms d_\phi}{\ms d_l}$ is just a power and thus monotone. Hence, the last summand of $L(t,R)$ is a telescoping
	sum. Using the facts that also $\frac{\ms c_l}{\ms c_\phi}$ is a power and that $T(R)$ is, in each case, bounded from
	below and above by some positive power, we have $L(T(R),R)\lesssim\log R$. 
\end{proof}

\noindent
\Cref{B66} already shows that our present results may improve drastically upon our previous work
\cite{pruckner.romanov.woracek:jaco}. The main feature that enables this is the use of $\ms c_\phi$: the present results exploit
convergence of angles more efficiently.

\begin{example}
\label{B6}
	Let $\alpha >1$ and $\beta \geq 0$, and let $H$ be the Hamburger Hamiltonian with lengths and angles 
	\[
		l_j\DE j^{-\alpha},\quad \phi_j\DE(-1)^jj^{-\beta}
		.
	\]
	Considering the expressions \cref{B13} and \cref{B14}, we see that \Cref{B66} can be applied with 
	\[
		\delta_l\DE\alpha,\ \delta_\phi\DE\beta,\ \gamma_l\DE \alpha-1,\ \gamma_\phi\DE\alpha+2\beta-1
		.
	\]
	Since $\delta=\alpha+\beta$ and $\gamma=\alpha+\beta-1$, this yields 
	\[
		\log\Big(\max{|z|=R}\|W_H(z)\|\Big)\lesssim R^{\frac 1{\alpha+\beta}}
		,
	\]
	and hence $\rho_H \leq \frac 1{\alpha+\beta}$. Now recall \cite[Example~2.23]{pruckner.romanov.woracek:jaco}, where we 
	saw that the order $\rho_H$ is at least $\frac 1{\alpha+\beta}$. Thus, $\rho_H=\frac 1{\alpha+\beta}$. 

	Let us compare this with what we can obtain from \cite{pruckner.romanov.woracek:jaco}. For $\alpha+\beta \geq 2$, the
	formula $\rho_H=\frac 1{\alpha+\beta}$ was already established in \cite[Example~2.23]{pruckner.romanov.woracek:jaco}. 
	For $\alpha+\beta<2$, we had obtained the upper bound $\frac{1-\beta}{\alpha-\beta}$ in that place. However, this 
	did not take into account that angles converge. If we want to take convergence into account, we should use
	\cite[Corollary~2.9]{pruckner.romanov.woracek:jaco}. The quantities used there identify as 
	\[
		\Delta_l^+=\alpha,\quad \Delta_\phi^*=\Lambda^*=\beta
		,
	\]
	and the bound for order so obtained thus is $\frac{1-\frac\beta 2}\alpha$. 
	If $\beta >0$, this lies strictly between $\frac{1}{\alpha+\beta}$ and $\frac{1-\beta}{\alpha-\beta}$, since
	\begin{align*}
		& \frac{1-\beta}{\alpha-\beta}-\frac{1-\frac\beta 2}\alpha=
		\frac{\frac\beta 2(2-\alpha-\beta)}{(\alpha-\beta)\alpha}>0
		,
		\\
		& \frac{1-\frac\beta 2}\alpha-\frac 1{\alpha+\beta}=
		\frac{\frac\beta 2(2-\alpha-\beta)}{\alpha(\alpha+\beta)}>0
		.
	\end{align*}
\end{example}

\subsection{Proof of \Cref{B8}}

Denote by $W(s,t;z)$ the unique solution of the initial value problem 
\[
	\left\{
	\begin{array}{l}
		\frac{\partial}{\partial t} W(s,t;z)J=zW(s,t;z)H(t),\qquad t\in(0,L)\text{ a.e.},
		\\[2mm]
		W(s,s;z)=I.
	\end{array}
	\right.
\]
and observe the multiplicativity property 
\[
	W(s,t;z)W(t,u;z)=W(s,u;z)
	.
\]
Our method for the proof of \Cref{B8} is to fix $t\geq 1$ and estimate $\|W(0,x_{\lceil t \rceil};z)\|$ and 
$\|W(x_{\lceil t\rceil},L;z)\|$ separately. The contribution of the first $\lceil t \rceil$ intervals is estimated using the 
explicit form of the fundamental solution on indivisible intervals. On the remaining part, we use Gr\"onwall's Lemma. 
This method is a refined version of an idea from \cite{romanov:2017} and its improved variant from \cite{pruckner.woracek:sinqA}. 

Set $W_j(z):=W(x_{j-1},x_j;z)$, then by a direct computation 
\[
	W_j(z)=I+zl_j \xi_{\phi_j} \xi_{\phi_j}^T J
	.
\]
Moreover, multiplicativity yields
\begin{equation}
\label{B18}
	W(0,x_N;z)=W_1(z) \cdot\ldots\cdot W_N(z), \qquad N \in \bb N
	.
\end{equation}
We often use the matrices introduced in \cite[Definition~4.4]{pruckner.woracek:sinqA}:
for $a > 0$ and $\psi \in \bb R$, denote
\[
	\Omega(a,\psi) \DE \begin{pmatrix}
	a & 0 \\
	0 & a^{-1}
	\end{pmatrix} \exp (-\psi J)= \begin{pmatrix}
	a & 0 \\
	0 & a^{-1}
	\end{pmatrix} \begin{pmatrix}
	\cos \psi & -\sin \psi \\
	\sin \psi & \cos \psi
	\end{pmatrix}.
\]
Recall the following properties from \cite[Lemma 4.6]{pruckner.woracek:sinqA}.

\begin{Lemma}
\label{B20}
	Let $a,b>0$ and $\phi,\psi \in \bb R$. 
	\begin{Enumerate}
	\item $\|\Omega (a,\psi) \|=\|\Omega (a,\psi)^{-1} \|=\max \{a,a^{-1}\}$;
	\item $\|\Omega (a,\psi) \xi_\phi \xi_\phi^T J \Omega (a,\psi)^{-1} \|=a^2 \cos^2 (\phi - \psi) 
		+\frac{1}{a^2} \sin^2 (\phi - \psi)$;
	\item $\|\Omega (a,\psi) \Omega (b,\phi )^{-1} \| \leq \max \big\{\frac{a}{b},\frac{b}{a} \big\}|\cos (\phi-\psi)| + 
		\max \big\{ab,\frac{1}{ab} \big\}|\sin (\phi-\psi)|$.
	\end{Enumerate}
\end{Lemma}

\noindent 
We use $\Omega (a_j,\psi_j)$ to rotate and dilate $W_j(z)$ such that the norm becomes small. In the following
lemma, we make the obvious choice $\psi_j=\phi_j$, but retain the freedom of choosing suitable parameters $a_j$.

\begin{lemma}
\label{B21}
	Let $N \in \bb N$ and let numbers $a_j \in (0,1]$, $j=1,\ldots,N$, be given. Then, for each $z \in \bb C$, we have 
	the estimate
	\begin{align}
		\label{B91}
		\|W(0&,x_N;z)\| \leq \frac{1}{a_1a_N}\cdot\prod_{j=1}^N \Big(1+|z|l_j a_j^2 \Big) 
		\\
		\nonumber
		\cdot &\, 
		\prod_{j=1}^{N-1} \bigg(
		\max\Big\{\frac{a_j}{a_{j+1}},\frac{a_{j+1}}{a_j}\Big\}\cdot\big|\cos\big(\phi_j-\phi_{j+1}\big)\big|
		+\frac{|\sin(\phi_j-\phi_{j+1})|}{a_ja_{j+1}}
		\bigg).
	\end{align}
\end{lemma}
\begin{proof}
	With $\Omega_j:=\Omega(a_j,\phi_j)$ we have 
	\[
		\Omega_j W_j(z) \Omega_j^{-1}=I+zl_j \Omega_j \xi_{\phi_j} \xi_{\phi_j}^T J \Omega_j^{-1},
	\]
	and we can use \Cref{B20}, $(ii)$, to estimate
	\[
		\|\Omega_j W_j(z) \Omega_j^{-1}\| \leq 1+|z|l_ja_j^2.
	\]
	Expanding on \cref{B18},
	\[
		W(0,x_N;z)=\Omega_1^{-1} \cdot (\Omega_1W_1(z)\Omega_1^{-1}) \cdot \Omega_1\Omega_2^{-1} \cdot\ldots\cdot
		(\Omega_NW_N(z)\Omega_N^{-1})\cdot \Omega_N.
	\]
	Using submultiplicativity of the norm and \Cref{B20}, $(i),(iii)$, we arrive at the desired estimate.
\end{proof}

\noindent 
The next lemma is the ingredient of \Cref{B8} that accounts for the term $\ms g(t,R)$ in \cref{B32}.

\begin{lemma}
\label{B25}
	Let $(l_j)_{j=1}^\infty$ be a summable sequence of positive numbers, and $(\phi_j)_{j=1}^\infty$ a sequence of real
	numbers. Let $H$ be the Hamburger Hamiltonian with these lengths and angles. 
	Let $\ms d_l\DF[1,\infty)\to(0,\infty)$ and $\ms d_\phi\DF[1,\infty)\to(0,1]$ be measurable and nonincreasing, and 
	assume that
	\[
		\forall j\in\bb N\DP 
		l_j\leq\ms d_l(j)\ \wedge\ |\sin(\phi_{j+1}-\phi_j)|\leq\ms d_\phi(j)
		.
	\]
	Then, for every $(N,R) \in \bb N \times \big[\frac{2}{(\ms d_l \ms d_\phi)(1)},\infty\big)$, 
	\[
		\log\Big(\max_{|z|=R}\|W(0,x_N;z)\|\Big)\leq
		\smfrac 72\cdot\Big[\ms g(N,R)+\lambda\big(\min\{N,\ms h(R)\},R\big)\Big]
	\]
	where
	\[
		\lambda(t,R)\DE 1+\log^+R+\log^+\frac{\ms d_l(1)}{\ms d_\phi(1)}
		+ \log^+ \frac{\ms d_l(\lfloor t \rfloor)}{\ms d_\phi (\lfloor t \rfloor)} 
		+\sum_{j=1}^{\lfloor t \rfloor -1} \bigg| \log \frac{\ms d_\phi(j)\ms d_l(j+1)}{\ms d_l(j)\ms d_\phi (j+1)} 
		\bigg|.
	\]
\end{lemma}
\begin{proof}
	Fix $z$ with $|z|=R$. In the first step we estimate $\log\|W(0,x_N;z)\|$ when $N\leq\ms h(R)$. This is done by an 
	application of \Cref{B21} with 
	\[
		a_j:= \Big(\frac{\ms d_\phi (j)}{R \ms d_l(j)} \Big)^{\frac 14}, \qquad j=1,\ldots , N.
	\]
	Since $N\leq\ms h(R)$ we have $a_j \in (0,1]$. The factors in \cref{B91} are treated separately:
	\[
		\frac{1}{a_1a_N}=R^{\frac 12}\cdot\Big(\frac{\ms d_l(1)}{\ms d_\phi(1)}\Big)^{\frac 14}\cdot
		\Big(\frac{\ms d_l(N)}{\ms d_\phi(N)}\Big)^{\frac 14}
		,
	\]
	\[
		1+|z|l_ja_j^2 \leq 1+R \ms d_l(j) a_j^2 = 1+\big(R(\ms d_l \ms d_\phi)(j) \big)^{\frac 12}
		,
	\]
	\begin{align*}
		\max &\, \Big\{\frac{a_j}{a_{j+1}},\frac{a_{j+1}}{a_j}\Big\}\big|\cos\big(\psi_j-\psi_{j+1}\big)\big|
		+\frac{|\sin(\psi_j-\psi_{j+1})|}{a_ja_{j+1}}
		\\
		= &\, \max\Big\{\frac{a_j}{a_{j+1}},\frac{a_{j+1}}{a_j}\Big\} \cdot
		\bigg(|\cos(\phi_{j+1}-\phi_j)|+\frac{|\sin(\phi_{j+1}-\phi_j)|}{\max \{a_j^2,a_{j+1}^2\}}\bigg)
		\\
		\leq &\,\max\Big\{\frac{a_j}{a_{j+1}},\frac{a_{j+1}}{a_j}\Big\}\cdot
		\bigg(1+\frac{\ms d_\phi (j)}{a_j^2}\bigg)
		\\
		= &\, \Big[\max\Big\{\frac{\ms d_\phi(j)\ms d_l(j+1)}{\ms d_l(j)\ms d_\phi (j+1)},
		\frac{\ms d_l(j)\ms d_\phi (j+1)}{\ms d_\phi(j)\ms d_l(j+1)}\Big\}\Big]^{\frac 14} 
		\cdot\Big(1+\big(R(\ms d_l \ms d_\phi)(j) \big)^{\frac 12}\Big)
		.
	\end{align*}
	Now \cref{B91} yields
	\begin{align*}
		\log\|W(0 & ,x_N;z)\|\leq \frac 12\log^+R+\frac 14\log^+\frac{\ms d_l(1)}{\ms d_\phi(1)}
		+\frac 14\log^+\frac{\ms d_l(N)}{\ms d_\phi(N)}
		\\
		& +\sum_{j=1}^N\log\Big(1+\big(R(\ms d_l\ms d_\phi)(j)\big)^{\frac 12}\Big)
		\\
		& +\frac 14\sum_{j=1}^{N-1}\bigg|\log\frac{\ms d_\phi(j)\ms d_l(j+1)}{\ms d_l(j)\ms d_\phi(j+1)}\bigg|
		+\sum_{j=1}^{N-1}\log\Big(1+\big(R(\ms d_l\ms d_\phi)(j)\big)^{\frac 12}\Big)
		.
	\end{align*}
	Since $\ms d_l,\ms d_\phi$ are nonincreasing, we can estimate sums by integrals as 
	\begin{multline*}
		\sum_{j=2}^N\log\bigg(1+\big(R(\ms d_l\ms d_\phi)(j)\big)^{\frac 12}\bigg)
		\leq\int\limits_1^N\log\Big(1+\big(R(\ms d_l\ms d_\phi)(s)\big)^{\frac 12}\Big)\DD s
		\\
		\leq\smfrac 32\mkern5mu\cdot\mkern-35mu\int_1\limits^{\min\{N,\ms k(R)\}}\mkern-30mu
		\log\big(R(\ms d_l\ms d_\phi)(s)\big)\DD s\mkern5mu+
		\mkern-30mu\int\limits_{\min\{N,\ms k(R)\}}^N\mkern-30mu\big(R(\ms d_l\ms d_\phi)(s)\big)^{\frac 12}\DD s 
		.
	\end{multline*}
	Noting that 
	\begin{align*}
		\log\Big(1+\big(R(\ms d_l\ms d_\phi)(1)\big)^{\frac 12}\Big)\leq &\, 
		\frac 32\Big(\log^+R+\log^+(\ms d_l\ms d_\phi)(1)\Big)
		\\
		\leq &\, 
		\frac 32\Big(\log^+R+\log^+\frac{\ms d_l(1)}{\ms d_\phi(1)}\Big)
		,
	\end{align*}
	leads to 
	\begin{align}
	\nonumber
		\log\| & W(0,x_N;z)\|\leq \frac 72\log^+R+\frac{13}4\log^+\frac{\ms d_l(1)}{\ms d_\phi(1)}
		+\frac 14\log^+\frac{\ms d_l(N)}{\ms d_\phi(N)}
		\\
	\nonumber
		& +\frac 14\sum_{j=1}^{N-1}\bigg|\log\frac{\ms d_\phi(j)\ms d_l(j+1)}{\ms d_l(j)\ms d_\phi(j+1)}\bigg| 
		\\
	\label{B1}
		&
		+3\mkern5mu\cdot\mkern-35mu\int_1\limits^{\min\{N,\ms k(R)\}}\mkern-30mu
		\log\big(R(\ms d_l\ms d_\phi)(s)\big)\DD s\mkern5mu+
		2\mkern5mu\cdot\mkern-35mu
		\int\limits_{\min\{N,\ms k(R)\}}^N\mkern-30mu\big(R(\ms d_l\ms d_\phi)(s)\big)^{\frac 12}\DD s
		.
	\end{align}
	We see that 
	\[
		\log\|W(0,x_N;z)\|\leq \smfrac 72\cdot\big[\ms g(N,R)+\lambda(N,R)\big]
		.
	\]
	Consider now the case that $N>\ms h(R)$. Since $\|\xi_{\phi_j} \xi_{\phi_j}^T J\|=1$, we see that
	\begin{align*}
		\|W(x_{\lfloor \ms h(R) \rfloor} & ,x_N;z)\| =
		\bigg\|\prod_{j=\lfloor \ms h(R) \rfloor+1}^N W(x_{j-1},x_j;z)\bigg\|
		\\
		\leq &\, \prod_{j=\lfloor \ms h(R) \rfloor+1}^N\Big\|I+zl_j \xi_{\phi_j} \xi_{\phi_j}^T J\Big\|
		\leq \prod_{j=\lfloor \ms h(R) \rfloor+1}^N \big(1+R\ms d_l(j)\big).
	\end{align*}
	Noting that 
	\[
		R\ms d_l\big(\lfloor\ms h(R)\rfloor+1\big)<1
		,
	\]
	and estimating the sum by an integral we arrive at
	\begin{align*}
		\log\|W(x_{\lfloor\ms h(R)\rfloor},x_N;z)\| 
		\mkern10mu\leq & \mkern-10mu\sum_{j=\lfloor\ms h(R)\rfloor+1}^N\mkern-20mu\log\Big(1+R\ms d_l(j)\Big)
		\\
		< &\, \log 2\mkern10mu+\mkern-10mu\int\limits_{\lfloor \ms h(R)\rfloor+1}^N\mkern-20mu R \ms d_l(s) \DD s.
	\end{align*}
	We combine this with \cref{B1} to obtain 
	\begin{align*}
		\log\|W(0,x_N;z)\|\leq &\, 
		\log\|W(0,x_{\lfloor\ms h(R)\rfloor};z)\|+\log\|W(x_{\lfloor\ms h(R)\rfloor},x_N;z)\|
		\\
		\leq &\, \smfrac 72\cdot\big[\ms g(N,R)+\lambda(\ms h(R),R)\big]
		.
	\end{align*}
\end{proof}

\noindent
The second lemma accounts for the contribution of the remaining intervals. 

\begin{lemma}
\label{B26}
	Let $(l_j)_{j=1}^\infty$ be a summable sequence of positive numbers, and $(\phi_j)_{j=1}^\infty$ a sequence of real
	numbers. Denote by $H$ the Hamburger Hamiltonian with these lengths and angles, and let $W_H$ be its monodromy matrix. 

	Let $\ms c_l,\ms c_\phi\DF[1,\infty)\to(0,\infty)$ be two functions with $\ms c_\phi\leq\ms c_l$. Choose $\psi\in\bb R $
	and assume that 
	\[
		\forall N\in\bb N\DP 
		\sum_{j=N+1}^\infty l_j\leq\ms c_l(N)\ \wedge\ \sum_{j=N+1}^\infty l_j\sin^2(\phi_j-\psi)\leq\ms c_\phi(N)
		.
	\]
	Then, for any $N \in \bb N$ and $z \in \bb C$,
	\[
		\|W(x_N,L;z)\| \leq \Big(\frac{\ms c_l(N)}{\ms c_\phi (N)}\Big)^{\frac 12} 
		\exp \bigg(2|z| \big(\ms c_l \ms c_\phi \big)^{\frac 12}(N) \bigg)
		.
	\]
\end{lemma}
\begin{proof}
	Let $\Omega \in \GL (2,\bb R)$ and consider the function $\tilde{W}(t) \DE \Omega W(x_N,t;z) \Omega^{-1} $ defined on 
	$[x_N,L]$. This function satisfies $\tilde{W}(x_N)=I$ and
	\[
		\tilde{W}'(t)=-z \tilde W(t) \Omega H(t) J \Omega^{-1}, \qquad t \in [x_N,L].
	\]
	By Gr\"onwall's Lemma,
	\[
		\|\tilde{W}(t)\| \leq \exp \bigg(|z| \int_{x_N}^t \|\Omega H(s) J \Omega^{-1}\| \DD s \bigg).
	\]
	Set $a \DE \Big(\frac{\ms c_\phi(N)}{\ms c_l(N)}\Big)^{\frac 14}$ and $\Omega = \Omega (a,\psi )$. 
	Using \Cref{B20}, $(i),(ii)$, we obtain
	\begin{align*}
		\|W(x_N &,L;z)\| \leq
		\|\Omega^{-1}\| \|\tilde{W}(L)\| \|\Omega\| \leq \frac{1}{a^2} \exp \bigg(|z| \mkern-5mu 
		\int_{x_N}^L \|\Omega H(s) J \Omega^{-1}\| \DD s \bigg) 
		\\
		= &\, 
		\Big(\frac{\ms c_l(N)}{\ms c_\phi (N)}\Big)^{\frac 12} 
		\exp \bigg(|z| \sum_{j=N+1}^\infty l_j \cdot \|\Omega \xi_{\phi_j}\xi_{\phi_j}^T J \Omega^{-1}\| \bigg) 
		\\
		= &\, 
		\Big(\frac{\ms c_l(N)}{\ms c_\phi (N)}\Big)^{\frac 12} 
		\exp \bigg(|z| 
		\sum_{j=N+1}^\infty l_j \big[a^2 \cos^2 (\phi_j - \psi) +\frac{1}{a^2} \sin^2 (\phi_j - \psi) \big] \bigg)
		\\
		\leq &\, 
		\Big(\frac{\ms c_l(N)}{\ms c_\phi (N)}\Big)^{\frac 12} 
		\exp \bigg(|z| 
		\bigg[a^2\sum_{j=N+1}^\infty l_j  +\frac{1}{a^2} \sum_{j=N+1}^\infty l_j\sin^2 (\phi_j - \psi) \bigg] \bigg)
		\\
		\leq &\, \Big(\frac{\ms c_l(N)}{\ms c_\phi (N)}\Big)^{\frac 12} 
		\exp \bigg(2|z| \big(\ms c_l \ms c_\phi \big)^{\frac 12}(N) \bigg)
	\end{align*}
\end{proof}

\noindent
The proof of \Cref{B8} is now easily completed. 

\begin{proof}[Proof of \Cref{B8}]
	Let $(t,R) \in [1,\infty ) \times \Big[\frac{2}{(\ms d_l \ms d_\phi)(1)} \Big)$ and set $N \DE \lceil t \rceil$. 
	We use \Cref{B25} and \Cref{B26} to estimate, for $|z|=R$,
	\begin{align*}
		\log\|W_H(z)\| \leq &\, \log \|W(0,x_N;z)\| + \log \|W(x_N,L;z)\| 
		\\
		\leq &\, \smfrac 72\cdot\Big[\ms g(N,R) + \lambda \big(\min \{N,\ms h(R)\},R \big) \Big] 
		\\
		& \, +\Big[\smfrac 12\log \frac{\ms c_l(N)}{\ms c_\phi (N)} + 2R\big(\ms c_l \ms c_\phi \big)^{\frac 12}(N) \Big] 
		\\
		\leq &\, \smfrac 72\cdot\Big[\ms g(t,R)+R\big(\ms c_l \ms c_\phi \big)^{\frac 12}(t)\Big] 
		+\smfrac 72\cdot\big|\ms g(N,R)-\ms g(t,R)\big|
		\\
		&\, +\smfrac 72\cdot\lambda\big(\min\{N,\ms h(R)\},R\big)+\smfrac 12\log \frac{\ms c_l(N)}{\ms c_\phi (N)}
		.
	\end{align*}
	Since $0\leq N-t\leq 1$ and 
	\begin{align*}
		g(s,R)\leq &\, 
		\begin{cases}
			\frac 32\big(\log^+R+\log^+\frac{\ms d_l(1)}{\ms d_\phi(1)}\big) \CAS 1\leq s<\ms k(R),
			\\
			\sqrt 2 \CAS \ms k(R)\leq s<\ms h(R),
			\\
			1 \CAS \ms h(R)\leq s
			,
		\end{cases}
		\\[3mm]
		\leq &\, \frac 32\cdot\lambda\big(\min\{N,\ms h(R)\},R\big)
		,
	\end{align*}
	we have 
	\[
		\big|\ms g(N,R)-\ms g(t,R)\big|\leq \frac 32\cdot\lambda\big(\min\{N,\ms h(R)\},R\big)
		,
	\]
	and obtain 
	\[
		\log\|W_H(z)\|\leq 9\cdot\Big[\max\Big\{\ms g\big(t,R\big),R(\ms c_l\ms c_\phi)^{\frac 12}(t)\Big\}+L(t,R)\Big].
	\]
	Since $t$ was arbitrary, we can pass to the infimum on the right side.
\end{proof}

\section{Properties of the upper bound}

In this section we study the expression on the right side of \cref{B32} independently of its meaning in the context of
\Cref{B8}. 

\begin{Definition}
\label{B53}
	Let $\ms d_l\DF[1,\infty)\to(0,\infty)$ and $\ms d_\phi\DF[1,\infty)\to(0,1]$ be measurable with 
	$\ms d_l\asymp 1\asymp\ms d_\phi$ locally. Then we denote
	\begin{equation}
	\label{B3}
		B(R)\DE
		\inf_{t\geq 1}\Big(\max\Big\{\ms g\big(t,R\big),R(\ms c_l\ms c_\phi)^{\frac 12}(t)\Big\}+L(t,R)\Big)
		.
	\end{equation}
\end{Definition}

\noindent
We consider the questions
\begin{Itemize}
\item Is $L(t,R)$ small compared to $\max\big\{\ms g\big(t,R\big),R(\ms c_l\ms c_\phi)^{\frac 12}(t)\big\}$ ?
\item Does $\inf_{t\geq 1}\max\big\{\ms g\big(t,R\big),R(\ms c_l\ms c_\phi)^{\frac 12}(t)\big\}$ depend monotonically on the 
	data $\ms d_l,\ms d_\phi,\ms c_l,\ms c_\phi$ ?
\end{Itemize}
It turns out that already very weak assumptions on the data, concerning monotonicity, continuity, variation, and power 
boundedness, are sufficient to ensure that the answers are affirmative. 
This tells us that usually $L(t,R)$ can be dropped from \cref{B3}, and hence $B(R)$ can be evaluated much more easily, and that 
usually the natural intuition which stems from the context of \Cref{B8}, that faster decay means smaller growth, is
indeed reflected in $B(R)$. 

\subsection{About monotonicity in the data}

The essence is the following result. 

\begin{proposition}
\label{B47}
	Let $\ms d_l,\widehat{\ms d}_l\DF[1,\infty)\to(0,\infty)$ and 
	$\ms d_\phi,\widehat{\ms d}_\phi\DF[1,\infty)\to(0,1]$ be measurable with 
	$\ms d_l\asymp\ms d_\phi\asymp\widehat{\ms d}_l\asymp\widehat{\ms d}_\phi\asymp 1$ locally, 
	and assume that $\widehat{\ms d}_l\widehat{\ms d}_\phi$ is $\asymp$ to some nonincreasing function.
	If $\widehat{\ms d}_l\lesssim \ms d_l$ and $\widehat{\ms d}_\phi\lesssim \ms d_\phi$ on $[1,\infty)$, then 
	(quantities $\widehat{\ms g}(t,R)$ etc.\ are defined correspondingly for $\widehat{\ms d}_l,\widehat{\ms d}_\phi$) 
	\begin{align}
	\label{B62}
		\widehat g(t,R)\lesssim g(t,R)
	\end{align}
	for $R \geq \frac{2}{(\ms d_l \ms d_\phi)(1)}$ and a.e. $(t,R)$ with $t < \widehat{\ms h}(R)$ or $t>\ms h(R)$. 
	If additionally $\frac{\widehat{\ms d}_\phi}{\widehat{\ms d}_l}$ is bounded or $\asymp$ on $[1,\infty)$ to a 
	nondecreasing function, then \cref{B62} holds for $R$ sufficiently large and a.e. $t \in [1,\infty)$.

	The constant that is implicit in the relation \cref{B62} depends on the constants implicit in the assumptions, 
	but not on the actual functions $\ms d_l,\ms d_\phi,\widehat{\ms d}_l,\widehat{\ms d}_\phi$. 
\end{proposition}

\noindent
We start with a preparatory lemma.

\begin{lemma}
\label{B72}
	Let $\Phi,\Psi\DF[1,\infty)\to(0,\infty)$ with $\Phi$ nondecreasing, and assume that $\alpha,\alpha'>0$ are such that 
	\[
		\forall t\in[1,\infty)\DP 
		\alpha\cdot\Phi(t)\leq\Psi(t)\leq\alpha'\cdot\Phi(t).
	\]
	Set
	\[
		t_0\DE\sup \bigg(\big\{t\in[1,\infty)\DS \sup_{1\leq s\leq t}\Psi(s)\leq 1\big\} \cup \{1\} \bigg) 
		\,\in [1,\infty]
		.
	\]
	Then, if $t_0 <\infty$,
	\[
		\forall t\in(t_0,\infty)\DP
		\Psi(t)>\frac\alpha{\alpha'}
		.
	\]
\end{lemma}
\begin{proof}
	Since $\Phi$ is nondecreasing, we have
	\begin{align*}
		\big\{t\in[1,\infty)\DS \alpha'\Phi(t)\leq 1\big\} \cup \{1\}
		= &\, \big\{t\in[1,\infty)\DS \sup_{1\leq s\leq t} \alpha'\Phi(s)\leq 1\big\} \cup \{1\}
		\\
		\subseteq &\, \big\{t\in[1,\infty)\DS \sup_{1\leq s\leq t} \Psi(s)\leq 1\big\} \cup \{1\}
		.
	\end{align*}
	Set
	\[
		t_1\DE\sup \bigg(\big\{t\in[1,\infty)\DS \alpha'\Phi(t)\leq 1\big\} \cup \{1\} \bigg) \, \in [1,\infty]
		,
	\]
	then $t_1\leq t_0$. Again by monotonicity, $\alpha'\Phi(t)>1$ for all $t>t_1$. and this yields that for $t>t_1$
	\[
		\Psi(t)\geq \alpha \Phi(t)>\frac\alpha{\alpha'}
		.
	\]
\end{proof}

\begin{proof}[Proof of \Cref{B47}]
	We have to trace constants, and thus make the constants from the assumptions explicit:
	choose a nondecreasing function $\ms u$ and $\alpha,\alpha'>0$, such that 
	\[
		\alpha\ms u(t)\leq \frac{2}{(\widehat{\ms d}_l \widehat{\ms d}_\phi)(t)} \leq \alpha'\ms u(t)
		, \qquad t \in [1,\infty)
		,
	\]
	and choose $\kappa_l,\kappa_\phi>0$, such that 
	\[
		\widehat{\ms d}_l(t)\leq\kappa_l\ms d_l(t),\quad
		\widehat{\ms d}_\phi(t)\leq\kappa_\phi \ms d_\phi(t)
		, \qquad t \in [1,\infty)
		.
	\]
	Moreover, set 
	\[
		\lambda\DE\max_{t\geq 2}\frac{\log t}{\sqrt t}
		,
	\]
	and note that 
	\[
		x\leq\gamma y\ \Rightarrow\ \log x\leq\Big(1+\frac{\log^+\gamma}{\log 2}\Big)\cdot\log y
		\quad \text{for}\ x>0,y\geq 2,\gamma>0
		.
	\]
	Next, observe that 
	\begin{equation}
	\label{B30}
		\big[R(\widehat{\ms d}_l \widehat{\ms d}_\phi)(t)\big]^{\frac 12}
		<\Big(\frac{2\alpha'}{\alpha}\Big)^{\frac 12}, \qquad t>\widehat{\ms k}(R).
	\end{equation}
	In order to see this, let $R$ be fixed and apply \Cref{B72} with $\Phi (t)=\frac{\ms u(t)}{R}$ and 
	$\Psi (t) = \frac{2}{R(\widehat{\ms d}_l \widehat{\ms d}_\phi)(t)}$. Then $t_0=\widehat{\ms k}(R)$, and hence 
	for $t>\widehat{\ms k}(R)$ we have $\Psi (t) > \frac{\alpha}{\alpha'}$, which implies \cref{B30}.

	In order to establish that $\widehat g(t,R)\lesssim g(t,R)$ we
	first assume $t < \widehat{\ms h}(R)$ and distinguish the following six cases according to the definitions of 
	$\widehat g(t,R)$ and $ g(t,R)$. 
	\begin{IIlist}
	\item $1\leq t<\widehat{\ms k}(R)\ \wedge\ 1\leq t<\ms k(R)$:
		\[
			\log\big[R(\widehat{\ms d}_l\widehat{\ms d}_\phi)(t)\big]
			\leq \Big(1+\frac{\log^+(\kappa_l\kappa_\phi)}{\log 2}\Big)\log\big[R(\ms d_l\ms d_\phi)(t)\big]
			.
		\]
	\item $1\leq t<\widehat{\ms k}(R)\ \wedge\ \ms k(R)<t<\ms h(R)$:
		\[
			\log\big[R(\widehat{\ms d}_l\widehat{\ms d}_\phi)(t)\big]
			\leq \lambda\big[R(\widehat{\ms d}_l\widehat{\ms d}_\phi)(t)\big]^{\frac 12}
			\leq\lambda(\kappa_l\kappa_\phi)^{\frac 12}\big[R(\ms d_l\ms d_\phi)(t)\big]^{\frac 12}
			.
		\]
	\item $1\leq t<\widehat{\ms k}(R)\ \wedge\ \ms h(R)<t$:
		\begin{align*}
			\log\big[R(\widehat{\ms d}_l\widehat{\ms d}_\phi)(t)\big]
			\leq &\, \lambda\big[R(\widehat{\ms d}_l\widehat{\ms d}_\phi)(t)\big]^{\frac 12}
			\leq\lambda\big[R(\widehat{\ms d}_l\widehat{\ms d}_\phi)(t)\big] 
			\\
			\leq &\, \lambda\big[R\widehat{\ms d}_l(t)\big] 
			\leq \lambda \kappa_l\big[R\ms d_l(t)\big] 
			.
		\end{align*}
	\item $\widehat{\ms k}(R)<t<\widehat{\ms h}(R)\ \wedge\ 1\leq t<\ms k(R)$. By \cref{B30},
		\[
			\big[R(\widehat{\ms d}_l\widehat{\ms d}_\phi)(t)\big]^{\frac 12}
			<\Big(\frac{2\alpha'}{\alpha}\Big)^{\frac 12}
			\leq\frac{1}{\log 2}\Big(\frac{2\alpha'}{\alpha}\Big)^{\frac 12}\log\big[R(\ms d_l\ms d_\phi)(t)\big]
			.
		\]
	\item $\widehat{\ms k}(R)<t<\widehat{\ms h}(R)\ \wedge\ \ms k(R)<t<\ms h(R)$:
		\[
			\big[R(\widehat{\ms d}_l\widehat{\ms d}_\phi)(t)\big]^{\frac 12}
			\leq(\kappa_l\kappa_\phi)^{\frac 12}\big[R(\ms d_l\ms d_\phi)(t)\big]^{\frac 12}
			.
		\]
	\item $\widehat{\ms k}(R)<t<\widehat{\ms h}(R)\ \wedge\ \ms h(R)<t$:
		\[
			\big[R(\widehat{\ms d}_l\widehat{\ms d}_\phi)(t)\big]^{\frac 12}
			\leq R\widehat{\ms d}_l(t)
			\leq \kappa_lR\ms d_l(t)
			.
		\]
	\end{IIlist}
	Since $(t,R)$ is assumed to satisfy $t < \widehat{\ms h}(R) \vee\ \ms h(R) < t$, and $t < \widehat{\ms h}(R)$ was 
	already covered, the only case left to treat is
	\begin{IIlist}
	\item $\widehat{\ms h}(R)<t\ \wedge\ \ms h(R)<t$:
		\[
			R\widehat{\ms d}_l(t)\leq\kappa_lR\ms d_l(t)
			.
		\]
	\end{IIlist}
	Now we assume the additional condition on $\frac{\widehat{\ms d}_\phi}{\widehat{\ms d}_l}$. 
	If this quotient is bounded, then $\widehat{\ms h}(R)=\infty$ for sufficiently large $R$, and there is nothing to show. 
	Otherwise, choose a nondecreasing function $\ms v$ and $\beta,\beta'>0$ such that
	\[
		\beta\ms v(t)\leq \frac{\widehat{\ms d}_\phi}{\widehat{\ms d}_l} \leq \beta'\ms v(t), \qquad t \in [1,\infty)
		.
	\]
	We claim that 
	\begin{align}
	\label{B34}
		R\widehat{\ms d}_l(t)\leq\frac{\beta'}{\beta}\widehat{\ms d}_\phi(t), \qquad t > \widehat{\ms h}(R).
	\end{align}
	In fact, if we fix $R$ and apply \Cref{B72} with $\Phi(t)=\frac{\ms v(t)}{R}$ and 
	$\Psi (t)=\frac{\widehat{\ms d}_\phi(t)}{R\widehat{\ms d}_l(t)}$, we find that $t_0=\widehat{\ms h}(R)$ and 
	$\Psi (t)>\frac{\beta}{\beta'}$ for $t>\widehat{\ms h}(R)$. This is equivalent to \cref{B34}.

	With these preparations, the assertion follows from the following simple inequalities:
	\begin{IIlist}
	\item $\widehat{\ms h}(R)<t\ \wedge\ 1\leq t<\ms k(R)$:
		\[
			R\widehat{\ms d}_l(t)\leq\frac{\beta'}{\beta}\widehat{\ms d}_\phi(t)\leq \frac{\beta'}{\beta}
			\leq\frac{\beta'}{\beta \log 2}\log\big[R(\ms d_l\ms d_\phi)(t)\big]
			.
		\]
	\item $\widehat{\ms h}(R)<t\ \wedge\ \ms k(R)<t<\ms h(R)$:
		\[
			R\widehat{\ms d}_l(t)\leq\frac{\beta'}{\beta}\widehat{\ms d}_\phi(t)
			\leq\frac{\beta'}{\beta}\kappa_\phi\ms d_\phi(t)
			\leq\frac{\beta'}{\beta}\kappa_\phi\big[R(\ms d_l\ms d_\phi)(t)\big]^{\frac 12}
			.
		\]
	\end{IIlist}
\end{proof}

\begin{corollary}
\label{B37}
	Let $\ms d_l,\ms d_\phi,\widehat{\ms d}_l,\widehat{\ms d}_\phi$ be as in \Cref{B47} and subject to the additional
	condition of this proposition. Further, let 
	$\ms c_l,\ms c_\phi,\widehat{\ms c}_l,\widehat{\ms c}_\phi\DF[1,\infty)\to(0,\infty)$.
	If 
	\[
		\widehat{\ms d}_l\lesssim \ms d_l,\ \widehat{\ms d}_\phi\lesssim \ms d_\phi,\ 
		\widehat{\ms c}_l\lesssim \ms c_l,\ \widehat{\ms c}_\phi\lesssim \ms c_\phi
	\]
	on $[1,\infty)$, then 
	\[
		\inf_{t\geq 1}\max\big\{\widehat{\ms g}\big(t,R\big),R(\widehat{\ms c}_l\widehat{\ms c}_\phi)^{\frac 12}(t)\big\}
		\lesssim\inf_{t\geq 1}\max\big\{\ms g\big(t,R\big),R(\ms c_l\ms c_\phi)^{\frac 12}(t)\big\}
	\]
	for sufficiently large $R$. 
	The constant that is implicit in the assertion depends on the constants implicit in the assumptions, 
	but not on the actual functions $\ms d_l,\widehat{\ms d}_l$ etc. 
\end{corollary}

\subsection{About smallness of $L(t,R)$}

To start with, let us observe that $\ms g(t,R)$ must grow at least logarithmically. 

\begin{lemma}
\label{B5} 
	Let $\ms d_l\DF[1,\infty)\to(0,\infty)$ and $\ms d_\phi\DF[1,\infty)\to(0,1]$ be measurable with 
	$\ms d_l\asymp 1\asymp\ms d_\phi$ locally. Then, for each $\varepsilon>0$, there exists a positive constant $c>0$ such 
	that, for $t \geq 1+\varepsilon$ and $R \geq \sup_{1 \leq s\leq 1+\varepsilon} \frac{2}{(\ms d_l \ms d_\phi)(s)}$ 
	we have 
	\[
		\ms g(t,R)\geq c\log R
		.
	\]
\end{lemma}
\begin{proof}
	Let $\eta \DE \sup_{1 \leq s \leq 1+\varepsilon} \frac{2}{(\ms d_l \ms d_\phi)(s)}$. 
	Then $R \geq \eta$ is equivalent to $1+\varepsilon \leq \ms k(R)$. Hence
	\begin{align*}
		\ms g(t,R)\geq &\, \int_1^{1+\varepsilon}\log\big(R(\ms d_l\ms d_\phi)(s)\big)\DD s
		\\
		= &\, \Big(\varepsilon+\frac 1{\log R}\int_1^{1+\varepsilon}\log(\ms d_l\ms d_\phi)(s)\DD s\Big)
		\cdot\log R
		.
	\end{align*}
	The term in the bracket tends to $\varepsilon$ for $R\to\infty$, and hence we find $R_0 \geq \eta$ such
	that $\ms g(t,R)\geq\frac\varepsilon 2\cdot\log R$ for all $t\geq 1+\varepsilon$ and $R\geq R_0$. The function $\ms g$
	is nonzero, and hence bounded away from zero, on the compact set $\{1+\varepsilon\}\times[\eta,R_0]$. By
	monotonicity, it is therefore also bounded away from zero on $[1+\varepsilon,\infty)\times[\eta,R_0]$. Together, we
	see that a constant $c>0$ can be
	chosen as required. 
\end{proof}

\noindent
The most cumbersome part of $L(t,R)$ is the sum written as the last of the six summands in \cref{B49}. 
Assuming monotonicity of the quotient $\frac{\ms d_\phi}{\ms d_l}$, this sum turns into a telescoping sum and can be estimated. 

\begin{lemma}
\label{B68}
	Let $\ms d_l,\ms d_\phi\DF[1,\infty)\to(0,\infty)$ and let $\ms h(R)$ be as in \cref{B45}. Assume that the quotient 
	$\frac{\ms d_\phi}{\ms d_l}$ is eventually monotone but not eventually constant. Then 
	\[
		\sum_{j=1}^{\min\{\lceil t\rceil,\lfloor\ms h(R)\rfloor\}-1} 
		\bigg| \log\bigg( \raisebox{4pt}{$\frac{\ms d_\phi(j)}{\ms d_l(j)}$}\mkern-3mu\bigg/\mkern-5mu
		\raisebox{-5pt}{$\frac{\ms d_\phi(j+1)}{\ms d_l(j+1)}$}\bigg) \bigg|
		\asymp 1+\bigg|\log \frac{\ms d_\phi(\min\{\lceil t\rceil,\lfloor\ms h(R)\rfloor\})}{
		\ms d_l(\min\{\lceil t\rceil,\lfloor\ms h(R)\rfloor\})}\bigg|
	\]
	for $(t,R)\in[1,\infty)\times[\frac{\ms d_\phi(1)}{\ms d_l(1)},\infty)$. 
\end{lemma}
\begin{proof}
	We abbreviate $T\DE\min \{\lceil t \rceil,\lfloor\ms h(R) \rfloor\}$. 
	Choose $M\in\bb N$ such that $\frac{\ms d_\phi}{\ms d_l}$ is monotone on $[M,\infty)$ and 
	$\frac{\ms d_\phi (M)}{\ms d_l(M)}\neq\frac{\ms d_\phi (M+1)}{\ms d_l(M+1)}$.
	For all $(t,R)$ such that $T\geq M+1$ we have
	\begin{align*}
		\sum_{j=M}^{T-1}
		\bigg| \log\bigg( \raisebox{4pt}{$\frac{\ms d_\phi(j)}{\ms d_l(j)}$}\mkern-3mu\bigg/\mkern-5mu
		\raisebox{-5pt}{$\frac{\ms d_\phi(j+1)}{\ms d_l(j+1)}$}\bigg) \bigg|
		= &\, 
		\bigg| \sum_{j=M}^{T-1}\log\bigg( \raisebox{4pt}{$\frac{\ms d_\phi(j)}{\ms d_l(j)}$}\mkern-3mu\bigg/\mkern-5mu
		\raisebox{-5pt}{$\frac{\ms d_\phi(j+1)}{\ms d_l(j+1)}$}\bigg) \bigg|
		\\
		= &\, 
		\bigg|\log \frac{\ms d_\phi (M)}{\ms d_l(M)}- \log \frac{\ms d_\phi(T)}{\ms d_l(T)} \bigg|
		.
	\end{align*}
	Clearly, this value is nondecreasing in $T$ and positive for $T\geq M+1$. 
	If $\big|\log\frac{\ms d_\phi(T)}{\ms d_l(T)}\big|$ is bounded, it is $\asymp 1$. 
	If $\big|\log\frac{\ms d_\phi(T)}{\ms d_l(T)}\big|$ is unbounded, it is 
	$\asymp\big|\log\frac{\ms d_\phi(T)}{\ms d_l(T)}\big|$. 

	The beginning part of the sum which we cut off, i.e., 
	\[
		\sum_{j=1}^{M-1}
		\bigg| \log\bigg( \raisebox{4pt}{$\frac{\ms d_\phi(j)}{\ms d_l(j)}$}\mkern-3mu\bigg/\mkern-5mu
		\raisebox{-5pt}{$\frac{\ms d_\phi(j+1)}{\ms d_l(j+1)}$}\bigg) \bigg|
	\]
	is some nonnegative number independent of $t$ and $R$. 
\end{proof}

\begin{corollary}
\label{B67}
	Let $\ms d_l,\ms d_\phi\DF[1,\infty)\to(0,\infty)$ and let $\ms h(R)$ be as in \cref{B45}.
	Assume that one of the following two assumptions holds.
	\begin{Enumerate}
	\item The quotient $\frac{\ms d_\phi}{\ms d_l}$ is eventually nondecreasing and 
		$\frac{\ms d_\phi}{\ms d_l}(t)\lesssim\sup_{1\leq s<t}\frac{\ms d_\phi}{\ms d_l}(s)$ 
		for sufficiently large $t$.
	\item The quotient $\frac{\ms d_\phi}{\ms d_l}$ is eventually nonincreasing and power bounded from below (i.e.\ 
		there exists $\alpha>0$ such that $\frac{\ms d_\phi(t)}{\ms d_l(t)}\gtrsim t^{-\alpha}$ 
		for sufficiently large $t$).
	\end{Enumerate}
	Then
	\begin{multline*}
		\sum_{j=1}^{ \min \{\lceil t \rceil,\lfloor\ms h(R) \rfloor\}-1 } 
		\bigg| \log\bigg( \raisebox{4pt}{$\frac{\ms d_\phi(j)}{\ms d_l(j)}$}\mkern-3mu\bigg/\mkern-5mu
		\raisebox{-5pt}{$\frac{\ms d_\phi(j+1)}{\ms d_l(j+1)}$}\bigg) \bigg|
		\\
		\lesssim 
		\begin{cases}
			1+\log^+ R\CAS \text{ {\rm(i)} holds}, 
			\\[2mm]
			1+\log t\CAS \text{ {\rm(ii)} holds},
		\end{cases}
	\end{multline*}
	for $(t,R) \in [1,\infty) \times [\frac{\ms d_\phi(1)}{\ms d_l(1)},\infty)$.
\end{corollary}
\begin{proof}
	Again abbreviate $T\DE\min \{\lceil t \rceil,\lfloor\ms h(R) \rfloor\}$, and let $M\in\bb N$ be such that 
	$\frac{\ms d_\phi}{\ms d_l}$ is monotone on $[M,\infty)$.
	Under the assumption {\rm(i)} we have 
	\[
		\frac{\ms d_\phi(T)}{\ms d_l(T)} 
		\leq \frac{\ms d_\phi(\ms h(R))}{\ms d_l(\ms h(R))}\lesssim R
		,
	\]
	and hence 
	\[
		\bigg|\log \frac{\ms d_\phi(T)}{\ms d_l(T)}\bigg|\lesssim 1+\log R
		.
	\]
	Under the assumption {\rm(ii)} we have 
	\[
		t^{-\alpha} \lesssim T^{-\alpha} \lesssim \frac{\ms d_\phi(T)}{\ms d_l(T)} \lesssim 1
		,
	\]
	and hence 
	\[
		\bigg|\log \frac{\ms d_\phi(T)}{\ms d_l(T)}\bigg|\lesssim 1+\log t
		.
	\]
\end{proof}

\noindent 
Another type of condition on the quotient $\frac{\ms d_\phi}{\ms d_l}$ that ensures that the sum 
can be estimated, is that its variation is not too fast.

\begin{lemma}
\label{B84}
	Let $\ms d_l,\ms d_\phi\DF[1,\infty)\to(0,\infty)$ and let $\ms h(R)$ be as in \cref{B45}.
	Assume that the quotient $\frac{\ms d_\phi}{\ms d_l}$ can be represented as 
	\[
		\frac{\ms d_\phi(t)}{\ms d_l(t)}=c(t) \cdot \exp \bigg(\int_1^t \epsilon (u) \frac{\DD u}{u} \bigg)
		, \qquad t \in [1,\infty),
	\]
	where $\epsilon\DF[1,\infty)\to\bb R$ is locally integrable and eventually bounded, and where 
	$c\DF[1,\infty)\to(0,\infty)$ is eventually constant. Then
	\[
		\sum_{j=1}^{ \min \{\lceil t \rceil,\lfloor\ms h(R) \rfloor\}-1 } 
		\bigg| \log\bigg( \raisebox{4pt}{$\frac{\ms d_\phi(j)}{\ms d_l(j)}$}\mkern-3mu\bigg/\mkern-5mu
		\raisebox{-5pt}{$\frac{\ms d_\phi(j+1)}{\ms d_l(j+1)}$}\bigg) \bigg|
		\lesssim 1+\log t
	\]
	for $(t,R) \in [1,\infty) \times [\frac{\ms d_\phi(1)}{\ms d_l(1)},\infty)$.
\end{lemma}
\begin{proof}
	We again abbreviate $T\DE\min \{\lceil t \rceil,\lfloor\ms h(R) \rfloor\}$. 
	For $M\in\bb N$ such that $\epsilon(t)$ is bounded and $c(t)$ is constant on $[M,\infty)$, we can estimate, for all 
	$(t,R)$ such that $T>M$,
	\begin{multline*}
		\sum_{j=M}^{T-1} 
		\bigg| \log\bigg( \raisebox{4pt}{$\frac{\ms d_\phi(j)}{\ms d_l(j)}$}\mkern-3mu\bigg/\mkern-5mu
		\raisebox{-5pt}{$\frac{\ms d_\phi(j+1)}{\ms d_l(j+1)}$}\bigg) \bigg|
		=\sum_{j=M}^{T-1} \bigg|\int_j^{j+1} \mkern-8mu \epsilon (u)\frac{\DD u}u\bigg|
		\\
		\leq \sum_{j=M}^{T-1}\int_j^{j+1} |\epsilon(u)|\frac{\DD u}u
		=\int_M^T |\epsilon(u)|\frac{\DD u}u \lesssim \log t
		.
	\end{multline*}
\end{proof}

\noindent
Now we give three sets of conditions, each of which ensures that $L(t,R)$ can be neglected.

\begin{proposition}
\label{B59}
	Let $\ms d_l\DF[1,\infty)\to(0,\infty)$ and $\ms d_\phi\DF[1,\infty)\to(0,1]$ be measurable with 
	$\ms d_l\asymp 1\asymp\ms d_\phi$ locally, let $\ms c_l,\ms c_\phi$
	be continuous and nonincreasing with $\lim\limits_{t\to\infty}(\ms c_l\ms c_\phi)(t)=0$, and let $T(R)$ be the 
	unique solution of \cref{B57}. 
	Assume that one of the following three sets of assumptions holds.
	\begin{Enumerate}
	\item The quotient $\frac{\ms d_\phi}{\ms d_l}$ is eventually nondecreasing and 
		$\frac{\ms d_\phi}{\ms d_l}(t)\lesssim\sup_{1\leq s<t}\frac{\ms d_\phi}{\ms d_l}(s)$ 
		for sufficiently large $t$. We have $\ms c_\phi(t)\lesssim\ms c_\phi(t+1)$ for sufficiently large $t$.
	\item The quotient $\frac{\ms d_\phi}{\ms d_l}$ is eventually nonincreasing and power bounded from below. 
		There exists $\alpha>0$ such that $(\ms c_l\ms c_\phi)(t)\lesssim t^{-\alpha}$ for sufficiently large $t$.
		The quotient $\frac{\ms c_l}{\ms c_\phi}$ is power bounded from above or $\ms c_\phi(t)\lesssim\ms c_\phi(t+1)$ 
		for sufficiently large $t$.
	\item The quotient $\frac{\ms d_\phi}{\ms d_l}$ can be represented as 
		\[
			\frac{\ms d_\phi(t)}{\ms d_l(t)}=c(t) \cdot \exp \bigg(\int_1^t \epsilon (u) \frac{\DD u}{u} \bigg)
			, \qquad t \in [1,\infty),
		\]
		where $\epsilon\DF[1,\infty)\to\bb R$ is locally integrable and eventually bounded, and where 
		$c\DF[1,\infty)\to(0,\infty)$ is eventually constant. 
		There exists $\alpha>0$ such that $(\ms c_l\ms c_\phi)(t)\lesssim t^{-\alpha}$ for sufficiently large $t$.
		The quotient $\frac{\ms c_l}{\ms c_\phi}$ is power bounded from above or $\ms c_\phi(t)\lesssim\ms c_\phi(t+1)$ 
		for sufficiently large $t$.
	\end{Enumerate}
	Let $B(R)$ be as in \cref{B3}. Then $L(T(R),R)\lesssim 1+\log^+R$, and 
	\[
		B(R)\asymp\ms g(T(R),R)=R(\ms c_l\ms c_\phi)^{\frac 12}(T(R))
		.
	\]
\end{proposition}
\begin{proof}
	Assume {\rm(i)}. The first sentence of (i) is nothing but the assumption (i) in \Cref{B67}, and hence the sum in the 
	definition of $L(t,R)$ is $\lesssim 1+\log^+R$ independently of $t$. The last but one summand in \cref{B49} is 
	$\lesssim 1$ by monotonicity of $\frac{\ms d_\phi}{\ms d_l}$. It remains to estimate the third summand in \cref{B49}, 
	and here we use the second sentence in the present assumption: using $\lim_{R \to \infty} T(R)=\infty$
	and \Cref{B5},
	\begin{equation}
	\label{B76}
	\begin{aligned}
		\log^+\frac{\ms c_l(\lceil T(R) \rceil)}{\ms c_\phi(\lceil T(R) \rceil)}
		\lesssim &\, 1+\log^+\frac{\ms c_l(T(R))}{\ms c_\phi(T(R))}
		\lesssim 1+\log^+ \frac{1}{(\ms c_l\ms c_\phi)(T(R))} 
		\\
		= &\, 1+\log^+\frac{R^2}{\ms g(T(R),R)^2} \lesssim 1+\log^+R.
	\end{aligned}
	\end{equation}
	Assume {\rm(ii)}. Let $\alpha>0$ be as in the second sentence of {\rm(ii)}, then 
	$R(\ms c_l\ms c_\phi)^{\frac 12}(R^{\frac 2\alpha})\lesssim 1\ll\ms g(R^{\frac 2\alpha},R)$, and hence 
	$T(R)\lesssim R^{\frac 2{\alpha}}$. 
	The first sentence of {\rm(ii)} is nothing but assumption {\rm(ii)} in \Cref{B67}. This yields that the sum in the 
	definition of $L(t,R)$ is $\lesssim 1+\log^+t$, and hence the sum in $L(T(R),R)$ is $\lesssim 1+\log^+R$. 
	The same holds for the last but one summand in \cref{B49} by power boundedness of $\frac{\ms d_\phi}{\ms d_l}$, 
	and for the third summand \cref{B49} in case that $\frac{\ms c_l}{\ms c_\phi}$ is power bounded. 
	If $\ms c_\phi(t)\lesssim\ms c_\phi(t+1)$ for sufficiently large $t$, then the third summand in \cref{B49} is 
	$\lesssim 1+\log^+R$ by \cref{B76}.

	Assume {\rm(iii)}. The second sentence again ensures that $T(R)$ is power bounded from above. The first sentence is the
	assumption of \Cref{B84}, and it follows that the sum in $L(T(R),R)$ is $\lesssim 1+\log^+R$. The form of 
	$\frac{\ms d_\phi}{\ms d_l}$ implies that this quotient is power bounded from below, and hence the last but one summand
	\cref{B49} is $\lesssim 1+\log^+R$. Concerning the third summand, just argue as above. 

	To finish the proof it remains to refer to \Cref{B5}.
\end{proof}

\section{Regularly varying decay}

In the preceding section, we studied properties of the function $B(R)$, in particular, we saw that the term $L(t,R)$ is usually
neglectable. Our next goal is to explicitly evaluate 
\[
	\min_{t\geq 1}\Big(\max\big\{\ms g(t,R),R(\ms c_l\ms c_\phi)^{\frac 12}(t)\big\}\Big)
\]
in the situation that $\ms d_l,\ms d_\phi,\ms c_l,\ms c_\phi$ are all regularly varying, cf.\ \Cref{B9}. 
As in the preceding section, this is a pure analysis of functions, and independent of the meaning in the context of \Cref{B8}. 
However, after having completed the proof of \Cref{B9}, we will return to the study of Hamburger Hamiltonians and combine 
\Cref{B9} with \Cref{B8}. This establishes an explicit upper bound for the growth of the monodromy matrix, cf.\ \Cref{B7}. 

The setup for \Cref{B9} is as follows:
\begin{Enumerate}
\item We are given regularly varying functions $\ms d_l,\ms d_\phi,\ms c_l,\ms c_\phi\DF[1,\infty)\to(0,\infty)$ with 
	nonpositive indices, $\ms d_\phi\leq 1$, and $\ms d_l\asymp 1\asymp\ms d_\phi$ locally. 
\item We denote
	\[
		\delta_l\DE-\Ind\ms d_l,\quad \delta_\phi\DE-\Ind\ms d_\phi,\quad \gamma_l\DE-\Ind\ms c_l,\quad
		\gamma_\phi\DE-\Ind\ms c_\phi
		,
	\]
	and further
	\begin{align}
	\label{B92}
		\ms D(t)\DE &\, \frac{1}{(\ms d_l\ms d_\phi)(t)},\quad \delta\DE\Ind\ms D=\delta_l+\delta_\phi
		,
		\\
		\label{B93}
		\ms C(t)\DE &\, \frac{1}{(\ms c_l\ms c_\phi)^{\frac 12}(t)},\quad \gamma\DE\Ind\ms C=\frac{\gamma_l+\gamma_\phi}2
		.
	\end{align}
\item The quantities of interest are
	\begin{equation}
	\label{B56}
		\begin{aligned}
			\ms B(t,R) &\, \DE \max\Big\{\ms g(t,R),\frac R{\ms C(t)}\Big\},
			\\
			\ms B(R) &\, \DE \min_{t\geq 1} \ms B(t,R).
		\end{aligned}
	\end{equation}
\end{Enumerate}
Four fundamentally different cases occur. They depend on the relative size of $\ms D(t)$ vs.\ $\ms C(t)$, and on the absolute 
size of $\ms D(t)$. Moreover, there are a few exceptional cases that have to be ruled out. 

\begin{theorem}
\label{B9}
	Assume we are given data as in {\rm(i)} above, and let notation be as in {\rm(ii)} and {\rm(iii)}.
	Assume further that $\delta>0$, that $\ms C$ is $\asymp$ on $[1,\infty)$ to a nondecreasing function, and that 
	$\lim_{t\to\infty}\ms C(t)=\infty$. 

	Then the following estimates for $\ms B(R)$ hold for sufficiently large $R$. 
	\begin{itemize}
	\item[\fbox{\rm A}] {\rm Case $\ms D(t)\lesssim t\ms C(t)$:}
		\\[1mm]
		Choose $\alpha\geq 4\cdot\sup\limits_{t\geq 1}\frac{\ms D(t)}{t\ms C(t)}$, and set 
		$\ms f(t)\DE t\ms C(t)\log\big[\alpha\frac{t\ms C(t)}{\ms D(t)}\big]$.
		Then an asymptotic inverse of $\ms f$ exists, and with $\tau (R) \DE \ms f^- \big( \frac{1}{\alpha}R \big)$ 
		we have
		\begin{equation}
		\label{B73}
			\ms B(R)\asymp \ms B \big(\tau (R),R \big) \asymp \frac R{\ms C(\ms f^-(R))}
			.
		\end{equation}
		The function on the right side is regularly varying with index $\frac 1{1+\gamma}$.

		If $\delta<1+\gamma$ we have $\ms f(t)\asymp t\ms C(t)\log t$, and $\ms B(R)$ is (up to $\asymp$)
		independent of $\ms d_l$ and $\ms d_\phi$. 
	\item[\fbox{\rm B}] {\rm Case $t\ms C(t)\lesssim\ms D(t)$, $\int\limits_1^\infty\ms D(s)^{-\frac 12}\DD s<\infty$, 
		and if $(\delta_l,\delta_\phi,\gamma)=(1,1,0)$ then $\frac {\ms d_\phi}{\ms d_l}$ is bounded or $\asymp$ to a 
		monotone function:}
		\\[1mm]
		Then
		\begin{equation}
		\label{B74}
			\ms k(R)\lesssim\ms B(R)\lesssim R^{\frac 12}\int\limits_{\ms k(R)}^\infty\ms D(s)^{-\frac 12}\DD s
			.
		\end{equation}
		Both bounds in \cref{B74} are regularly varying with index $\frac 1{\delta}$. 

		If $\delta>2$, then $\asymp$ holds throughout \cref{B74} and $\ms B(R)$ is (up to $\asymp$) independent 
		of $\ms c_l$ and $\ms c_\phi$. 
	\item[\fbox{\rm C}] {\rm Case $\frac 1{\ms d_l(t)}\lesssim t\ms C(t)\lesssim\ms D(t)$,
		$\int\limits_1^\infty\ms D(s)^{-\frac 12}\DD s=\infty$, and $(\delta,\gamma)\neq(2,0)$:}
		\\[1mm]
		Set $\ms f_0(t)\DE\frac{t^2\ms C(t)^2}{\ms D(t)}$ and 
		$\ms f_1(t)\DE\Big(\ms C(t)\int\limits_1^t\ms D(s)^{-\frac 12}\DD s\Big)^2$.
		Then asymptotic inverses $\ms f_0^-,\ms f_1^-$ exist, and with
		\[
			\tau (R)\DE\min\big\{\max\big\{\ms f_1^-(R),\ms k(R)\big\},\ms h(R)\big\}
		\]
		we have
		\begin{equation}
		\label{B75}
			\frac R{\ms C(\ms f_0^-(R))}\lesssim\ms B(R)\lesssim\ms B(\tau(R),R)\lesssim\frac R{\ms C(\ms f_1^-(R))}
			.
		\end{equation}
		Both bounds in \cref{B75} are regularly varying with index $\frac{2-\delta+\gamma}{2-\delta+2\gamma}$. 

		If $\delta<2$, then $\asymp$ holds throughout \cref{B75}.
	\item[\fbox{\rm D}] {\rm Case $t\ms C(t)\lesssim\frac 1{\ms d_l(t)}$, 
		$\int\limits_1^\infty\ms D(s)^{-\frac 12}\DD s=\infty$, and $\int_1^\infty\ms d_l(s)\DD s<\infty$:}
		\\[1mm]
		Then $\frac{\ms d_\phi}{\ms d_l}$ is unbounded and 
		\begin{equation}
		\label{B86}
			R\ms h(R)\ms d_l(\ms h(R))\lesssim\ms B(R)\lesssim 
			R^{\frac 12}\int\limits_1^{\ms h(R)}\ms D(s)^{-\frac 12}\DD s+
			R\int\limits_{\ms h(R)}^\infty\ms d_l(s)\DD s
			.
		\end{equation}
		If $\delta_l>\delta_\phi$, then both bounds in \cref{B86} are regularly varying with index 
		$\frac{1-\delta_\phi}{\delta_l-\delta_\phi}$. 

		If $\delta<2$ and $\delta_l>1$, then $\asymp$ holds throughout \cref{B86}.
	\end{itemize}
\end{theorem}

\subsection{Proof of \Cref{B9}}

We start with two lemmata that determine the function $\ms g(t,R)$ on the intervals $[1,\ms k(R)]$ and $[\ms k(R),\ms h(R)]$,
respectively. 

\begin{lemma}
\label{B77}
	Let $\ms D\DF[1,\infty)\to(0,\infty)$ be regularly varying with $\Ind\ms D>0$ and $\ms D\asymp 1$ locally. Set
	\begin{equation}
	\label{B80}
		\ms k(R)\DE
		\sup\Big\{t\in[1,\infty)\DS \sup_{1\leq s\leq t}\frac{2\ms D(s)}R\leq 1\Big\}, \qquad R\geq 2\ms D(1)
		.
	\end{equation}
	Then we have 
	\begin{equation}
	\label{B85}
		\int_1^t\log\frac R{\ms D(s)}\DD s\asymp t\log\frac R{\ms D(t)}
	\end{equation}
	for $R \in [2\ms D(1),\infty)$ and $t\in\big[2\log\frac R{\ms D(1)},\ms k(R)\big]$. 
\end{lemma}
\begin{proof}
	In the first part of the proof we establish the assertion under the additional assumption that $\ms D$ is 
	continuously differentiable with $\ms D'(t)>0$ for all $t\geq 1$. 
	Let $R \in [2\ms D(1),\infty)$ and $t\in\big[2\log\frac R{\ms D(1)},\ms k(R)\big]$. Integration by parts yields that 
	\begin{equation}
	\label{B81}
		\int_1^t\log\frac R{\ms D(s)}\DD s
		=t\log\frac R{\ms D(t)}-\log\frac R{\ms D(1)}+\int_1^t s\frac{\ms D'(s)}{\ms D(s)}\DD s
		.
	\end{equation}
	The relation ``$\gtrsim$'' in \cref{B85} readily follows: The integral on the right side is positive, and by the definition of $\ms k(R)$ we have $\frac R{\ms D(t)}\geq 2$. For $t\geq 2\log\frac R{\ms D(1)}$ it follows
	that 
	\begin{align*}
		\int_1^t\log\frac R{\ms D(s)}\DD s\geq &\, 
		\Big(\frac 12\cdot t\log\frac R{\ms D(t)}+\frac 12\cdot 2\log\frac R{\ms D(1)}\cdot\log 2\Big)
		-\log\frac R{\ms D(1)}
		\\
		\geq &\, \frac 12\cdot t\log\frac R{\ms D(t)}
		.
	\end{align*}
	To show that ``$\lesssim$'' holds in \cref{B85} we make a change of variable and use Karamata's theorem. 
	The inverse function $\ms D^{-1}$ exists and is continuously differentiable with positive derivative, and is
	regularly varying with positive index. We obtain that, for $v\to\infty$, 
	\[
		\int_1^v s\frac{\ms D'(s)}{\ms D(s)}\DD s=\int_{\ms D(1)}^{\ms D(v)}\ms D^{-1}(u)\frac{\DD u}u
		\sim(\Ind\ms D)\cdot v
		.
	\]
	As a function of $v \in [2\log 2,\infty)$, both $v$ and the integral above are continuous and nonzero. Hence
	\[
	\int_1^v s\frac{\ms D'(s)}{\ms D(s)}\DD s \asymp v, \qquad v \in [2\log 2,\infty).
	\]	
	Referring to \cref{B81} and again using that $\frac R{\ms D(t)}\geq 2$, we obtain 
	\[
		\int_1^t\log\frac R{\ms D(s)}\DD s\leq
		t\log\frac R{\ms D(t)}+\int_1^t s\frac{\ms D'(s)}{\ms D(s)}\DD s
		\lesssim t\log\frac R{\ms D(t)}
		.
	\]
	The second part of the proof is to reduce to the above settled special case. 
	Hence, assume that $\ms D$ is given as in the statement of the lemma. \Cref{B31} gives a regularly varying function $\ms D_1$ which is continuously differentiable with $\ms D_1'>0$, such that $\ms D \asymp \ms D_1$ on $[1,\infty)$. Dividing $\ms D_1$ by a sufficiently large constant, we obtain a function $\tilde{\ms D}$ with the properties listed above and 
	\begin{align*}
		\frac 1\alpha\ms D(t)\leq\tilde{\ms D}(t)\leq\ms D(t)\quad\text{for }t\geq 1
	\end{align*}
	for some $\alpha\geq 1$. 

Clearly, $\ms k(R)\leq\tilde{\ms k}(R)$ (where $\tilde{\ms k}(R)$ is defined by \cref{B80} with
	$\tilde{\ms D}$ in place of $\ms D$). It follows that for $t < \ms k(R)$
	\[
		\log\frac R{\ms D(t)}\leq\log\frac R{\tilde{\ms D}(t)}\leq\log\Big(\alpha\frac R{\ms D(t)}\Big)
		\leq \log\frac R{\ms D(t)}+\log\alpha
		\leq\Big(1+\frac{\log\alpha}{\log 2}\Big)\log\frac R{\ms D(t)}
		.
	\]
	We see that 
	\[
		\int_1^t\log\frac R{\ms D(s)}\DD s\asymp\int_1^t\log\frac R{\tilde{\ms D}(s)}\DD s
		\asymp t\log\frac R{\tilde{\ms D(t)}}\asymp t\log\frac R{\ms D(t)}
		.
	\]
\end{proof}

\begin{lemma}
\label{B78}
	Let $\ms D\DF[1,\infty)\to(0,\infty)$ be regularly varying with $\Ind\ms D>0$ and $\ms D\asymp 1$ locally. Then we have 
	\[
		R^{\frac 12}t\ms D(t)^{-\frac 12}\lesssim\ms k(R)+R^{\frac 12}\int_{\ms k(R)}^t\ms D(s)^{-\frac 12}(s)\DD s
		\lesssim R^{\frac 12}\int_1^t\ms D(s)^{-\frac 12}\DD s
	\]
	for $R\in[2\ms D(1),\infty)$ and $t \in [\ms k(R),\infty)$. 
\end{lemma}
\begin{proof}
	The asserted estimate from above is easy to see: by Karamata's theorem we have 
	\[
		\ms k(R)\asymp R^{\frac 12}\ms k(R)\ms D(\ms k(R))^{-\frac 12}\lesssim
		R^{\frac 12}\int_1^{\ms k(R)}\ms D(s)^{-\frac 12}\DD s
		.
	\]
	To obtain an estimate from below, we use \Cref{B31}. This gives a regularly varying function $\ms D_1$ which is 
	continuously differentiable with $\ms D_1'>0$ and satisfies $\ms D \asymp \ms D_1$. By multiplying $\ms D_1$ with a sufficiently large constant, we obtain a function $\tilde{\ms D}$ with the properties listed above and
	\[
		\ms D(t)\leq\tilde{\ms D}(t)\leq\alpha\ms D(t)
	\]
	with some $\alpha\geq 1$. 

	Consider the functions 
	\[
		\ms f_1(t)\DE R^{\frac 12}t\tilde{\ms D}(t)^{-\frac 12},\quad
		\ms f_2(t)\DE \sqrt 2\,\ms k(R)+R^{\frac 12}\int_{\ms k(R)}^t\tilde{\ms D}(s)^{-\frac 12}\DD s
		.
	\]
	Then 
	\begin{multline*}
		\ms f_1(\ms k(R))
		=R^{\frac 12}\ms k(R)\tilde{\ms D}(\ms k(R))^{-\frac 12}
		\\
		\leq R^{\frac 12}\ms k(R)\ms D(\ms k(R))^{-\frac 12}
		=\sqrt 2\,\ms k(R)=\ms f_2(\ms k(R))
		,
	\end{multline*}
	and, for $t>\ms k(R)$, 
	\[
		\ms f_1'(t)=R^{\frac 12}\tilde{\ms D}(t)^{-\frac 12}
		+R^{\frac 12}t\underbrace{\frac{\DD}{\DD t}\Big(\tilde{\ms D}(t)^{-\frac 12}\Big)}_{<0}
		\leq R^{\frac 12}\tilde{\ms D}(t)^{-\frac 12}=\ms f_2'(t)
		.
	\]
	We see that 
	\[
		\ms f_1(t)\leq\ms f_2(t)\quad\text{for }t\geq\ms k(R)
		.
	\]
	It remains to note that 
	\[
		\ms f_1(t)\gtrsim R^{\frac 12}t\ms D(t)^{-\frac 12},\quad 
		\ms f_2(t)\lesssim \ms k(R)+R^{\frac 12}\int_{\ms k(R)}^t\ms D(s)^{-\frac 12}\DD s
		.
	\]
\end{proof}

\begin{lemma}
\label{B51}
	Consider regularly varying functions $\ms d_l,\ms d_\phi, \ms D$ as introduced in the discussion preceding \Cref{B9}.
	Suppose that $\int_1^\infty \ms D(s)^{-\frac 12} \DD s <\infty$. If $\delta_\phi=\delta_l=1$, we also assume that
	$\frac{\ms d_\phi}{\ms d_l}$ is bounded or $\asymp$ to a nondecreasing function. 
	Then for sufficiently large $R$ we have, independently of $t$, that 
	\[
		\ms g(t,R) \lesssim R^{\frac 12}\int\limits_{\ms k(R)}^\infty\ms D(s)^{-\frac 12}\DD s
		.
	\]
\end{lemma}
\begin{proof}
	We distinguish cases. 
	\begin{Ilist}
	\item Assume that $\frac{\ms d_\phi}{\ms d_l}$ is bounded, i.e., $\ms h(R)=\infty$ for all sufficiently large $R$. 
		This is surely the case if $\delta_l<\delta_\phi$.
		
		By Karamata's theorem 
		\[
			\ms k(R)=\sqrt 2\,R^{\frac 12}\ms k(R)\ms D(\ms k(R))^{-\frac 12}\lesssim
			R^{\frac 12}\int_{\ms k(R)}^\infty\ms D(s)^{-\frac 12}\DD s
			.
		\]
		Now it follows that, for all $t\geq\ms k(R)$,
		\[
			\ms g(t,R)\asymp \ms k(R)+R^{\frac 12}\int\limits_{\ms k(R)}^t \ms D(s)^{-\frac 12}\DD s \lesssim 
			R^{\frac 12}\int\limits_{\ms k(R)}^\infty \ms D(s)^{-\frac 12}\DD s
			.
		\]
	\item Assume that $\frac{\ms d_\phi}{\ms d_l}$ is $\asymp$ to a nondecreasing function. 
		This is surely the case if $\delta_l>\delta_\phi$.
		
		If $R$ is large enough and $t > \ms h(R)$,
		\[
			R\ms d_l(t)=\bigg(\frac{R\ms d_l(t)}{\ms d_\phi(t)}\bigg)^{\frac 12} \cdot 
			\big[ R (\ms d_l \ms d_\phi)(t) \big]^{\frac 12} \lesssim 
			\big[ R (\ms d_l \ms d_\phi)(t) \big]^{\frac 12}.
		\]
		We see that also in this case,
		\begin{align*}
			\ms g(t,R)\asymp &\, \ms k(R)
			+R^{\frac 12}\mkern-20mu\int\limits_{\ms k(R)}^{\min\{t,\ms h(R)\}}\mkern-20mu\ms D(s)^{-\frac 12}\DD s
			+R\mkern-20mu\int\limits_{\min\{t,\ms h(R)\}}^t\mkern-20mu\ms d_l(s)\DD s
			\\
			\lesssim &\, R^{\frac 12}\int\limits_{\ms k(R)}^\infty \ms D(s)^{-\frac 12}\DD s.
		\end{align*}
	\item Assume that $\delta_l=\delta_\phi>1$.

		By Karamata's theorem,
		\[
			R\int_{\ms h(R)}^\infty \ms d_l(s) \DD s  \asymp R \ms h(R) \ms d_l(\ms h(R)) \asymp 
			\ms h(R) \ms d_\phi (\ms h(R)).
		\]
		Since $\delta_\phi>1$, the function $t \mapsto t \ms d_\phi (t)$ is eventually decreasing. 
		Since it is locally bounded, it is even bounded on $[1,\infty)$. Again we see that 
		$\ms g(t,R) \lesssim R^{\frac 12}\int\limits_{\ms k(R)}^\infty \ms D(s)^{-\frac 12}\DD s$ independently of $t$.
	\end{Ilist}
\end{proof}

\begin{lemma}
\label{B82}
	In the situation of \Cref{B9}, the inequality
	\begin{equation}
	\label{B94}
		\min \Big\{\ms g(s,R),\frac R{\ms C(s)}\Big\} \lesssim \ms B(R) \leq \ms B(t,R)= 
		\max \Big\{\ms g(t,R),\frac R{\ms C(t)}\Big\} 
	\end{equation}
	holds for $(s,t,R)\in[1,\infty)^2\times[\frac 2{(\ms d_l\ms d_\phi)(1)},\infty)$. 
	In particular, if $\tau (R)$ is a function satisfying $\ms g(\tau (R),R) \asymp \frac{R}{\ms C(\tau (R))}$, then
	\[
		\ms B(R) \asymp \ms g(\tau (R),R) \asymp \frac{R}{\ms C(\tau (R))}.
	\]
\end{lemma}
\begin{proof}
	The upper bound in \cref{B94} is obvious from the definition of $\ms B(t,R)$ and $\ms B(R)$. 
	For the lower bound, we observe that
	\begin{align*}
		s & \leq t \quad \Rightarrow  \quad \ms g(s,R)\leq \ms g(t,R), 
		\\
		s & >t  \quad \Rightarrow \quad \frac{R}{\ms C(s)}\lesssim \frac{R}{\ms C(t)}
	\end{align*}
	and hence
	\[
		\min \Big\{\ms g(s,R),\frac R{\ms C(s)}\Big\} \lesssim \ms B(t,R) \qquad (s,t) \in [1,\infty)^2.
	\]
	Taking the infimum over $t$ proves the lower bound in \cref{B94}.
\end{proof}

\noindent
Now we are ready for the proof of \Cref{B9}.

\begin{proof}[Proof of Case \fbox{A}]
	We have $\Ind\ms f=1+\gamma>0$, and hence $\ms f^-$ exists. Set 
	\[
		\tau (R)\DE\ms f^-\big(\frac 1\alpha R\big)
		.
	\]
	Then, for $R\to\infty$, 
	\[
		\frac 1\alpha R\sim\ms f(\tau (R))=\tau (R)\ms C(\tau (R))
		\log\Big[\alpha\frac{\tau (R)\ms C(\tau (R))}{\ms D(\tau (R))}\Big]
		,
	\]
	and hence 
	\[
		\frac R{\ms D(\tau (R))}\sim
		\Big(\alpha\frac{\tau (R)\ms C(\tau (R))}{\ms D(\tau (R))}\Big)
		\log\Big[\alpha\frac{\tau (R)\ms C(\tau (R))}{\ms D(\tau (R))}\Big]
		.
	\]
	Our choice of $\alpha$ implies that $2\ms D(\tau (R))\leq R$ for all sufficiently large values of $R$, and therefore that 
	$\tau (R)\leq\ms k(R)$ for such $R$. 

	It holds that 
	\[
		\frac t{\log t}\circ t\log t\asymp t\quad\text{for }t\geq 4
		,
	\]
	and we further obtain 
	\[
		\frac{\frac R{\ms D(\tau (R))}}{\log\big(\frac R{\ms D(\tau (R))}\big)}\asymp
		\frac{\tau (R)\ms C(\tau (R))}{\ms D(\tau (R))}
		.
	\]
	Using this and \Cref{B77}, thus
	\[
		\frac R{\ms C(\tau (R))}\asymp \tau (R)\log\Big(\frac R{\ms D(\tau (R))}\Big)\asymp\ms g(\tau (R),R)
		.
	\]
	In view of \Cref{B82}, this completes the proof of \cref{B73}. 

	To see the additional statements, note first that 
	\[
		\Ind\frac R{\ms C(\ms f^-(R))}=1-\gamma\cdot\frac 1{1+\gamma}=\frac 1{1+\gamma}
		.
	\]
	Further $\delta<1+\gamma$ implies $\Ind\frac{t\ms C(t)}{\ms D(t)}>0$, and thus
	$\log\big(\alpha\frac{t\ms C(t)}{\ms D(t)}\big)\asymp\log t$. 
\end{proof}

\begin{proof}[Proof of Case \fbox{B}]
	To obtain the bound from below, we use \Cref{B77} 
	and our assumption that $t\ms C(t)\lesssim\ms D(t)$. Setting $t =\ms k(R) \asymp \ms D^-(R)$, we find that
	\[
		\ms g(\ms k(R),R)\asymp \ms k(R),\qquad 
		\frac R{\ms C(\ms k(R))}\asymp \frac R{\ms C(\ms D^-(R))}\gtrsim\ms D^-(R) \asymp \ms k(R)
		.
	\]
	Now we prove the upper bound for $\ms B(R)$. If \Cref{B51} is applicable, then
	\[
		\ms g(t,R) \lesssim R^{\frac 12}\int\limits_{\ms k(R)}^\infty\ms D(s)^{-\frac 12}\DD s
		.
	\]
	Choosing $\tau(R)$ so large that $\frac{R}{\ms C(\tau (R))}$ is less than the right hand side of this relation, 
	we see that
	\[
		\ms B(R) \leq \ms B(\tau (R),R) \lesssim R^{\frac 12}\int\limits_{\ms k(R)}^\infty\ms D(s)^{-\frac 12}\DD s
		.
	\]
	If \Cref{B51} is not applicable, we must have $\delta_l=\delta_\phi=1$ and $\gamma>0$. Choose $\rho$ so large that 
	\[
		\Ind \ms k=\frac 1\delta >1-\gamma \rho= \Ind \bigg(\frac{R}{\ms C (R^\rho)} \bigg)
		.
	\]
	Since $\frac{\ms d_\phi}{\ms d_l}$ is slowly varying, we have $R^\rho \lesssim \ms h(R)$, cf. \Cref{B42}. Therefore,
	\[
		\ms g(\ms h(R),R) \gtrsim \ms k(R) \gg \frac{R}{\ms C \big(R^{\rho} \big)} \geq  \frac{R}{\ms C (\ms h(R))}
		.
	\]
	We thus obtain
	\[
		\ms B(R) \lesssim \ms B(\ms h(R),R) \asymp \ms g(\ms h(R),R) 
		\lesssim R^{\frac 12}\int\limits_{\ms k(R)}^\infty\ms D(s)^{-\frac 12}\DD s
		.
	\]
	To see the additional statements, note first that $\Ind\ms k=\frac 1\delta$ and also 
	\[
		\Ind\Big(R^{\frac 12}\int_{\ms k(R)}^\infty\ms D(s)^{-\frac 12}\DD s\Big)=
		\frac 12+\frac 1\delta\cdot\Big(1-\frac\delta 2\Big)=\frac 1\delta
		.
	\]
	Further, $\delta>2$ implies that 
	\[
		R^{\frac 12}\int_{\ms k(R)}^\infty\ms D(s)^{-\frac 12}\DD s
		\asymp R^{\frac 12}\ms k(R)\ms D(\ms k(R))^{-\frac 12}\asymp \ms k(R)
		.
	\]
\end{proof}

\begin{proof}[Proof of Case \fbox{C}]
	We have 
	\[
		\Ind\ms f_0=2-\delta+2\gamma,\quad
		\Ind\ms f_1=2\Big(\gamma+\big(1-\frac\delta 2\big)\Big)=2-\delta+2\gamma
		.
	\]
	Since $\ms D^{-\frac 12}$ is not integrable, we must have $\delta\leq 2$. The case that $\delta=2$ and $\gamma=0$ is
	ruled out by assumption, and it follows that $2-\delta+2\gamma>0$. Thus, asymptotic inverses $\ms f_0^-$ and 
	$\ms f_1^-$ exist. 

	To prove the bound from below, we use $\ms f_0^-(R)$ pushed into the interval $[\ms k(R),\ms h(R)]$, namely we set 
	\[
		\tau (R)\DE\min\big\{\max\big\{\ms f_0^-(R),\ms k(R)\big\},\ms h(R)\big\}
		.
	\]
	Our assumption that $\frac 1{\ms d_l(t)}\lesssim t\ms C(t)\lesssim\ms D(t)$ gives 
	\[
		\frac{\ms d_\phi(t)}{\ms d_l(t)}\lesssim\ms f_0(t)\lesssim\ms D(t)
		,
	\]
	and hence $\ms k(R)\lesssim\ms f_0^-(R)\lesssim\ms h(R)$. Thus, $\tau (R)\asymp\ms f_0^-(R)$. 

	Using \Cref{B78} we obtain 
	\[
		\ms g(\tau (R),R)\gtrsim R^{\frac 12}\tau (R)\ms D(\tau (R))^{-\frac 12}
		\asymp R^{\frac 12}\ms f_0^-(R)\ms D(\ms f_0^-(R))^{-\frac 12}
		.
	\]
	On the other hand, the definition of $\ms f_0$ yields for $R \to \infty$ that
	\[
		R\sim\ms f_0(\ms f_0^-(R))=\frac{\ms f_0^-(R)^2\ms C(\ms f_0^-(R))^2}{\ms D(\ms f_0^-(R))}
		,
	\]
	and hence 
	\[
		\frac{R^{\frac 12}}{\ms C(\ms f_0^-(R))}\sim \ms f_0^-(R)\ms D(\ms f_0^-(R))^{-\frac 12}
		.
	\]
	Together, thus, 
	\[
		\min\Big\{\ms g(\tau (R),R),\frac R{\ms C(\tau (R))}\Big\}\gtrsim
		R^{\frac 12}\ms f_0^-(R)\ms D(\ms f_0^-(R))^{-\frac 12}\asymp \frac R{\ms C(\ms f_0^-(R))}
		.
	\]
	For the proof of the upper bound we use 
	\[
		\tau (R)\DE\min\big\{\max\big\{\ms f_1^-(R),\ms k(R)\big\},\ms h(R)\big\}
		.
	\]
	Since $\ms f_0\lesssim\ms f_1$ we have $\ms f_1^-\lesssim\ms f_0^-\lesssim\ms h$, and therefore 
	$\ms f_1^-(R)\lesssim \tau (R)$. This readily implies that 
	\[
		\frac R{\ms C(\tau (R))}\lesssim\frac R{\ms C(\ms f_1^-(R))}
		.
	\]
	To estimate $\ms g(\tau (R),R)$ we distinguish two cases. First assume that $\ms f_1^-(R)\leq\ms k(R)$. Then 
	$\tau (R)=\ms k(R)$, and hence 
	\[
		\ms g(\tau (R),R)\asymp\ms k(R)\lesssim\frac R{\ms C(\ms k(R))}\leq\frac R{\ms C(\ms f_1^-(R))}
		.
	\]
	Second assume that $\ms f_1^-(R)\geq\ms k(R)$. Then $\tau (R)\leq\ms f_1^-(R)$. We use \Cref{B78} to obtain 
	\[
		\ms g(\tau (R),R)\lesssim R^{\frac 12}\int_1^{\tau (R)}\ms D(s)^{-\frac 12}\DD s
		\lesssim R^{\frac 12}\int_1^{\ms f_1^-(R)}\ms D(s)^{-\frac 12}\DD s
		.
	\]
	The last integral can be evaluated from the definition of $\ms f_1$: for $R \to \infty$
	\[
		R\sim\ms f_1(\ms f_1^-(R))=\Big(\ms C(\ms f_1^-(R))\int_1^{\ms f_1^-(R)}\ms D(s)^{-\frac 12}\DD s\Big)^2
		,
	\]
	and hence 
	\[
		\int_1^{\ms f_1^-(R)}\ms D(s)^{-\frac 12}\DD s\sim\frac{R^{\frac 12}}{\ms C(\ms f_1^-(R))}
		.
	\]
	Altogether we see that always 
	\[
		\max\Big\{\ms g(\tau (R),R),\frac R{\ms C(\tau (R))}\Big\}\lesssim\frac R{\ms C(\ms f_1^-(R))}
		.
	\]
	This completes the proof of \cref{B75}. 

	To see the additional statements, recall that $\Ind\ms f_0=\Ind\ms f_1=2-\delta+2\gamma$. This implies that 
	\[
		\Ind\frac R{\ms C(\ms f_0^-(R))}=\Ind\frac R{\ms C(\ms f_1^-(R))}=
		1-\gamma\cdot\frac 1{2-\delta+2\gamma}=\frac{2-\delta+\gamma}{2-\delta+2\gamma}
		.
	\]
	Further, $\delta<2$ implies that $\int_1^t\ms D(s)^{-\frac 12}\DD s\asymp t\ms D(t)^{-\frac 12}$, 
	and hence $\ms f_0\asymp\ms f_1$. 
\end{proof}

\begin{proof}[Proof of Case \fbox{D}]
	The first thing to show is that $\frac{\ms d_\phi}{\ms d_l}$ is unbounded: if we had $\ms d_\phi(t)\lesssim\ms d_l(t)$, 
	it would follow that 
	\[
		\ms D(t)^{-\frac 12}=\sqrt{\ms d_l(t)\ms d_\phi(t)}\lesssim\ms d_l(t)
		,
	\]
	and this contradicts the fact that $\ms d_l(t)$ is integrable while $\ms D^{-\frac 12}(t)$ is not. 

	To show the bound from below asserted in \cref{B86}, we note that 
	\[
		\ms g(\ms h(R),R)\gtrsim R^{\frac 12}\ms h(R)\ms D(\ms h(R))^{-\frac 12}
	\]
	by \Cref{B78}, and that 
	\[
		\frac R{\ms C(\ms h(R))}\gtrsim R\ms h(R)\ms d_l(\ms h(R))
	\]
	by our assumption that $t\ms C(t)\lesssim\frac 1{\ms d_l(t)}$. Since $R\asymp\frac{\ms d_\phi}{\ms d_l}(\ms h(R))$ for 
	$R \to \infty$, we obtain 
	\[
		R\ms h(R)\ms d_l(\ms h(R))\asymp R^{\frac 12}\ms h(R)
		\Big(\frac{\ms d_\phi(\ms h(R))}{\ms d_l(\ms h(R))}\Big)^{\frac 12}\ms d_l(\ms h(R))
		=R^{\frac 12}\ms h(R)\ms D(\ms h(R))^{-\frac 12}
		,
	\]
	and see that 
	\[
		\min\Big\{\ms g(\ms h(R),R),\frac R{\ms C(\ms h(R))}\Big\}\gtrsim R\ms h(R)\ms d_l(\ms h(R))
		.
	\]
	For the proof of the bound from above, we show that $\ms g(t,R)$ is bounded independently of $t\geq\ms h(R)$. 
	Using again \Cref{B78}, we find 
	\begin{align*}
		\ms g(t,R)= &\, \ms g(\ms h(R),R)+R\int_{\ms h(R)}^t\ms d_l(s)\DD s
		\\
		\lesssim &\, R^{\frac 12}\int_1^{\ms h(R)}\ms D(s)^{-\frac 12}\DD s
		+R\int_{\ms h(R)}^\infty\ms d_l(s)\DD s
		.
	\end{align*}
	It remains to choose $\tau (R)$ sufficiently large, so that $\frac R{\ms C(\tau (R))}\leq\ms g(\tau (R)),R)$, 
	which is clearly possible. This completes the proof of \cref{B86}.

	To see the additional statements, assume first that $\delta_l>\delta_\phi$. Then $\ms h(R)$ is regularly varying with 
	index $\frac 1{\delta_l-\delta_\phi}$. We obtain 
	\[
		\Ind\Big(R\ms h(R)\ms d_l(\ms h(R))\Big)=1+\frac 1{\delta_l-\delta_\phi}-\frac{\delta_l}{\delta_l-\delta_\phi}
		=\frac{1-\delta_\phi}{\delta_l-\delta_\phi}
		.
	\]
	By Karamata's theorem 
	\[
		\Ind\Big(R^{\frac 12}\int_1^{\ms h(R)}\ms D(s)^{-\frac 12}\DD s\Big)
		=\frac 12+\big(1-\frac\delta 2\big)\frac 1{\delta_l-\delta_\phi}=\frac{1-\delta_\phi}{\delta_l-\delta_\phi}
		,
	\]
	\[
		\Ind\Big(R\int_{\ms h(R)}^\infty\ms d_l(s)\DD s\Big)
		=1+(1-\delta_l)\frac 1{\delta_l-\delta_\phi}=\frac{1-\delta_\phi}{\delta_l-\delta_\phi}
		.
	\]
	Further, again referring to Karamata's theorem, $\delta<2$ implies 
	\[
		R^{\frac 12}\int_1^{\ms h(R)}\ms D(s)^{-\frac 12}\DD s\asymp R^{\frac 12}\ms h(R)\ms D(\ms h(R))^{-\frac 12}
		\asymp R\ms h(R)\ms d_l(\ms h(R))
		,
	\]
	and $\delta_l>1$ implies 
	\[
		R\int_{\ms h(R)}^\infty\ms d_l(s)\DD s\asymp R\ms h(R)\ms d_l(\ms h(R))
		.
	\]
\end{proof}

\subsection{Bound for the monodromy matrix}

We combine \Cref{B8} with \Cref{B9} to obtain a bound for the growth of $W_H$ when the lengths and angle differences of $H$ are
bounded by regularly varying functions. This yields a far reaching generalisation of \Cref{B66}. 

\begin{theorem}
\label{B7}
	Let $(l_j)_{j=1}^\infty$ be a summable sequence of positive numbers, and $(\phi_j)_{j=1}^\infty$ a sequence of real
	numbers. Denote by $H$ the Hamburger Hamiltonian with these lengths and angles, and let $W_H$ be its monodromy matrix. 

	Let $\psi\in\bb R$ and let $\ms d_l,\ms d_\phi,\ms c_l,\ms c_\phi$ be regularly varying functions that are $\approx$ 
	to some nonincreasing functions, such that $\ms d_l \asymp 1 \asymp \ms d_\phi$ locally, 
	$\ms c_\phi(t)\lesssim\ms c_l(t)$ for sufficiently large $t$, $\Ind\ms d_l+\Ind\ms d_\phi<0$, and 
	$\lim_{t \to \infty} \big(\ms c_l \ms c_\phi \big)(t)=0$. Assume that
	\begin{align*}
		l_j &\lesssim\ms d_l(j),& |\sin(\phi_{j+1}-\phi_j)|&\lesssim\ms d_\phi(j),  
		&&j \text{ sufficiently large}
		,
		\\
		\sum_{j=N+1}^\infty l_j&\lesssim\ms c_l(N),& \sum_{j=N+1}^\infty l_j\sin^2(\phi_j-\psi)&\lesssim\ms c_\phi(N), 
		&&N \text{ sufficiently large}
		.
	\end{align*}
	Denote
	\[
		\ms D(t)\DE\frac{1}{(\ms d_l\ms d_\phi)(t)},\ \delta\DE\Ind\ms D,\quad
		\ms C(t)\DE\frac{1}{(\ms c_l\ms c_\phi)^{\frac 12}(t)},\ \gamma\DE\Ind\ms C
		.
	\]
	Then we have the following bounds for $W_H$ and its order $\rho_H$. 
	\\[3mm]
	\scalebox{0.8}{
	\[
		\begin{array}{lll|l|l}
			\multicolumn{3}{l|}{\text{\rm Data satisfies}} 
			& \log\big(\max_{|z|=R}\|W_H(z)\|\big)\text{\rm\ is }\lesssim
			&\rho_H \leq
			\raisebox{-6pt}{\rule{0pt}{17pt}}
			\\
			\hline
			\hline
			\multicolumn{3}{l|}{\ms D(t)\lesssim t\ms C(t)}
			& \frac R{\ms C(\ms f^-(R))}\quad\text{\rm where}
			\raisebox{-9pt}{\rule{0pt}{25pt}}
			\\
			& 
			& 
			& \ms f(t)\DE t\ms C(t)\log\big[\alpha\frac{t\ms C(t)}{\ms D(t)}\big]
			\raisebox{-6pt}{\rule{0pt}{1pt}}
			& \raisebox{-2pt}{\mbox{\LARGE\(\frac{1}{1+\gamma}\)}}
			\\
			& 
			& 
			& \alpha\DE 4\sup_{t\geq 1}\frac{\ms D(t)}{t\ms C(t)}
			\raisebox{-13pt}{\rule{0pt}{1pt}}
			\\
			\hline
			t\ms C(t)\lesssim\ms D(t),
			& \int\limits_1^\infty\ms D(s)^{-\frac 12}\DD s<\infty,
			& 
			& \raisebox{-10pt}{$R^{\frac 12}\int\limits_{\ms k(R)}^\infty\ms D(s)^{-\frac 12}\DD s$}
			\raisebox{-13pt}{\rule{0pt}{30pt}}
			& \raisebox{-12pt}{\mbox{\LARGE\(\frac{1}{\delta}\)}}
			\\[-8pt]
			& \multicolumn{2}{l|}{\parbox{55mm}{\rm ($\gamma>0$ or $\frac{\ms d_\phi}{\ms d_l}$ 
			$\approx$ to nondecreasing)}}
			& 
			\raisebox{-10pt}{\rule{0pt}{1pt}}
			\\
			\cline{2-5}
			& \int\limits_1^\infty\ms D(s)^{-\frac 12}\DD s=\infty,
			& \frac 1{\ms d_l(t)}\lesssim t\ms C(t),
			& \frac R{\ms C(\ms f_1^-(R))}\quad\text{\rm where}
			\raisebox{-13pt}{\rule{0pt}{30pt}}
			& \raisebox{-10pt}{\mbox{\large\(\frac{2-\delta+\gamma}{2-\delta+2\gamma}\)}}
			\\
			& 
			& \parbox{24mm}{\rm $(\delta,\gamma)\neq(2,0)$}
			& \ms f_1(t)\DE[\ms C(t)\int_1^t\ms D(s)^{-\frac 12}\DD s]^2
			\raisebox{-10pt}{\rule{0pt}{1pt}}
			\\
			\cline{3-5}
			& 
			& t\ms C(t)\lesssim\frac 1{\ms d_l(t)},
			& \raisebox{-6pt}{$R^{\frac 12}\int\limits_1^{\ms h(R)}\ms D(s)^{-\frac 12}\DD s$}
			& \raisebox{-6pt}{\mbox{\Large\(\frac{1-\delta_\phi}{\delta_l-\delta_\phi}\)}}
			\raisebox{-12pt}{\rule{0pt}{31pt}}
			\\[-12pt]
			& 
			& \parbox{25.5mm}{\rm $\int\limits_1^\infty\ms d_l(s)\DD s<\infty$, \\
			$\delta_l\!>\!\delta_\phi$}
			& 
			\raisebox{-2pt}{$\quad \,\,\, +R\mkern-5mu\int\limits_{\ms h(R)}^\infty\mkern-5mu\ms d_l(s)\DD s$}
			& 
			\\[7mm]
			\hline
		\end{array}
	\]
	}
	\\[3mm]
	\noindent In each case the relation $\lesssim$ holds for $R$ sufficiently large.
\end{theorem}

\noindent
The main assumptions on the data $\ms d_l,\ms d_\phi,\ms c_l,\ms c_\phi$ put in this theorem are that those functions are
regularly varying and that $\Ind\ms d_l+\Ind\ms d_\phi<0$ (i.e., $\neq 0$). 
Monotonicity assumptions are only a minor restriction: for example they are automatically fulfilled whenever the 
function under consideration is not slowly varying. The assumption that $\ms c_\phi\lesssim\ms c_l$ is no loss in generality,
since replacing $\ms c_\phi$ by $\min\{\ms c_\phi,\ms c_l\}$ does not affect validity of any of the other assumptions. 

\begin{proof}[Proof of \Cref{B7}] 
\phantom{}
\begin{Elist}
\item Observe that neither the assumptions of the theorem nor the case distinction in the assertion of the theorem depends on the
	equivalence class modulo $\approx$ of the functions $\ms d_l,\ms d_\phi,\ms c_l,\ms c_\phi$. Further, using \Cref{B69}, 
	we see that the functions
	written in the second column of the table change only up to $\approx$ when we pass to other data equivalent modulo
	$\approx$ to $\ms d_l,\ms d_\phi,\ms c_l,\ms c_\phi$.

	Hence, it is enough to prove the theorem for suitable modifications $\widehat{\ms d}_l,\widehat{\ms d}_\phi,
	\widehat{\ms c}_l,\widehat{\ms c}_\phi$ of $\ms d_l,\ms d_\phi,\ms c_l,\ms c_\phi$ which differ only up to $\approx$.
\item The task is to define $\widehat{\ms d}_l,\widehat{\ms d}_\phi,\widehat{\ms c}_l,\widehat{\ms c}_\phi$ in such a way 
	that \Cref{B8}, \Cref{B59}, and \Cref{B9} become applicable. To this end we use \Cref{B54} and the freedom of choice of
	nonincreasing smoothenings mentioned in \Cref{B55}.

	We choose $\ms S[\ms c_\phi]$ such that 
	\[
		\forall N\in\bb N\DP \sum_{j=N+1}^\infty l_j\sin^2(\phi_j-\psi)\leq\ms S[\ms c_\phi](N)
		,
	\]
	which is clearly possible by first choosing an arbitrary nonincreasing smoothening of $\ms c_\phi$ and then multiplying 
	it with a sufficiently large positive constant. Next, we choose $\ms S[\ms c_l]$ such that 
	\[
		\ms S[\ms c_\phi]\leq\ms S[\ms c_l],\qquad 
		\forall N\in\bb N\DP \sum_{j=N+1}^\infty l_j\leq\ms S[\ms c_l](N)
		.
	\]
	A suitable modification of $\ms d_\phi$ is found as follows. If $\liminf_{t\to\infty}\ms d_\phi(t)>0$, we 
	choose $\ms S[\ms d_\phi]\DE 1$ (which corresponds to the choices $\epsilon[\ms d_\phi]\DE 0$ and 
	$\kappa[\ms d_\phi]\DE 1$). If $\lim_{t\to\infty}\ms d_\phi(t)=0$, we choose $\ms S[\ms d_\phi]$ such that 
	\[
		\ms S[\ms d_\phi]\leq 1,\qquad 
		\forall j\in\bb N\DP |\sin(\phi_{j+1}-\phi_j)|\leq\ms S[\ms d_\phi](j)
		,
	\]
	which is possible by choosing an arbitrary nonincreasing smoothening of $\ms d_\phi$, then multiplying 
	it with a sufficiently large positive constant, and then cutting it off at $1$. 

	It remains to define a modification of $\ms d_l$. In the generic case that $\gamma>0$ we proceed just the same
	as above and choose $\ms S[\ms c_l]$ such that
	\begin{equation}
	\label{B61}
		\forall j\in\bb N\DP l_j\leq\ms S[\ms d_l](j)
		.
	\end{equation}
	In the boundary case that $\gamma=0$ we make a further case distinction. 
	If we are in the situation of the 1\textsuperscript{st} or the 3\textsuperscript{rd} row of the table in the theorem,
	we do just the same as above. If we are in the situation of the 2\textsuperscript{nd} or the 4\textsuperscript{th} row,
	the additional assumption ensures that $\frac{\ms d_l}{\ms d_\phi}$ is $\approx$ to some nonincreasing function, and
	we choose 
	\[
		\ms S[\ms d_l]\DE\ms S[\ms d_\phi]\cdot\ms S\Big[\frac{\ms d_l}{\ms d_\phi}\Big]
		,
	\]
	where $\ms S\big[\frac{\ms d_l}{\ms d_\phi}\big]$ is sufficiently large to ensure that \cref{B61} holds. 

	Now we set 
	\[
		\widehat{\ms d}_l\DE\ms S[\ms d_l],\ \widehat{\ms d}_\phi\DE\ms S[\ms d_\phi],\ 
		\widehat{\ms c}_l\DE\ms S[\ms c_l],\ \widehat{\ms c}_\phi\DE\ms S[\ms c_\phi]
		.
	\]
\item We apply our previous results with the data $\widehat{\ms d}_l,\widehat{\ms d}_\phi,\widehat{\ms c}_l,
	\widehat{\ms c}_\phi$. Here we denote by $\widehat B(R)$, $\widehat L(t,R)$, $\widehat{\ms B}(R)$, etc.\ the
	correspondingly defined functions. 

	\Cref{B8} implies that 
	\[
		\log\Big(\max_{|z|=R}\|W_H(z)\|\Big)\lesssim\widehat B(R)
		,
	\]
	and we face the task to control $\widehat L(t,R)$. In almost all cases \Cref{B59} takes care of this:
	\begin{Itemize}
	\item \Cref{B59}\,(i) applies if $\delta_\phi>\delta_l$, or if $\gamma=0$ and we are in the situation of the 
		2\textsuperscript{nd} row of the table. 
	\item \Cref{B59}\,(ii) applies if $\gamma>0$ and $\delta_\phi<\delta_l$. 
	\item \Cref{B59}\,(iii) applies if $\gamma>0$ and $\delta_\phi=\delta_l$. 
	\end{Itemize}
	Thus, in all of these cases, $\widehat B(R)\asymp\widehat{\ms B}(R)$, and the bounds asserted in the second column of
	the table follows from \Cref{B9}. 

	It remains to study the situation that $\gamma=0$, $\delta_\phi\leq\delta_l$, and we are in the 
	1\textsuperscript{st} or 3\textsuperscript{rd} row of the table. Applying \Cref{B67}\,(ii) if
	$\delta_\phi<\delta_l$ and \Cref{B84} if $\delta_\phi=\delta_l$, yields 
	\[
		\widehat L(t,R)\lesssim 1+\log t
		.
	\]
	Let $\tau(R)$ be the power bounded function used in \Cref{B9}\,\fbox{\rm A}\,\fbox{\rm C} to estimate 
	$\widehat{\ms B}(R)$. Then 
	\[
		\widehat B(R)\leq\widehat{\ms B}(\tau(R),R)+\widehat L(\tau(R),R)
		\lesssim\widehat{\ms B}(\tau(R),R)
		.
	\]
	The bounds asserted in the second column of the table thus follow from \Cref{B9}. 
\item The bounds for $\rho_H$ arise simply by taking the indices of the regularly varying functions in the second column. 
\end{Elist}
\end{proof}

\noindent
In the following corollary we revisit the setting of \Cref{B7}, except we are not given functions $\ms c_l,\ms c_\phi$.

\begin{corollary}
\label{B96}
	Let $H$ be a Hamburger Hamiltonian and let $\ms d_l,\ms d_\phi$ be regularly varying. Assume that $\ms d_\phi(t)$ is
	$\sim$ to a nonincreasing function as $t \to \infty$ and that 
	$\ms d_l\asymp 1 \asymp \ms d_\phi$ locally.

	Assume that the lenghts and angles of $H$ are bounded as 
	\[
		 l_j\lesssim\ms d_l(j),\ |\sin(\phi_{j+1}-\phi_j)|\lesssim\ms d_\phi(j) \text{ for sufficiently large } j
		 ,
	\]
	and that $\ms d_l \in L^1([1,\infty))$. Then the following statements hold.
	\begin{Itemize}
	\item If $\delta>2$, then $\log\big(\max_{|z|=R}\|W_H(z)\|\big)\lesssim\ms k(R)$. 
	\item If $0<\delta<2$, then $\log\big(\max_{|z|=R}\|W_H(z)\|\big)\lesssim 
		R\int\limits_{\ms h(R)}^\infty \ms d_l (x) \DD x$.
	\item If $0<\delta\leq 2$, $\delta_l>1$, and there exists $\psi \in \bb R$ such that 
		$|\sin(\phi_j-\psi)| \lesssim |\sin(\phi_{j+1}-\phi_j)|$, then again $\log\big(\max_{|z|=R}\|W_H(z)\|\big)\lesssim\ms k(R)$.
	\item If $\delta=2$ and $(\delta_l,\delta_\phi)\neq(1,1)$, then $\rho_H \leq \frac 12$.
	\end{Itemize}
\end{corollary}
\begin{proof}
	Our goal is to apply \Cref{B7}, and in order to do that we need to construct suitable functions $\ms c_l,\ms c_\phi$.
	Note that the assumption $\ms d_l \in L^1 ([1,\infty))$ implies $\delta_l \geq 1$.
\begin{Elist}
\item Without any a priori assumption, we can set
	\[
		\ms c_l(t) \DE \ms c_\phi (t) \DE \int_t^\infty \ms d_l(x) \DD x \gtrsim t \ms d_l(t).
	\]
	Then $\ms d_l,\ms d_\phi$ together with $\ms c_l,\ms c_\phi$ satisfy the general assumptions of \Cref{B7}. We have
	\[
		t \ms C(t) \lesssim \frac{1}{\ms d_l(t)} \lesssim \ms D(t)
	\]
	and $\gamma=\delta_l-1 \geq 0$. This choice of $\ms c_l,\ms c_\phi$ is sufficient to prove the asserted bounds in the
	following cases:
	\begin{Itemize}
	\item $\delta>2 \,\wedge \, \delta_l>1$. Since $\gamma>0$, \Cref{B7} gives the upper bound 
		$R^{\frac 12}\int\limits_{\ms k(R)}^\infty \ms D(s)^{-\frac 12}\DD s \asymp \ms k(R)$. 
	\item $\delta <2$. Then $\delta_\phi=\delta-\delta_l<1 \leq \delta_l$, hence $\frac{\ms d_\phi}{\ms d_l}$ has 
		positive index and is eventually nondecreasing. Since
		\[
			R^{\frac 12}\int\limits_1^{\ms h(R)}\ms D(s)^{-\frac 12}\DD s \asymp 
			R^{\frac 12} \ms h(R) \ms D(\ms h(R))^{-\frac 12} \lesssim R \int_{\ms h(R)}^\infty \ms d_l(s) \DD s,
		\] 
		\Cref{B7} gives the asserted upper bound.
	\item $\delta=2 \, \wedge \, \delta_l>1$. Then we are either in the second or fourth row of \Cref{B7}, but in both cases 
		the upper bound for the order is equal to $\frac 12$.
	\end{Itemize}
\item Assume $\delta>2$ and $\delta_\phi>1$. 
	This implies that $\sum_{j=1}^\infty |\sin (\phi_{j+1}-\phi_j)| <\infty$. We want to choose 
	$\psi \DE \lim_{j \to \infty} \phi_j$, so we need to prove that this limit actually exists
	(at least if the angles $\phi_j$ are all modified by adding integer multiples of $\pi$ which leaves $H$ unchanged).
	Start by choosing $n_0$ so large that $\sum_{j=n_0}^\infty |\sin (\phi_{j+1}-\phi_j)| \leq \frac{\pi}{8}$. 
	Adding to $\phi_j$ an integer multiple of $\pi$, we may assume $|\phi_j-\phi_{n_0}|
	\leq \frac{\pi}{2}$ for all $j \in \mathbb{N}$. Since $|x|\leq 2|\sin (x)|$ for $|x|\leq \frac{\pi}{2}$, and $|\sin
	(x+y)| \leq |\sin (x)|+|\sin (y)|$, we have for $j>k \geq n_0$ 
	\begin{align*}
		|\phi_j-\phi_k| &\leq |\phi_j-\phi_{n_0}|+|\phi_{n_0}-\phi_k| 
		\leq 2 \big(|\sin(\phi_j-\phi_{n_0})|+|\sin (\phi_{n_0}-\phi_k)| \big) 
		\\
		&\leq 4\sum_{n=n_0}^\infty |\sin (\phi_{n+1}-\phi_n)| \leq \frac{\pi}{2}.
	\end{align*}
	Therefore, 
	\[
		|\phi_j-\phi_k| \leq 2 |\sin (\phi_j-\phi_k)| \leq 2\sum_{n=k}^{j-1} |\sin (\phi_{n+1}-\phi_n)|
	\]
	and thus $(\phi_j)_{j=1}^\infty$ is a Cauchy sequence. Let $\psi$ be its limit. We set
	\[
		\ms c_l(t):= \int_t^{\infty} \ms d_l(x) \DD x, \qquad 
		\ms c_\phi (t):= \ms c_l(t) \cdot \Big(\int_t^\infty \ms d_\phi (x) \DD x \Big)^2
	\]
	and observe that $\ms d_l,\ms d_\phi$ together with $\ms c_l,\ms c_\phi$ satisfy the assumptions of \Cref{B7}. A
	calculation shows that $\gamma=\delta-2>0$. Hence $\Ind [t\ms C(t)|=\delta-1<\delta$ and in particular 
	$t\ms C(t) \ll\ms D(t)$.  Again \Cref{B7} provides the desired upper bound. 
\item Assume that we have $\psi \in \bb R$ such that $|\phi_j-\psi| \lesssim |\phi_{j+1}-\phi_j|$, and $\delta_l >1$.
	Set 
	\begin{align*}
		\ms c_l(t)&:= \int_t^{\infty} \ms d_l(x) \DD x \asymp t\ms d_l(t), 
		\\
		\ms c_\phi (t) &:= \int_t^{\infty} \ms d_l(x)\ms d_\phi (x)^2 \DD x \asymp t \ms d_l(t) \ms d_\phi (t)^2.
	\end{align*}
	Then $t \ms C(t) \asymp \ms D(t)$ and we are in the first row of \Cref{B7}. Now we note that the bound given in the
	theorem is $\asymp$ to $\ms k(R)$.
\end{Elist}
\end{proof}

\section{Additions and examples}
\subsection{Combining with a bound from below}

In the same way that upper bounds for lengths and angle differences lead to upper bounds for the growth of the monodromy matrix, 
lower bounds lead to lower bounds. We recall a result obtained in \cite[Corollary~2.5]{pruckner.woracek:sinqA}.

\begin{proposition}
\label{B35}
	Let $H$ be a Hamburger Hamiltonian with lengths $(l_j)_{j=1}^\infty$ and angles $(\phi_j)_{j=1}^\infty$, and let 
	$\ms f$ be regularly varying with positive index. If 
	\begin{equation}
	\label{B33}
		l_{j+1}l_j \sin^2 (\phi_{j+1} - \phi_j) \gtrsim \frac{1}{\ms f(j)}, \qquad j \in \bb N,
	\end{equation}
	then
	\begin{align*}
		\log\Big(\max_{|z|=R}\|W_H(z)\|\Big) \gtrsim \ms f^-(R^2)
	\end{align*}
	for sufficiently large $R$.
\end{proposition}

\noindent
Let us translate \Cref{B35} to the setting where we compare the lengths and the angle differences to regularly varying functions
$\ms d_l, \ms d_\phi$.

\begin{corollary}
\label{B41}
	Let $(l_j)_{j=1}^\infty$ be a summable sequence of positive numbers, and let $(\phi_j)_{j=1}^\infty$ be a sequence of 
	real numbers. Denote by $H$ the Hamburger Hamiltonian with these lengths and angles, 
	and let $W_H$ be its monodromy matrix. 
	Assume that $\ms d_l,\ms d_\phi$ are regularly varying and satisfy
	\[
		\forall j \in \bb N. \quad l_j \gtrsim \ms d_l(j) \qquad |\sin (\phi_{j+1}-\phi_j )| \gtrsim \ms d_\phi (j).
	\]
	Then
	\[
		\log\Big(\max_{|z|=R}\|W_H(z)\|\Big) \gtrsim \Big[\frac{1}{\ms d_l \ms d_\phi} \Big]^-(R)
	\]
	for sufficiently large $R$. 

	In particular, the order of $W_H$ is at least $\frac 1\delta$, where $\delta\DE-(\Ind\ms d_l+\Ind\ms d_\phi)$.
\end{corollary}
\begin{proof}
	Since $\ms d_l$ is regularly varying, we have $\ms d_l(t+1) \sim \ms d_l(t)$. Setting 
	$\ms D(t) \DE \frac{1}{(\ms d_l\ms d_\phi)(t)}$, we see that \cref{B33} is satisfied for $\ms f(t)\DE\ms D(t)^2$. 
	Since $\ms f^-(t)=\ms D^- (t^{\frac 12})$, we obtain 
	\[
		\log\Big(\max_{|z|=R}\|W_H(z)\|\Big) \gtrsim \ms f^-(R^2)=\ms D^-(R)=\Big[\frac{1}{\ms d_l \ms d_\phi} \Big]^-(R).
	\]
\end{proof}

\noindent 
If the lengths and angle differences are well-behaved and summable, the growth of $W_H$ can be determined up to $\asymp$. 
Note that no functions $\ms c_l,\ms c_\phi$ appear in the formulation of the following theorem.

\begin{theorem}
\label{B38}
	Let $(l_j)_{j=1}^\infty$ be a summable sequence of positive numbers, and let $(\phi_j)_{j=1}^\infty$ be a sequence of 
	real numbers. Denote by $H$ the Hamburger Hamiltonian with these lengths and angles.
	Consider regularly varying functions $\ms d_l,\ms d_\phi$ with $\ms d_l \asymp 1 \asymp \ms d_\phi$ locally. If
	\begin{Enumerate}
	\item $l_j \asymp \ms d_l(j)$ and $|\sin (\phi_{j+1}-\phi_j )| \asymp \ms d_\phi (j)$ for sufficiently large $j$,
	\item $\delta \DE - \big(\Ind \ms d_l+\Ind \ms d_\phi \big)>2$,
	\end{Enumerate}
	then
	\[
		\log\Big(\max_{|z|=R}\|W_H(z)\|\Big) \, \asymp \,  \Big[\frac{1}{\ms d_l \ms d_\phi} \Big]^-(R)
		\qquad\text{for sufficiently large }R,
	\]
	and 
	\[
		\rho_H=\frac{1}{\delta}.
	\]
\end{theorem}
\begin{proof}
	Due to \Cref{B41}, we only need to show that 
	\[
		\log\big(\max_{|z|=R}\|W_H(z)\|\big) \lesssim \big[\frac{1}{\ms d_l \ms d_\phi} \big]^-(R)
		.
	\]
	We notice that $\Ind \ms d_l \leq -1$ since $(l_j)_{j=1}^\infty$ is summable. Hence $\ms d_l$ is $\sim$ to an
	eventually monotone function, and thus
	\[
	\int_M^\infty \ms d_l(x) \DD x \lesssim \sum_{j=M}^\infty \ms d_l(j) \lesssim  \sum_{j=M}^\infty l_j <\infty.
	\]
	Since $\ms d_l \asymp 1$ locally, this shows that $\ms d_l \in L^1 ([1,\infty))$. 
	By \Cref{B96},
	\[
		\log\Big(\max_{|z|=R}\|W_H(z)\|\Big) \lesssim \ms k(R) \asymp \Big[\frac{1}{\ms d_l \ms d_\phi} \Big]^-(R).
	\]
\end{proof}

\subsection{Power-log-majorisations and exceptional cases}

By considering $\ms d_l,\ms d_\phi,\ms c_l,\ms c_\phi$ consisting of a power times a power of a logarithm, we can gain some 
insight into the exceptional cases of our results. For data of this form all functions occurring in our results can in principle
be computed explicitly. We do not aim at giving a complete picture, but rather give a couple of illustrative examples. It should
also be added that not all phenomena (counterexamples, exceptional cases, or similar) can be illustrated with functions of this
form; one would have to admit an additional double-logarithmic factor. 

The facts presented below are shown by straightforward, yet tedious and elaborate, computations. 
It is practical to use the lexicographic order on $\bb R^2$, and we denote it by $\preceq$. Explicitly, thus
\[
	(\alpha,\beta)\preceq(\alpha',\beta')\DI \alpha<\alpha'\vee\big(\alpha=\alpha'\wedge\beta\leq\beta'\big)
\]
and, as ususal, $\prec$ stands for ``$\preceq$ but not $=$''.

\begin{GenericDefinition}{Setting and Notation}
	Assume we are given parameters 
	\[
		(\delta_l,\alpha_l),(\delta_\phi,\alpha_\phi),(\gamma_l,\beta_l),(\gamma_\phi,\beta_\phi)
		\in[0,\infty)\times\bb R
		,
	\]
	and denote 
	\[
		\delta\DE\delta_l+\delta_\phi,\ \alpha\DE\alpha_l+\alpha_\phi, 
		\gamma\DE\frac 12(\gamma_l+\gamma_\phi),\ \beta\DE\frac 12(\beta_l+\beta_\phi)
		.
	\]
	Assume that these parameters satisfy 
	\begin{Itemize}
	\item ${\displaystyle 
		(\delta_l,\alpha_l),(\delta_\phi,\alpha_\phi),(\gamma_l,\beta_l),(\gamma_\phi,\beta_\phi)\succeq(0,0)
		}$,
	\item ${\displaystyle 
		\delta>0,(\gamma,\beta)\succ(0,0)
		}$,
	\item ${\displaystyle 
		(\gamma_l,\beta_l)\preceq(\gamma_\phi,\beta_\phi)
		}$.
	\end{Itemize}
	Let $\ms d_l,\ms d_\phi,\ms c_l,\ms c_\phi$ be continuous and nonincreasing functions, such that $\ms d_\phi\leq 1$, 
	$\ms c_\phi\leq\ms c_l$, and that (for sufficiently large $t$)
	\begin{align*}
		& \ms d_l(t)=t^{-\delta_l}(\log t)^{-\alpha_l},\quad 
		\ms d_\phi(t)=t^{-\delta_\phi}(\log t)^{-\alpha_\phi}
		,
		\\
		& \ms c_l(t)=t^{-\gamma_l}(\log t)^{-\beta_l},\quad 
		\ms c_\phi(t)=t^{-\gamma_\phi}(\log t)^{-\beta_\phi}
		.
	\end{align*}
\end{GenericDefinition}

\noindent
We compute some of the basic ingredients. Here, and throughout the following, all formulas are understood to hold for
sufficiently large $t$ or $R$.

\begin{lemma}
\label{B2}
	We have 
	\begin{Enumerate}
	\item ${\displaystyle
		\ms D(t)=\PL\delta\alpha t,\quad \ms C(t)=\PL\gamma\beta t
		}$,
		\\[2mm]
		${\displaystyle
		\frac{\ms d_\phi(t)}{\ms d_l(t)}=\PL{\delta_l-\delta_\phi}{\alpha_l-\alpha_\phi}t,\quad
		\frac{\ms c_l(t)}{\ms c_\phi(t)}=\PL{\gamma_\phi-\gamma_l}{\beta_\phi-\beta_l}t
		}$.
	\item ${\displaystyle
		\ms k(R)\approx \PL{\frac 1\delta}{-\frac\alpha\delta}R
		}$,
		\\[2mm]
		${\displaystyle
			\ms h(R)=\infty\ \Leftrightarrow\ (\delta_l,\alpha_l)\preceq(\delta_\phi,\alpha_\phi)
		}$,
		\\[1mm]
		${\displaystyle
		\ms h(R)\approx
		\begin{cases}
			\PL{\frac 1{\delta_l-\delta_\phi}}{-\frac{\alpha_l-\alpha_\phi}{\delta_l-\delta_\phi}}R
			\CAS \delta_l>\delta_\phi,
			\\[1mm]
			\exp\Big(R^{\frac 1{\alpha_l-\alpha_\phi}}\Big)
			\CAS \delta_l=\delta_\phi\wedge\alpha_l>\alpha_\phi.
		\end{cases}
		}$
	\end{Enumerate}
\end{lemma}
\begin{proof}
	The formulas stated in {\rm(i)} follow directly from the definitions. To see {\rm(ii)} we first use \Cref{B19}:
	the function $\ms k(R)$ is $\approx$ to an asymptotic inverse of $\ms D(t)$, and if $\delta_l>\delta_\phi$ then 
	$\ms h(R)$ is an asymptotic inverse of $\frac{\ms d_\phi(t)}{\ms d_l(t)}$. 

	If $(\delta_l,\alpha_l)\preceq(\delta_\phi,\alpha_\phi)$, then the quotient $\frac{\ms d_\phi}{\ms d_l}$ is bounded and
	hence $\ms h(R)=\infty$. If $\delta_l=\delta_\phi$ and $\alpha_l>\alpha_\phi$, we solve the equation 
	$\big(\frac{\ms d_\phi}{\ms d_l}\circ\ms h)(R)=R$. 
\end{proof}

\begin{lemma}
\label{B23}
	We have 
	\[
		\int_1^\infty\ms D(s)^{-\frac 12}\DD s<\infty\ \Leftrightarrow\ (\delta,\alpha)\succ(2,2)
	\]
	\[
		\int_t^\infty\ms D(s)^{-\frac 12}\DD s=
		\begin{cases}
			\frac 1{\frac\delta2-1}\cdot\PL{1-\frac\delta 2}{-\frac\alpha 2}t
			\CAS \delta>2
			\\
			\frac 1{\frac\alpha 2-1}\cdot\pL{1-\frac\alpha 2}t
			\CAS \delta=2\wedge\alpha>2
		\end{cases}
	\]
	\[
		\int_1^t\ms D(s)^{-\frac 12}\DD s\sim
		\begin{cases}
			\frac 1{1-\frac\delta2}\cdot\PL{1-\frac\delta 2}{-\frac\alpha 2}t
			\CAS \delta<2,
			\\
			\frac 1{1-\frac\alpha 2}\cdot\pL{1-\frac\alpha 2}t
			\CAS \delta=2\wedge\alpha<2,
			\\
			\log\log t
			\CAS \delta=\alpha=2.
		\end{cases}
	\]
\end{lemma}
\begin{proof}
	The equivalence stated in the first line is clear, and the stated formulas for $\delta>2$ and $\delta<2$ follow for
	example from Karamata's theorem (or explicit calculation). Assume that $\delta=2$. A primitive of 
	$\frac 1t(\log t)^{-\frac\alpha 2}$ is given by 
	\[
		\begin{cases}
			\frac 1{1-\frac\alpha 2}\cdot\pL{1-\frac\alpha 2}t
			\CAS \alpha\neq 2,
			\\
			\log\log t
			\CAS \alpha=2,
		\end{cases}
	\]
	and also in this case the assertion follows. 
\end{proof}

\noindent
In our first example we discuss the role of the term $L(t,R)$; this fits the context of \Cref{B59}. We show that there are
situations where $L(t,R)$ cannot be neglected, but also that the assumptions in \Cref{B59} are only sufficient for 
$B(R)\asymp\ms B(R)$. 

\begin{example}
\label{B24}
	Assume that 
	\begin{equation}
	\label{B27}
		(\delta_l,\alpha_l)\prec(\delta_\phi,\alpha_\phi),\quad \gamma=0,\quad(\delta,\alpha)\succ(2,2)
		.
	\end{equation}
	The first assumption is there to rule out applicability of \Cref{B59}\,(i), the second to rule out applicability of 
	\Cref{B59}\,(ii),(iii), and the third to reduce computational effort (the facts we want to illustrate occur already under
	this additional assumption). 

	For parameters subject to \cref{B27} it holds that 
	\begin{equation}
	\label{B28}
		\ms B(R)\asymp
		\begin{cases}
			\PL{\frac 12}{1-\frac\alpha2}R
			\CAS \delta=2,
			\\
			\PL{\frac 1\delta}{-\frac\alpha\delta}R
			\CAS \delta>2.
		\end{cases}
	\end{equation}
	\begin{equation}
	\label{B29}
		B(R)
		\begin{cases}
			\asymp\ms B(R)
			\CAS \delta_l=\delta_\phi\vee\beta>1\vee(\beta=\delta-1>1\wedge\alpha<0),
			\\
			\asymp R^{\frac 1{1+\beta}}\gg\ms B(R)
			\CASO
			.
		\end{cases}
	\end{equation}
\end{example}
\begin{proof}
	We start with showing \cref{B28}. If $\delta>2$, the stated formula follows from \Cref{B9}\,\fbox{B}. 

	Assume that $\delta=2$, then $\alpha>2$. By \Cref{B51} we have 
	\[
		\ms B(R)\leq\sup_{t\geq 1}\ms g(t,R)\lesssim R^{\frac 12}\int_{\ms k(R)}^\infty\ms D(s)^{-\frac 12}\DD s
		\asymp\PL{\frac 12}{1-\frac\alpha2}R
		.
	\]
	Set $\tau(R)\DE\exp\big(R^{\frac 1{2\beta}}\big)$. Since 
	\[
		\frac R{\ms C(\tau(R))}=R^{\frac 12}\gg\PL{\frac 12}{1-\frac\alpha2}R\gtrsim\ms g(\tau(R),R)
		,
	\]
	and 
	\[
		\int_{\tau(R)}^\infty\ms D(s)^{-\frac 12}\DD s\asymp R^{\frac{1-\frac\alpha2}{2\beta}}
		,
	\]
	and $1-\frac\alpha2<0$, we have 
	\begin{align*}
		\ms B(R)\gtrsim &\, \ms g(\tau(R),R)\gtrsim R^{\frac 12}\int_{\ms k(R)}^{\tau(R)}\ms D(s)^{-\frac 12}\DD s
		\\
		\gtrsim &\, \PL{\frac 12}{1-\frac\alpha2}R-R^{\frac 12}R^{\frac{1-\frac\alpha2}{2\beta}}
		\asymp\PL{\frac 12}{1-\frac\alpha2}R
		.
	\end{align*}
	For the proof of \cref{B29} we have to include the term $L(t,R)$ into the discussion. We have 
	\[
		B(R)\asymp\inf_{t\geq 1}\max\Big\{\ms g(t,R),\frac R{\ms C(t))},L(t,R)\Big\}
		,
	\]
	and hence 
	\begin{multline*}
		\max\Big\{\ms B(R),\inf_{t\geq 1}\Big\{\frac R{\ms C(t))},L(t,R)\Big\}\Big\}\lesssim B(R)
		\\
		\lesssim \max\Big\{\sup_{t\geq 1}\ms g(t,R),\inf_{t\geq 1}\Big\{\frac R{\ms C(t))},L(t,R)\Big\}\Big\}
		.
	\end{multline*}
	By what we already showed $\ms B(R)\asymp\sup_{t\geq 1}\ms g(t,R)$, and it remains to evaluate the written infimum. 
	By \Cref{B68} and the fact that 
	\[
		\log^+\frac{\ms c_l(t)}{\ms c_\phi(t)}\lesssim\log\log t
		,
	\]
	we obtain that 
	\[
		L(t,R)\asymp 1+\log^+R+
		\begin{cases}
			\log t\CAS \delta_l<\delta_\phi,
			\\
			\log\log t\CAS \delta_l=\delta_\phi.
		\end{cases}
	\]
	Note that therefore $t\mapsto L(t,R)$ is $\asymp$ to a nondecreasing continuous function. 

	Consider first the case that $\delta_l<\delta_\phi$. Then we set $\tau(R)\DE\exp\big(R^{\frac 1{1+\beta}}\big)$, 
	and obtain 
	\[
		\frac R{\ms C(\tau(R))}=R^{1-\frac\beta{1+\beta}}=R^{\frac 1{1+\beta}}\asymp L(\tau(R),R)
		.
	\]
	and hence 
	\[
		\min_{t\geq 1}\max\Big\{\frac R{\ms C(t)},L(t,R)\Big\}\asymp R^{\frac 1{1+\beta}}
		.
	\]
	Assume now that $\delta_l=\delta_\phi$. Then we set 
	$\tau(R)\DE\exp\Big(\big(\frac R{\log R}\big)^{\frac 1\beta}\Big)$, and obtain 
	\[
		\frac R{\ms C(\tau(R))}=\log R\asymp L(\tau(R),R)
		,
	\]
	which leads to 
	\[
		\min_{t\geq 1}\max\Big\{\frac R{\ms C(t)},L(t,R)\Big\}\asymp\log R
		.
	\]
	Putting together, the relation \cref{B29} follows. 
\end{proof}

\noindent
In our second example we discuss the bounds from \Cref{B9}\,\fbox{C}. We show that they are sharp but need not necessarily be
attained. 

\begin{example}
\label{B36}
	Assume that 
	\begin{equation}
	\label{B50}
		(\delta_l,\alpha_l)\preceq(1+\gamma,\beta)\preceq(\delta,\alpha),\quad\delta=2,\alpha\leq 2,\quad\gamma>0
		.
	\end{equation}
	These assumptions ensure that we are in the situation of \Cref{B9}\,\fbox{C} and that in \cref{B75} we do not
	automatically have equality. 

	For parameters subject to \cref{B50} it holds that 
	\begin{equation}
	\label{B52}
		\frac R{\ms C(\ms f_0^-(R))}\asymp\PL{\frac 12}{-\frac\alpha2}R\asymp\ms k(R)
		,
	\end{equation}
	\begin{equation}
	\label{B60}
		\frac R{\ms C(\ms f_1^-(R))}\asymp
		\begin{cases}
			\PL{\frac 12}{1-\frac\alpha2}R\CAS \alpha<2
			,
			\\
			R^{\frac 12}\log\log R\CAS \alpha=2
			.
		\end{cases}
	\end{equation}
	\begin{equation}
	\label{B63}
		\ms B(R)\asymp
		\begin{cases}
			\PL{\frac 12}{1-\frac\alpha2}R\CAS \gamma<1
			,
			\\
			\PL{\frac 12}{-\frac\alpha2}R\log\log R\CAS \gamma=1\wedge\alpha>\beta
			,
			\\
			\PL{\frac 12}{-\frac\alpha2}R\CAS \gamma=1\wedge\alpha=\beta
			.
		\end{cases}
	\end{equation}
	We see that 
	\begin{Itemize}
	\item $\ms B(R)\asymp\frac R{\ms C(\ms f_1^-(R))}$ if $\gamma<1\wedge\alpha<2$,
	\item $\frac R{\ms C(\ms f_0^-(R))}\ll\ms B(R)\asymp\frac R{\ms C(\ms f_1^-(R))}$ if $\gamma<1\wedge\alpha=2$ 
		or $\gamma=1\wedge\alpha>\beta$,
	\item $\ms B(R)\asymp\frac R{\ms C(\ms f_0^-(R))}$ if $\gamma=1\wedge\alpha=\beta$.
	\end{Itemize}
\end{example}
\begin{proof}
	Plugging the definitions and using \Cref{B19} yields
	\[
		\ms f_0(t)=\PL{2\gamma}{2\beta-\alpha}t,\quad
		\ms f_0^-(R)\asymp\bigg(\frac R{(\log R)^{2\beta-\alpha}}\bigg)^{\frac 1{2\gamma}},
	\]
	\[
		\ms f_1(t)\sim
		\begin{cases}
			\frac 1{2^{-\alpha}}\PL{2\gamma}{2-\alpha+2\beta}t\CAS \alpha<2,
			\\
			\PL{2\gamma}{2\beta}t(\log\log t)^2\CAS \alpha=2,
		\end{cases}
	\]
	\[
		\ms f_1^-(R)\approx
		\begin{cases}
			\bigg(\frac R{(\log R)^{2-\alpha+2\beta}}\bigg)^{\frac 1{2\gamma}}\CAS \alpha<2,
			\\
			\bigg(\frac R{(\log R)^{2\beta}(\log\log R)^2}\bigg)^{\frac 1{2\gamma}}\CAS \alpha=2.
		\end{cases}
	\]
	From this \cref{B52} and \cref{B60} follow immediately. 

	Consider the case that $\gamma<1$. Set 
	\[
		\tau(R)\DE\bigg(\frac R{(\log R)^{2-\alpha+2\beta}}\bigg)^{\frac 1{2\gamma}}
		.
	\]
	Then 
	\[
		\ms k(R)\ll\tau(R)\ll\ms f_0^-(R)\lesssim\ms h(R)
		,
	\]
	and hence 
	\[
		\ms g(\tau(R),R)\asymp\ms k(R)+R^{\frac 12}\int_{\ms k(R)}^{\tau(R)}\ms D(s)^{-\frac 12}\DD s
		.
	\]
	Using \Cref{B23}, we obtain 
	\begin{align*}
		\int_{\ms k(R)}^{\tau(R)}\ms D(s)^{-\frac 12}\DD s\approx &\,
		\begin{cases}
			(\log\tau(R))^{1-\frac\alpha2}-(\log\ms k(R))^{1-\frac\alpha2}\CAS \alpha<2
			\\
			\log\log\tau(R)-\log\log\ms k(R)\CAS \alpha=2
		\end{cases}
		\\
		\asymp &\, 
		\begin{cases}
			(\log R)^{1-\frac\alpha2}\CAS \alpha<2
			\\
			1 \CAS \alpha=2
		\end{cases}
		\qquad\asymp(\log R)^{1-\frac\alpha2}
		.
	\end{align*}
	Hence, 
	\[
		\ms g(\tau(R),R)\asymp\PL{\frac 12}{1-\frac\alpha2}R
		.
	\]
	Plugging the definitions shows that also 
	\[
		\frac R{\ms C(\tau(R))}\asymp\PL{\frac 12}{1-\frac\alpha2}R
		,
	\]
	and it follows that $\ms B(R)\asymp\PL{\frac 12}{1-\frac\alpha2}R$. 

	Next, consider the case that $\gamma=1$ and $\alpha>\beta$. Set 
	\[
		\tau(R)\DE\frac{\PL{\frac 12}{\frac\alpha2-\beta}R}{\log\log R}
		.
	\]
	Again
	\[
		\ms k(R)\ll\tau(R)\ll\ms f_0^-(R)\ll\ms h(R)
		,
	\]
	and therefore
	\[
		\ms g(\tau(R),R)\asymp\ms k(R)+R^{\frac 12}\int_{\ms k(R)}^{\tau(R)}\ms D(s)^{-\frac 12}\DD s
		.
	\]
	The 1\textsuperscript{st} mean value theorem provides us with $\xi(R)\in[\ms k(R),\tau(R)]$ such that 
	\begin{align*}
		\int_{\ms k(R)}^{\tau(R)}\ms D(s)^{-\frac 12}\DD s
		= &\, \int_{\ms k(R)}^{\tau(R)}s^{-1}(\log s)^{-\frac\alpha2}\DD s
		\\
		= &\, (\log\xi(R))^{-\frac\alpha2}\int_{\ms k(R)}^{\tau(R)}s^{-1}\DD s
		=(\log\xi(R))^{-\frac\alpha2}\log\frac{\tau(R)}{\ms k(R)}
		.
	\end{align*}
	We have $\log\tau(R)\sim\log\ms k(R)\approx\log R$, and hence $\log\xi(R)\asymp\log R$. Moreover, 
	\[
		\frac{\tau(R)}{\ms k(R)}=\frac{(\log R)^{\alpha-\beta}}{\log\log R}
		.
	\]
	It follows that
	\[
		\int_{\ms k(R)}^{\tau(R)}\ms D(s)^{-\frac 12}\DD s\approx(\log R)^{-\frac\alpha2}\log\log R
		,
	\]
	and in turn 
	\[
		\ms g(\tau(R),R)\asymp\PL{\frac 12}{-\frac\alpha2}R\log\log R
		.
	\]
	Again simply plugging the definitions shows that also 
	\[
		\frac R{\ms C(\tau(R))}\asymp\PL{\frac 12}{-\frac\alpha2}R\log\log R
		,
	\]
	and thus $\ms B(R)\asymp\PL{\frac 12}{-\frac\alpha2}R\log\log R$.

	It remains to settle the case that $\gamma=1$ and $\alpha=\beta$, but this is easily done. Simply plug the definitions
	to obtain 
	\[
		\frac R{\ms C(\ms k(R))}\asymp\ms k(R)\asymp\ms g(\ms k(R),R)
		,
	\]
	and therefore $\ms B(R)\asymp\ms k(R)\asymp\PL{\frac 12}{-\frac\alpha2}R$. 
\end{proof}

\noindent
In our third example we discuss the exceptional case ``$(\delta,\gamma)=(2,0)$'' in the third row of the table in \Cref{B7}. 
We show that in some cases (interpreting $\ms f_1^-$ appropriately) the written bound still holds and is even attained by
$B(R)$, while in others $B(R)$ is strictly larger due to domination of $L(t,R)$. 

\begin{example}
\label{B11}
	Assume that 
	\begin{equation}
	\label{B65}
		(\delta_l,\alpha_l)\preceq(1,1+\beta),\quad (\delta,\gamma)=(2,0),\quad \alpha\leq 2
		.
	\end{equation}
	These assumptions ensure that the exceptional case from the third row of the table in \Cref{B7} (and also of 
	\Cref{B9}\,\fbox{C}) takes place. Note that \cref{B65} implies that $\delta_l\leq\delta_\phi$. Moreover, 
	$\beta>0$ since $(\gamma,\beta)\succ(0,0)$. 

	For parameters subject to \cref{B65} it holds that 
	\[
		\ms B(R)\asymp
		\begin{cases}
			R^{\frac{2-\alpha+\beta}{2-\alpha+2\beta}}\CAS \alpha<2,
			\\
			R^{\frac 12}\log R \CAS \alpha=2,
		\end{cases}
	\]
	and 
	\[
		B(R)\ 
		\begin{cases}
			\asymp\ms B(R)
			\CAS \delta_l=\delta_\phi\vee(\delta_l<\delta_\phi\wedge\alpha\leq 1+\beta),
			\\
			\asymp R^{\frac 1{1+\beta}}\gg\ms B(R)
			\CAS \delta_l<\delta_\phi\wedge\alpha>1+\beta.
		\end{cases}
	\]
	Observe, moreover, that the bound for order in the third row of the table in \Cref{B7} has no continuous extension to
	$(2,0)$; its directional limits vary from $\frac 12$ to $1$. The above formula shows that the actual order of the bound
	$B(R)$ has nothing to do with $(\delta,\gamma)$ (being equal to $(2,0)$). Yet, it is sometimes given by the same
	formula, only with the ``logarithmic exponents'' $\alpha,\beta$ instead of $\delta,\gamma$. Also note that the 
	exponent $\frac 1{1+\beta}$ also occurred in \Cref{B24}.
\end{example}
\begin{proof}
	We have 
	\[
		\ms f_1(t)\sim
		\begin{cases}
			\frac 1{2-\alpha}(\log t)^{2-\alpha+2\beta}\CAS \alpha<2,
			\\
			(\log t)^{2\beta}(\log\log t)^2\CAS \alpha=2.
		\end{cases}
	\]
	The function $\ms f_1\circ\exp$ is regularly varying with positive index, and has approximate inverse
	\[
		(\ms f_1\circ\exp)^-(R)=
		\begin{cases}
			R^{\frac 1{2-\alpha+2\beta}}\CAS \alpha<2,
			\\
			\Big(\frac{R^{\frac 12}}{\log R}\Big)^{\frac 1\beta}\CAS \alpha=2.
		\end{cases}
	\]
	We define 
	\[
		\ms f_1^-(R)\DE(\exp\circ(\ms f_1\circ\exp)^-)(R)=
		\begin{cases}
			\exp\Big(R^{\frac 1{2-\alpha+2\beta}}\Big)\CAS \alpha<2,
			\\
			\exp\Big(\Big(\frac{R^{\frac 12}}{\log R}\Big)^{\frac 1\beta}\Big)\CAS \alpha=2.
		\end{cases}
		,
	\]
	so that
	\[
		(\ms f_1\circ\ms f_1^-)(R)\asymp R,\quad \log\big[(\ms f_1^-\circ\ms f_1)(R)\big]\asymp\log R
		.
	\]
	For the sake of consistency we use also here the notation of a function $\tau(R)$, and set $\tau(R)\DE\ms f_1^-(R)$. 
	Then 
	\[
		\ms k(R)\ll\tau(R)\ll\ms h(R)
		,
	\]
	where the second relation is seen as follows: If $\ms h(R)=\infty$, there is nothing to prove. Otherwise, we must have
	$\delta_l=\delta_\phi$ and $\alpha_l>\alpha_\phi$. Since $\delta=2$, thus $\delta_l=\delta_\phi=1$, and it follows that
	$\alpha_l\leq 1+\beta$. This implies that 
	\[
		\alpha_l-\alpha_\phi\leq 2-\alpha+2\beta
		.
	\]
	We conclude that 
	\[
		\ms g(\tau(R),R)\asymp\ms k(R)+R^{\frac 12}\int_{\ms k(R)}^{\tau(R)}\ms D(s)^{-\frac 12}\DD s
		.
	\]
	Using \Cref{B23}, we obtain
	\begin{align*}
		\int_{\ms k(R)}^{\tau(R)}\ms D(s)^{-\frac 12}\DD s\approx &\,
		\begin{cases}
			(\log\tau(R))^{1-\frac\alpha2}-(\log\ms k(R))^{1-\frac\alpha2}\CAS \alpha<2,
			\\
			\log\log\tau(R)-\log\log\ms k(R)\CAS \alpha=2,
		\end{cases}
		\\
		\asymp &\, 
		\begin{cases}
			R^{\frac{1-\frac\alpha2}{2-\alpha+2\beta}}\CAS \alpha<2,
			\\
			\log R \CAS \alpha=2.
		\end{cases}
	\end{align*}
	Hence, 
	\[
		\ms g(\tau(R),R)\asymp
		\begin{cases}
			R^{\frac{2-\alpha+\beta}{2-\alpha+2\beta}}\CAS \alpha<2,
			\\
			R^{\frac 12}\log R \CAS \alpha=2.
		\end{cases}
	\]
	Plugging the definitions yields 
	\[
		\frac R{\ms C(\tau(R))}\asymp
		\begin{cases}
			R^{\frac{2-\alpha+\beta}{2-\alpha+2\beta}}\CAS \alpha<2,
			\\
			R^{\frac 12}\log R \CAS \alpha=2,
		\end{cases}
	\]
	and the formula asserted for $\ms B(R)$ follows. 

	In order to show the asserted formulas for $B(R)$, we have to include the term $L(t,R)$ into the discussion. 
	We distinguish several cases.

	Assume that $\delta_l=\delta_\phi$ and $\alpha_l\geq\alpha_\phi$. The \Cref{B59}\,(i) shows that 
	$B(R)\asymp\ms B(R)$. 
	Assume that $\delta_l=\delta_\phi$ and $\alpha_l<\alpha_\phi$. We have 
	\[
		\log^+\frac{\ms c_l(\tau(R))}{\ms c_\phi(\tau(R))}=
		\log^+\Big[(\log\tau(R))^{\beta_\phi-\beta_l}\Big]\lesssim\log R
		,
	\]
	\[
		\log^+\frac{\ms d_l(\tau(R))}{\ms d_\phi(\tau(R))}=
		\log^+\Big[(\log\tau(R))^{\alpha_\phi-\alpha_l}\Big]\approx\log R
		.
	\]
	Now \Cref{B68} implies that $L(\tau(R),R)\asymp\log R$, and we see that $B(R)\asymp\ms B(R)$.
	Assume that $\delta_l<\delta_\phi$ and $\alpha\leq 1+\beta$. Then 
	\[
		\log^+\frac{\ms d_l(\tau(R))}{\ms d_\phi(\tau(R))}\asymp\log\tau(R)\asymp
		\begin{cases}
			R^{\frac 1{2-\alpha+2\beta}}\CAS \alpha<2,
			\\
			\Big(\frac{R^{\frac 12}}{\log R}\Big)^{\frac 1\beta}\CAS \alpha=2,
		\end{cases}
	\]
	and hence, again referring to \Cref{B68}, also 
	\[
		L(\tau(R),R)\asymp
		\begin{cases}
			R^{\frac 1{2-\alpha+2\beta}}\CAS \alpha<2,
			\\
			\Big(\frac{R^{\frac 12}}{\log R}\Big)^{\frac 1\beta}\CAS \alpha=2.
		\end{cases}
	\]
	This shows that $L(\tau(R),R)\lesssim\ms B(R)$, and therefore again $B(R)\asymp\ms B(R)$. 

	Assume now that $\delta_l<\delta_\phi$ and $\alpha>1+\beta$. This case is different, and we also use a different 
	function $\tau(R)$, namely, we now set 
	\[
		\tau(R)\DE\exp\Big(R^{\frac 1{1+\beta}}\Big)
		.
	\]
	We have $\log\frac{\ms d_l(t)}{\ms d_\phi(t)}\asymp\log t$, and hence obtain 
	\[
		L(\tau(R),R)\asymp R^{\frac 1{1+\beta}}
		.
	\]
	Further,
	\[
		\frac R{\ms C(\tau(R))}=R^{\frac 1{1+\beta}},\quad 
		\ms g(\tau(R),R)\asymp 
		\begin{cases}
			R^{\frac 12+\frac{1-\frac\alpha2}{1+\beta}} \CAS \alpha<2,
			\\
			R^{\frac 12}\log R \CAS \alpha=2.
		\end{cases}
	\]
	Since $\alpha>1+\beta$, we have 
	\[
		\frac 12+\frac{1-\frac\alpha2}{1+\beta}<\frac 1{1+\beta}
		,
	\]
	and see that 
	\[
		B(R)\asymp\min_{t\geq 1}\max\Big\{\ms g(t,R),\frac R{\ms C(t)},L(t,R)\Big\}\asymp R^{\frac 1{1+\beta}}
		.
	\]
\end{proof}

\subsection{Two corollaries given in terms of Jacobi parameters}

We present two applications of \Cref{B38} in which we return to the regime of power moment problems. The first is a supplement
to a result from \cite{pruckner:blubb}, and in the second we give examples where the Nevanlinna matrix has prescribed growth.

At this point we need the concrete formulae relating Jacobi parameters with Hamiltonian parameters. They read as 
\begin{equation}
\label{B71}
	\begin{aligned}
		\frac 1{b_n}= &\, \sin(\phi_{n+1}-\phi_n)\sqrt{l_{n+1}l_n},
		\\
		a_n= &\, -\frac 1{l_n}\big[\cot(\phi_{n+1}-\phi_n)+\cot(\phi_n-\phi_{n-1})\big],
	\end{aligned}
\end{equation}
where the angles are chosen such that $\phi_{n+1}-\phi_n\in[0,\pi)$, cf.\ \cite{kac:1999}. 
Given the Jacobi parameters, it is in general hard to solve the equations \cref{B71} for the Hamiltonian parameters. Under the
assumption that $b_n,a_n$ have a certain power-like asymptotic, the approximate size of $l_n$ and $|\sin(\phi_{n+1}-\phi_n)|$
can be determined. 

\begin{Corollary}
\label{B79}
	Let $a_n\in\bb R$ and $b_n>0$ be sequences which have asymptotics 
	\begin{align*}
		b_n= &\, n^\sigma\bigg(\frac{|y_0|}2+\frac{x_1}n+\frac{x_2}{n^2}+\BigO\Big(\frac 1{n^{2+\epsilon}}\Big)\bigg),
		\\
		a_n= &\, n^\sigma\bigg(y_0+\frac{y_1}n+\frac{y_2}{n^2}+\BigO\Big(\frac 1{n^{2+\epsilon}}\Big)\bigg),
	\end{align*}
	where 
	\[
		\sigma>2,\ y_0\neq 0,\quad x_1,x_2,y_1,y_2\in\bb R,\quad \epsilon>0
		.
	\]
	Assume that the moment problem with these Jacobi parameters is indeterminate, and let $W(z)$ be its Nevanlinna matrix. 
	Then 
	\[
		\log\Big(\max_{|z|=R}\|W(z)\|\Big)\asymp R^{\frac 1\sigma}
		.
	\]
\end{Corollary}

\noindent
Before we come to the proof, let us put this statement in the right context. It deals with the critical situation that 
off-diagonal and diagonal of the Jacobi matrix are comparable with ratio $2$. This setting was considered in 
\cite[Theorem~2]{pruckner:blubb}. In that theorem occurrence of limit circle case was characterised in terms of the data of the
expansions, and it was shown that $W(z)$ is of order $\frac 1\sigma$ with positive type. The significance of \Cref{B79} is that
now we know that $W(z)$ is also of finite type. 

\begin{proof}
	In the proof of \cite[Theorem~2]{pruckner:blubb} it was shown that 
	\[
		l_n\asymp\lambda(n)^2,\quad |\sin(\phi_{n+1}-\phi_n)|\asymp\frac 1{n^\sigma\lambda(n)^2}
		,
	\]
	where $\lambda$ is a function of the form $\lambda(t)=n^\tau$ or $\lambda(t)=n^\tau\log t$ with some 
	$\tau\in[-\frac\sigma 2,-\frac 12)$. The assumptions of \Cref{B38} are thus satisfied with 
	$\ms d_l(t)\DE\lambda(t)^2$ and $\ms d_\phi(t)\DE\frac 1{n^\sigma\lambda(n)^2}$. Applying this theorem yields 
	\[
		\log\Big(\max_{|z|=R}\|W(z)\|\Big)\asymp\Big[\frac 1{\ms d_l\ms d_\phi}\Big]^-(R)\asymp R^{\frac 1\sigma}
		.
	\]
\end{proof}

\noindent
We come to our second corollary, where we produce a variety of examples with prescribed growth of the Nevanlinna matrix (slower
than the threshold $R^{\frac 12}$). Thereby, the speed of growth is always determined by the off-diagonal, and the diagonal can
be as large or as small as we please. 

\begin{Corollary}
\label{B83}
	Let $\ms g$ be regularly varying with $\Ind\ms g\in(0,\frac 12)$. 
	\begin{Enumerate}
	\item Let $\omega\in[-2,2]$. Then there exist Jacobi parameters $b_n$ and $a_n$, such that 
		\[
			b_n\sim\ms g^-(n),\ \frac{a_n}{b_n}\to\omega,\quad
			\log\Big(\max_{|z|=R}\|W(z)\|\Big)\asymp\ms g(R)
			.
		\]
	\item Let $\omega_n\neq 0$, $\omega_n\to 0$, be such that $\lim_{n\to\infty}\frac{\omega_{n-1}}{\omega_n}$ exists in 
		$(-1,\infty)$. Then there exist Jacobi parameters $b_n$ and $a_n$, such that 
		\[
			b_n\sim\ms g^-(n),\ \frac{a_n}{b_n}\sim\omega_n,\quad
			\log\Big(\max_{|z|=R}\|W(z)\|\Big)\asymp\ms g(R)
			.
		\]
	\end{Enumerate}
\end{Corollary}

\noindent
The condition in (i) that $|\omega|\leq 2$ is no restriction, since otherwise we could not have limit circle case by 
Wouk's theorem.

\begin{proof}
	We specify lengths and angles, which is done differently in different cases. 
	\begin{Itemize}
	\item Assume that $\omega\in(-2,2)$. Let $\psi\in(0,\pi)$ be such that $\cos\psi=-\frac\omega2$, and set 
		\[
			l_n\DE\frac 1{\sin\psi\cdot\ms g^-(n)},\qquad 
			\phi_{n+1}\DE n\psi,\ n\in\bb N,
		\]
		so that $\phi_{n+1}-\phi_n=\psi$. 
	\item Assume that $\omega=-2$. Set 
		\[
			l_n\DE\frac n{\ms g^-(n)},\qquad 
			\phi_{n+1}\DE\sum_{k=1}^n\frac 1k,\ n\in\bb N,
		\]
		so that $\phi_{n+1}-\phi_n=\frac 1n$.
	\item Assume that $\omega=2$. Set 
		\[
			l_n\DE\frac n{\ms g^-(n)},\qquad 
			\phi_{n+1}\DE n\pi-\sum_{k=1}^n\frac 1k,\ n\in\bb N,
		\]
		so that $\phi_{n+1}-\phi_n=\pi-\frac 1n$.
	\item Assume that $\omega_n\neq 0$, $\omega_n\to 0$, and 
		$\gamma\DE\lim_{n\to\infty}\frac{\omega_{n-1}}{\omega_n}$ exists in $(-1,\infty)$. Set 
		\[
			l_n\DE\frac 1{\ms g^-(n)},\qquad 
			\phi_{n+1}\DE n\frac\pi2+\frac 1{2(1+\gamma)}\sum_{\substack{k=1\\|\omega_n|<\pi}}^n\omega_n,\ n\in\bb N,
		\]
		so that $\phi_{n+1}-\phi_n=\frac\pi2+\frac{\omega_n}{2(1+\gamma)}$ for all sufficiently large $n$.
	\end{Itemize}
	Let $b_n$ and $a_n$ be the Jacobi parameters given by \cref{B71}. Then, in all cases, $b_n\sim\ms g^-(n)$. 
	Multiplying the two equations from \cref{B71}, shows that 
	\[
		\frac{a_n}{b_n}=-\sqrt{\frac{l_{n+1}}{l_n}}\cdot\bigg(\cos(\phi_{n+1}-\phi_n)+\cos(\phi_n-\phi_{n-1})
		\frac{\sin(\phi_{n+1}-\phi_n)}{\sin(\phi_n-\phi_{n-1})}\bigg)
		,
	\]
	and this implies the asserted property of $\frac{a_n}{b_n}$. Finally, we apply \Cref{B38} with the obvious choices for 
	$\ms d_l,\ms d_\phi$ to obtain that 
	\[
		\log\Big(\max_{|z|=R}\|W(z)\|\Big)\asymp\ms g(R)
		.
	\]
\end{proof}

%
%
\addtocontents{toc}{\protect\contentsline{section}{\protect\numberline{}\textcolor{Sepia}{\kern-15pt Appendix}}{}{}}
\appendix
\makeatletter
\DeclareRobustCommand{\@seccntformat}[1]{%
  \def\temp@@a{#1}%
  \def\temp@@b{section}%
  \ifx\temp@@a\temp@@b
  Appendix\ \csname the#1\endcsname.\quad%
  \else
  \csname the#1\endcsname\quad%
  \fi
}
\makeatother
\renewcommand{\thelemma}{\Alph{section}.\arabic{lemma}}
%
%

%
%
%
\section{Regularly varying functions}

In complex analysis the growth of an analytic function is compared with functions of the form $\exp(\ms a(r))$. The most 
classical comparison functions are powers $\ms a(r)=r^\rho$, and this leads to the notions of order and type. 
A refined comparison scale was introduced already at a very early stage by E.Lindel\"of \cite{lindeloef:1905} who considered 
comparison functions behaving for $r\to\infty$ like
\[
	r^\alpha\cdot\bigl(\log r\bigr)^{\beta_1}\cdot\bigl(\log\log r\bigr)^{\beta_2}
	\cdot\ldots\cdot
	\bigl(\underbrace{\log\cdots\log}_{\text{\footnotesize$m$\textsuperscript{th} iterate}}r\bigr)^{\beta_m},
\]
where $\alpha>0$ and $\beta_1,\ldots,\beta_m\in\bb R$.
Functions which are nowadays commonly used as comparison functions are regularly varying functions in Karamata sense, 
cf.\ \cite[Chapter~7]{bingham.goldie.teugels:1989} (for other levels of generality see also \cite{levin:1980,rubel:1996}). 
Lindel\"of's comparison functions are examples of functions of that kind. 

Let us now recall Karamata's definition of regular variation.

\begin{Definition}
\label{B4}
	A function $\ms a\DF[1,\infty)\to(0,\infty)$ is called \emph{regularly varying} at $\infty$ with 
	\emph{index} $\alpha\in\bb R$, if it is measurable and 
	\[
		\forall \lambda\in(0,\infty)\DP \lim_{r\to\infty}\frac{\ms a(\lambda r)}{\ms a(r)}=\lambda^\alpha
		.
	\]
	We write $\Ind\ms a$ for the index of regular variation of function $\ms a$. 
	A regularly varying function with index $0$ is also called \emph{slowly varying}.
\end{Definition}

\noindent
Regularly varying function $\ms a$ are used to quantify growth for $r\to\infty$, and hence the values of $\ms a(r)$ for small
$r$ are irrelevant. This allows to change $\ms a$ on any finite interval without changing the essence of results, and this
freedom can often be used to assume $\ms a$ has some additional practical properties. 

We cite a number of fundamental theorems on regularly varying functions. 
Proofs can be found, e.g., in \cite{bingham.goldie.teugels:1989} or \cite{seneta:1976}.
We start with an \emph{representation theorem}.

\begin{theorem}[Representation theorem]
\label{B39}
	Let $\alpha \in \bb R$. A function $\ms a\DF[1,\infty)\to(0,\infty)$ is regularly varying with index $\alpha$ 
	if and only if it has a representation of the form 
	\[
		\ms a(r)=r^\alpha \cdot c(r) \exp \bigg(\int_1^r \epsilon (u) \frac{\DD u}{u} \bigg), \qquad r \in [1,\infty),
	\] 
	where $c,\epsilon$ are measurable, $\lim_{r \to \infty} c(r)=c \in (0,\infty)$, and $\lim_{r \to \infty} \epsilon(r)=0$. 

	If $\ms a$ is slowly varying (i.e., $\alpha=0$) and eventually nondecreasing (nonincreasing), then $\epsilon$ may be 
	taken eventually nonnegative (nonpositive). 
\end{theorem}

\noindent
It is a legitimate intuition that regularly varying functions fill in the scale of powers, and that a regularly varying 
function with index $\alpha$ behaves roughly like the power $r^\alpha$. The following results, which we will use frequently,
express this intuition very clearly. The first is a variant of the \emph{Potter bounds}, and the second is the classical 
\emph{Karamata Theorem} about asymptotic integration. 

\begin{theorem}[Potter bounds (variant)]
\label{B10}
	Let $\ms a$ be regularly varying with index $\alpha\in\bb R$. 
	\begin{Enumerate}
	\item ${\displaystyle
		\forall \epsilon>0\DP r^{\alpha-\epsilon}\ll\ms a(r)\ll r^{\alpha+\epsilon}
		}$
		,		
	\item $\lim_{r \to \infty} \frac{\log \ms a(r)}{\log r}=\alpha$,
	\item For all $\epsilon>0$ the quotients $\frac{\ms a(r)}{r^{\alpha-\epsilon}}$ and 
		$\frac{r^{\alpha+\epsilon}}{\ms a(r)}$ are $\sim$ to an eventually increasing function.
	\end{Enumerate}
\end{theorem}

\begin{theorem}[Karamata's Theorem]
\label{B12}
	Let $\ms a$ be regularly varying with index $\alpha\in\bb R$. 
	\begin{Enumerate}
	\item Assume that $\alpha\geq -1$. Then the function $x\mapsto\int_1^x \ms a(t)\DD t$ is regularly
		varying with index $\alpha+1$, and 
		\[
			\lim_{x\to\infty}\bigg(
			\raisebox{3pt}{$x\ms a(x)$}\Big/\,\raisebox{-2pt}{$\int\limits_1^x \ms a(t)\DD t$}
			\bigg)=\alpha+1
			.
		\]
	\item Assume that $\alpha\leq -1$ and $\int_1^\infty \ms a(t)\DD t<\infty$. Then the function 
		$x\mapsto\int_x^\infty \ms a(t)\DD t$ is regularly varying with index $\alpha+1$, and 
		\[
			\lim_{x\to\infty}\bigg(
			\raisebox{3pt}{$x\ms a(x)$}\Big/\,\raisebox{-2pt}{$\int\limits_x^\infty \ms a(t)\DD t$}
			\bigg)=-(\alpha+1)
			.
		\]
	\end{Enumerate}
\end{theorem}

\noindent
A regularly varying function $\ms a$ with positive index is -- at least asymptotically -- invertible. In fact, if
\[
	\ms a^-(x) \DE \sup \big\{t\in[1,\infty)\DS \ms a(t)<x\big\},
\]
we have the following result, cf.\ \cite[Theorem~1.5.12]{bingham.goldie.teugels:1989}.

\begin{theorem}
\label{B16}
	Let $\ms a$ be regularly varying with index $\alpha>0$. Then $\ms a^-$ is regularly varying with index 
	$\frac{1}{\alpha}$, and 
	\begin{equation}
	\label{B40}
		(\ms a\circ\ms a^-)(x)\sim(\ms a^-\circ\ms a)(x)\sim x
		.
	\end{equation}
\end{theorem}

\noindent
Any regularly varying function $\ms a^-$ with the property \cref{B40} is called an \emph{asymptotic inverse} of $\ms a$, and 
asymptotic inverses are determined uniquely up to $\sim$.
We recall a useful formula for computing asymptotic inverses for functions of a particular form.

\begin{Remark}
\label{B19}
	Assume that $\rho>0$ and that $\ms f$ is regularly varying. Set $\ms g\DE\ms f\circ\log$ (for sufficiently large $t$). 
	Then 
	\[
		\ms a(t)\DE t^\rho\ms g(t),\quad \ms a^-(t)\DE \rho^{\frac{\Ind\ms f}\rho}\cdot
		\Big(\frac t{\ms g(t)}\Big)^{\frac 1\rho}
	\]
	are asymptotic inverses of each other. 
\end{Remark}

\noindent
Another practical observation is the following.

\begin{Remark}
\label{B69}
	Let $\ms f$ be regularly varying with $\Ind\ms f>0$, and assume we have a function $\ms g$ with $\ms g\approx\ms f$. 
	Then $\ms g$ is regularly varying with $\Ind\ms g=\Ind\ms f$ and $\ms g^-\approx\ms f^-$. 
\end{Remark}

\noindent
By the Potter bounds every regularly varying function $\ms a$ is bounded and bounded away from zero on every interval
$[r_1,r_2]$ sufficiently far to the right. Sometimes it is needed for technical reasons to assume this property for all compact 
intervals in the domain of $\ms a$. Of course, this is no loss in generality; remember that modification on a finite
interval does not change the essence of the function $\ms a$. 

\begin{lemma}
\label{B42}
	Let $\ms a \DF [1,\infty) \to (0,\infty)$ be slowly varying and assume that $\ms a\lesssim 1$ locally. 
	For $R \geq \ms a(1)$ set 
	\[
		\ms b(R) \DE \sup \Big\{t\in[1,\infty)\DS \sup_{1\leq s\leq t}\frac{\ms a(s)}{R}\leq 1\Big\} \, \in [1,\infty].
	\]
	Then $\ms b$ grows faster than any power, i.e., $R^\rho \leq \ms b(R)$ for every $\rho>0$ and $R$ sufficiently large.
\end{lemma}
\begin{proof}
	Let $\rho>0$ and set $\epsilon \DE \frac 1\rho$. By \Cref{B10}, there exists $M>0$ such that $\ms a(r) \leq r^\epsilon$ 
	for all $r \geq M$. If $R \geq \sup_{r \in [1,M]} \ms a(r)$, this means that
	\[
		\Big\{t\in[1,\infty)\DS \sup_{1\leq s\leq t}\frac{s^\epsilon}{R} \leq 1\Big\} \subseteq 
		\Big\{t\in[1,\infty)\DS \sup_{1\leq s\leq t}\frac{\ms a(s)}{R}\leq 1\Big\}.
	\]
	The assertion follows if we take suprema of both sets.
\end{proof}

\begin{lemma}
\label{B31}
	Let $\ms a$ be regularly varying with index $\alpha>0$ and assume that $\ms a \asymp 1$ locally. Then there exists a
	continuously differentiable and regularly varying function $\ms s$, where $\ms s'(t)>0$ for $t \in [1,\infty)$,
	\begin{align*}
		\ms a(t) &\asymp \ms s(t), \qquad t \in [1,\infty),
		\\
		\ms a(t) &\sim s(t), \qquad t \to \infty.
	\end{align*}
\end{lemma}
\begin{proof}
	We use the smooth variation theorem \cite[Theorem~1.8.2]{bingham.goldie.teugels:1989}. This gives a function $\ms s$
	that is continuously differentiable, whose derivative is positive for all sufficiently large $t$, and such that $\ms a
	\sim \ms s$ for $t \to \infty$. W.l.o.g. we assume $\ms s'(t)>0$ for all $t \in [1,\infty)$. 

	Choose $t_0 \geq 1$ such that $\frac 12 \ms a(t) \leq \ms s(t) \leq \frac 32 \ms a(t)$ for all $t \geq t_0$. 
	By our assumption, $\ms a \asymp 1$ on $[1,t_0]$. Since $\ms s$ is continuous, $\ms s \asymp 1 \asymp \ms a$ on 
	$[1,t_0]$. Summing up, we have $\ms s \asymp\ms a$ on $[1,\infty)$.
\end{proof}

\noindent
Another lemma in a similar direction is the following. The proof is immediate from the representation theorem and we do not go
into details.

\begin{lemma}
\label{B54}
	Let $\ms a$ be regularly varying, and assume that $\ms a$ is $\approx$ to some nonincreasing function. Then there exist
	\[
		\epsilon[\ms a],\kappa[\ms a]\DF[1,\infty)\to(0,\infty)
	\]
	such that 
	\begin{Itemize}
	\item $\epsilon[\ms a]$ is locally integrable, $\lim_{t\to\infty}\epsilon[\ms a](t)=0$, and 
		$\epsilon[\ms a]\leq 0$ if $\Ind\ms a=0$,
	\item $\kappa[\ms a]$ is eventually constant,
	\item The function 
		\[
			\ms S[\ms a](t)\DE\kappa[\ms a](t)\cdot t^{\Ind\ms a}\cdot
			\exp\bigg(\int_1^r\epsilon[\ms a](u)\frac{\DD u}u\bigg)
		\]
		is nonincreasing, continuous, and $\approx\ms a$. 
	\end{Itemize}
\end{lemma}

\noindent
We speak of any function $\ms S[\ms a]$ as in the lemma as a \emph{nonincreasing smoothening} of $\ms a$. 

Of course, $\kappa[\ms a]$ and $\epsilon[\ms a]$ are far from unique. We mention two particular instances of the freedom of
choice in $\ms S[\ms a]$. 

\begin{Remark}
\label{B55}
	\phantom{}
	\begin{Enumerate}
	\item If $\ms S[\ms a]$ is some nonincreasing smoothening of $\ms a$ and $\alpha>0$, then also 
		$\alpha\cdot\ms S[\ms a]$ is a nonincreasing smoothening of $\ms a$. This corresponds to multiplying 
		$\kappa[\ms a]$ by $\alpha$. 
	\item If $\ms S[\ms a]$ is some nonincreasing smoothening of $\ms a$ and $\alpha>\lim_{t\to\infty}\ms S[\ms a](t)$, 
		then also $\min\{\ms S[\ms a](t),\alpha\}$ is a nonincreasing smoothening of $\ms a$. This is seen by modifying 
		$\kappa[\ms a]$ on a finite interval. 
	\end{Enumerate}
\end{Remark}

\noindent
One more property of this construction is as follows.
Assume we have two functions $\ms a_1,\ms a_2$ that are both regularly varying and $\approx$ to some nonincreasing function. 
If $\ms S[\ms a_1]$ and $\ms S[\ms a_2]$ are nonincreasing smoothenings of $\ms a_1$ and $\ms a_2$, respectively, then 
$\ms S[\ms a_1]\cdot\ms S[\ms a_2]$ is a nonincreasing smoothening of $\ms a_1\cdot\ms a_2$. This corresponds to taking 
\[
	\kappa[\ms a_1\ms a_2]\DE\kappa[\ms a_1]\cdot\kappa[\ms a_2],\quad 
	\epsilon[\ms a_1\ms a_2]\DE\epsilon[\ms a_1]+\epsilon[\ms a_2]
	.
\]


{\footnotesize
\begin{flushleft}
	R.~Pruckner\\
	Fachgruppe Patent- und Lizenzmanagement\\
	Vienna University of Technology\\
	Resselgasse~3/058\\
	1040~Wien\\
	AUSTRIA\\
	email: \texttt{raphael.pruckner@tuwien.ac.at}\\[5mm]
\end{flushleft}
\begin{flushleft}
	J.~Reiffenstein\\
	Department of Mathematics\\
	University of Vienna\\
	Oskar-Morgenstern-Platz~1\\
	1090~Wien\\
	AUSTRIA\\
	email: \texttt{jakob.reiffenstein@univie.ac.at}\\[5mm]
\end{flushleft}
\begin{flushleft}
	H.~Woracek\\
	Institute for Analysis and Scientific Computing\\
	Vienna University of Technology\\
	Wiedner Hauptstra{\ss}e\ 8--10/101\\
	1040~Wien\\
	AUSTRIA\\
	email: \texttt{harald.woracek@tuwien.ac.at}\\[5mm]
\end{flushleft}
}

\end{document}